\documentclass[12pt]{report}
\usepackage{nmtthes2000}
\usepackage[intlimits]{amsmath}
\usepackage{amssymb}
\usepackage{amsthm}
\usepackage{pst-all}

\newcommand{\ds}{\displaystyle}
\newcommand{\tld}{\widetilde}

\newcommand{\eps}{\epsilon}

\renewcommand{\rho}{\varrho}
\renewcommand{\epsilon}{\varepsilon}
\renewcommand{\phi}{\varphi}
\renewcommand{\kappa}{\varkappa}
\renewcommand{\emptyset}{\varnothing}

\newcommand{\Cinf}{C^\infty}
\newcommand{\Cinfcomp}{C^\infty_c}

\newcommand{\Ltwo}{\ensuremath{L^2}}

\newcommand{\N}{\ensuremath{\mathbb{N}}}
\newcommand{\Z}{\ensuremath{\mathbb{Z}}}

\newcommand{\R}{\ensuremath{\mathbb{R}}}
\newcommand{\CC}{\ensuremath{\mathbb{C}}}
\newcommand{\PP}{\ensuremath{\mathbb{P}}}

\newcommand{\mf}[1]{\mathfrak{#1}} 

\newcommand{\abs}[1]{\left|#1\right|}
\newcommand{\bigO}[1]{\mathrm{O}\left(#1\right)}

\newcommand{\norm}[1]{\left\|#1\right\|}
\newcommand{\set}[1]{\left\{#1\right\}}

\newcommand{\e}[1]{\exp\left(#1\right)}
\newcommand{\ip}[1]{\left(#1\right)}
\newcommand{\la}{\langle}
\newcommand{\ra}{\rangle}
\newcommand{\mi}[1]{\mathbf{#1}} 
\newcommand{\pd}{\partial}
\newcommand{\pdf}[3][]{\dfrac{\pd^{#1} {#2}}{\pd {#3}^{#1}}}
\newcommand{\sig}[1]{\sigma\left(#1\right)}
\newcommand{\tensorbundle}[2]{\mathcal{T}^{#1}_{#2}}
\newcommand{\vect}[1]{\mathbf{#1}}
\newcommand{\clearfrac}[2]{\genfrac{}{}{0pt}{}{#1}{#2}}

\DeclareMathOperator{\Alt}{Alt}
\DeclareMathOperator{\Det}{Det}

\DeclareMathOperator{\Ker}{Ker}
\DeclareMathOperator{\RE}{Re}
\DeclareMathOperator{\Sym}{Sym}
\DeclareMathOperator{\Trace}{Tr}

\DeclareMathOperator{\dist}{dist}
\DeclareMathOperator{\dvol}{dvol}

\DeclareMathOperator{\len}{len}
\DeclareMathOperator{\ord}{ord}

\DeclareMathOperator{\sign}{sign}

\DeclareMathOperator{\tr}{tr}

\DeclareMathOperator{\TwoPointTensorBundle}{\mathcal{W}}
\newcommand{\tptb}[4]{\sideset{^{#1}_{#2}}{^{#3}_{#4}}\TwoPointTensorBundle}

\newtheorem{corollary}{Corollary}[chapter]
\newtheorem*{definition}{Definition}
\newtheorem{lemma}{Lemma}[chapter]
\newtheorem*{remark_env}{Remark} \newenvironment{remark}{\begin{remark_env}\rm}{\end{remark_env}}
\newtheorem{theorem}{Theorem}[chapter]

\newcommand{\keyterm}[1]{\textbf{#1}}
\def\sideremark#1{\ifvmode\leavevmode\fi\vadjust{\vbox to0pt{\vss
\hbox to 0pt{\hskip\hsize\hskip1em
\vbox{\hsize2cm\tiny\raggedright\pretolerance10000
\noindent #1\hfill}\hss}\vbox to8pt{\vfil}\vss}}}

\author{Andrey Novoseltsev}

\title{SPECTRAL GEOMETRY OF RIEMANNIAN SUBMANIFOLDS} 

\degree{Master of Science in Mathematics}

\graduationdate{July, 2005}

\address{poselok Neytrino dom 1 kvartira 5\\ Elbrusskiy rayon\\ Kabardino-Balkariya\\ 361609, Russian Federation}

\supervisor{Ivan G. Avramidi}

\committeesize{3}

\begin{document}

	%
	%

\titlepage

\begin{abstract}
In this thesis we study the geometry of the fixed point set $\Sigma$ of a smooth mapping $\Phi:M\to M$ on an $n$-dimensional Riemannian manifold $(M,g)$ by computing the asymptotic expansion as $t\to 0^+$ of the trace of the deformed heat kernel $U\big(t;x,\Phi(x)\big)$ of the Laplace operator on $M$. We assume that the fixed point set is a union of connected components $\Sigma_i$, $i=1,\ldots,s$, each of which is a smooth compact submanifold of $M$ with dimension $m_i$, $0\leqslant m_i \leqslant n-1$.

There exists the asymptotic expansion as $t\to 0^+$
\begin{gather*}
\int_M \dvol(x)\, U\big(t;x,\Phi(x)\big) \sim
\sum_{i=1}^s  \sum_{k=0}^\infty t^{k-m_i/2} \int_{\Sigma_i} \dvol(y)\, a^{(i)}_k(y) \,,
\end{gather*}
where $a^{(i)}_k(y)$ are scalar invariants on $\Sigma_i$ depending on covariant derivatives of the curvature of the metric $g$ and symmetrized covariant derivatives of the differential $d\Phi$ of the mapping $\Phi$, evaluated on $\Sigma_i$.

Due to the localization principle, it is possible to compute the expansion for each component separately. We develop a generalized Laplace method for computing the coefficients $a^{(i)}_k(y)$ in this expansion and compute the coefficients $a_0$, $a_1$, and $a_2$ explicitly in the following cases:
\begin{enumerate}
\item $M$ is a flat two-dimensional manifold, $\Sigma$ is a zero- or one-dimensional component of the fixed point set.
\item $M$ is an $n$-dimensional curved manifold, $\Sigma$ is a zero-dimensional component of the fixed point set.
\end{enumerate}

We also develop algorithms for computing further coefficients and expressing them in a symmetrization-free form.
\end{abstract}

\begin{acknowledgments}
I would like to express my gratitude to my scientific advisor Prof. Ivan~G. Avramidi for arousing my interest in differential geometry and his area of work. For endless hours he spent with me during my years at New Mexico Tech exploring possible research directions and helping to improve my background where it lacked. He was always ready to answer my questions and it was just impossible to wish for a better mentor!

Thanks to my scientific advisor Prof. Igor~B. Simonenko at Rostov State University, Russia. His inexhaustible enthusiasm inspired me a lot during all my six years there, and our prolonged discussions helped me to develop a mathematical style of thinking.

Thanks to Don Clewett and especially Rachael Defibaugh-Ch\'avez for reading drafts of this thesis and making a lot of useful comments and corrections.
\end{acknowledgments}

\tableofcontents



\chapter{Introduction}

\section{Differential Geometry}

Here we provide a brief review of differential geometry concepts that we will use further. For more detailed treatment of the subject see, for example,~\cite{Frankel,Rosenberg}.

\begin{definition}
A topological Hausdorff space $M$ with a countable base is called an $n$-dimensional \keyterm{manifold} (without boundary), $n\in\N$, if, for any point $p\in M$, there exists an open neighborhood $U\subset M$ of the point $p$, open subset $V\subset\R^n$, and homeomorphism $\phi:U\to V$. The pair $(U,\phi)$ is called a \keyterm{coordinate chart} or \keyterm{local coordinates}.
\end{definition}

We always assume that all manifolds we deal with are smooth. That is, given an $n$-dimensional manifold $M$ and two coordinate patches $(U,\phi)$ and $(V,\psi)$, such that $U\cap V\neq\emptyset$, the corresponding \keyterm{transition map}
\begin{gather*}
\psi\circ\phi^{-1}:\phi(U\cap V)\to \R^n
\end{gather*}
is smooth (infinitely differentiable).

\begin{definition} Let $M$ be a manifold and $f:M\to\R$ be a function on it. The function $f$ is \keyterm{continuously differentiable} on $M$ if, for any point $p\in M$ and any coordinate patch $(U,\phi)$, containing the point $p$, the function $f\circ\phi^{-1}: \phi(U)\to\R$ is continuously differentiable. If $f$ is continuously differentiable on $M$, we will write $f\in C^1(M)$.
\end{definition}

In the same manner one can give definitions of $k$-times differentiable functions $C^k(M)$ and infinitely differentiable functions $\Cinf(M)$, as well as differentiable mappings $f:M\to \R^n$. We will denote the space of smooth functions with compact support  by $\Cinfcomp(M)$.

\begin{definition}
Let $M$ be an $n$-dimensional manifold. A \keyterm{vector} $\vect{X}$ at a point $p\in M$ is a linear differential operator of the first order $\vect{X}:C^1(M)\to\R$. In local coordinates it can be written as
\begin{gather*}
\vect{X}(f)=\sum_{i=1}^n X^i \left[\pdf{f}{x^i}\right](p) \,,
\end{gather*}
where $X^i$ are components of $\vect{X}$. The vector $\vect{X}$ is also called the \keyterm{tangent} or \keyterm{contravariant vector}.
The set of all tangent vectors to the manifold $M$ at the point $p$ is called the \keyterm{tangent space} and is denoted by $T_p M$.
The union of all tangent spaces to the manifold $M$ at all points is called the \keyterm{tangent bundle} and is denoted by $TM$.
A \keyterm{vector field} (or a \keyterm{section of the tangent bundle}) $\vect{v}$ on the manifold $M$ is a mapping $\vect{v}: M\to TM$, such that, for any point $p\in M$, there holds $\vect{v}(p)\in T_p M$.
\end{definition}

Below we follow the standard Einstein convention of summation over repeated indices from one to the dimension of the manifold. Where it is necessary, we will specify limits of the summation explicitly. In local coordinates a vector field $\vect{v}$ may be written as
\begin{gather*}
\vect{v}(p)=v^i(p)\,\pd_i \,,
\end{gather*}
with $\pd_i$ being an $i$-th vector of the coordinate basis in the tangent space $T_p M$, defined by
\begin{gather*}
\pd_i(f)=\pdf{f}{x^i} \,,
\end{gather*}
where $x^i$ are local coordinates.
We will say that the vector field $\vect{v}$ is smooth and will write $\vect{v}\in\Cinf(TM)$, if all its components $v^i(x)$ are smooth functions in all coordinate patches.

\begin{remark} One can show that, although we use local coordinates to introduce various geometrical objects, they have invariant geometrical meaning which does not depend on the choice of a coordinate patch. All our results will be formulated in the invariant form.
\end{remark}

\begin{definition}
Let $M$ be a manifold. The \keyterm{cotangent space} $T^*_p M$ at a point $p\in M$ is the dual space to $T_p M$, i.e. the space of all linear functionals on $T_p M$. Its elements are called \keyterm{cotangent vectors}, \keyterm{covariant vectors} or \keyterm{covectors}.
The union of all cotangent spaces at all points is called the \keyterm{cotangent bundle} and is denoted by $T^*M$.
A \keyterm{covector field} (also called a \keyterm{$1$-form} or a \keyterm{section of the cotangent bundle}) $\vect{v}$ on the manifold $M$ is a mapping $\vect{v}: M\to T^* M$, such that, for any point $p\in M$, there holds $\vect{v}(p)\in T^*_p M$.
\end{definition}

The basis in a cotangent space, dual to the coordinate basis $\pd_i$, is denoted by $dx^i$ and is defined by the condition $dx^i(\pd_j)=\delta^i_j$, where $\delta^i_j$ is the Kronecker symbol. As well as for vector fields, we define the set of all smooth covector fields $\Cinf(T^* M)$.

\begin{definition}
Let $\Phi:M\to M$ be a smooth mapping on a manifold $M$. Let $p\in M$ and $p'=\Phi(p)$. Let $x^i$ be local coordinates in the point $p$, $x^{i'}$ be local coordinates in the point $p'$, and $\Phi^{i'}$ be the coordinate functions of the mapping $\Phi$. The \keyterm{differential} of the mapping $\Phi$ at the point $p$ is a linear mapping
\begin{gather*}
d\Phi:T_p M\to T_{p'} M
\end{gather*}
with components
\begin{gather*}
(d\Phi)^{i'}{}_j=\pdf{\Phi^{i'}}{x^j} \,.
\end{gather*}
\end{definition}

The \keyterm{tensor product} of two functionals $f:X\to\R$ and $g:Y\to\R$ is a functional $f\otimes g:X\times Y\to\R$ defined by $(f\otimes g)(x,y)=f(x)g(y)$ for any $x\in X$ and $y\in Y$. 

\begin{definition}
A \keyterm{tensor $T$ of the type $(s,q)$} at a point $p$ of a manifold $M$ is a linear functional
\begin{gather*}
T: \underbrace{T_p M\times\ldots\times T_p M}_q \times \underbrace{T^*_p M\times\ldots\times T^*_p M}_s \to\R \,.
\end{gather*}
In the coordinate basis it is given by
\begin{gather*}
T=T^{i_1 \ldots i_s}_{j_1 \ldots j_q}\, dx^{j_1}\otimes\ldots\otimes dx^{j_q} \otimes \pd_{i_1}\otimes\ldots\otimes\pd_{i_s} \,.
\end{gather*}
The set of all tensors of the type $(s,q)$ at all points of the manifold $M$ is called the \keyterm{$(s,q)$-tensor bundle} and is denoted by $\tensorbundle{s}{q} M$.
A smooth \keyterm{$(s,q)$-tensor field} is a smooth mapping assigning to each point $p\in M$ a tensor of the type $(s,q)$.
The set of all smooth $(s,q)$-tensor fields is denoted by $\Cinf(\tensorbundle{s}{q}M)$.
\end{definition}
When $s$ or $q$ in the above definition equals zero, we will not write it, that is, $\tensorbundle{p}{0}=\tensorbundle{p}{}$ and $\tensorbundle{0}{q}=\tensorbundle{}{q}$.

For any $m\in\N$ let $S_m$ be the set of all permutations of the numbers $1,\ldots,m$.
A transposition is a permutation, which exchanges only two numbers, leaving all others fixed.
Every permutation can be represented as a composition of transpositions.
Such a representation is not unique, but the number of transpositions in it is an invariant $\!\!\mod 2$.
The signature of a permutation $\phi\in S_m$, denoted by $\sign(\phi)$, is $1$, if it is represented by a composition of even number of transpositions, and $-1$ otherwise.
For any permutation $\phi\in S_m$ and any tensor $T\in\tensorbundle{}{m} M$ we define the tensor $\phi T\in\tensorbundle{}{m} M$ by
\begin{gather*}
(\phi T)_{i_1\ldots i_m}=T_{i_{\phi(1)}\dots i_{\phi(m)}} \,.
\end{gather*}

We will say that $T\in\tensorbundle{}{m} M$ is a \keyterm{symmetric tensor}, if $\phi T=T$ for any $\phi\in S_m$.
Let $V_m M \subset \tensorbundle{}{m} M$ be the space of symmetric $(0,m)$-tensors.
We define the \keyterm{operator of symmetrization} $\Sym: \tensorbundle{}{m} M\to V_m M$ by
\begin{gather*}
\Sym(T)=\frac{1}{m!}\sum_{\phi\in S_{m}} \phi T \,.
\end{gather*}
We will also use parentheses to indicate symmetrization:
\begin{gather*}
(\Sym(T))_{i_1\ldots i_m}=T_{(i_1\dots i_m)}.
\end{gather*}
If it is necessary to exclude some indices from symmetrization, we will delimit them by $|\ldots|$, that is,
\begin{gather*}
T_{(i_1\dots i_k | j_1\ldots j_l | i_{k+1} \ldots i_m)}=
\frac{1}{m!}\sum_{\phi\in S_{m}}
T_{(i_{\phi(1)}\dots i_{\phi(k)} | j_1\ldots j_l | i_{\phi(k+1)} \ldots i_{\phi(m)})} \,.
\end{gather*}

We will say that $T\in\tensorbundle{}{m} M$ is an \keyterm{antisymmetric tensor}, if $\phi T=\sign(\phi) T$ for any $\phi\in S_m$. 
Let $\Lambda^m M \subset \tensorbundle{}{m} M$ be the space of antisymmetric $(0,m)$-tensors (also called \keyterm{$m$-forms}).
It follows from the definition, that, for any $n$-dimensional manifold $M$, only spaces $\Lambda^m M$ with $0\leqslant m\leqslant n$ are non-trivial, that is, they contain non-zero tensors.
We define the \keyterm{operator of antisymmetrization} (alternation) ${\Alt: \tensorbundle{}{m} M\to \Lambda^m M}$ by
\begin{gather*}
\Alt(T)=\frac{1}{m!}\sum_{\phi\in S_{m}} \sign(\phi) \phi T \,.
\end{gather*}
We will also use square brackets to indicate alternation:
\begin{gather*}
(\Alt(T))_{i_1\ldots i_m}=T_{[i_1\dots i_m]}.
\end{gather*}
The \keyterm{exterior} or \keyterm{wedge product} $\wedge:\Lambda^p M \times\Lambda^q M \to\Lambda^{p+q} M$ is defined by
\begin{gather*}
\alpha\wedge\beta=\frac{(p+q)!}{p!\, q!}\Alt(\alpha\otimes\beta) \,.
\end{gather*}

When it is not confusing, we will write $\tensorbundle{p}{q}$, $V_n$, and $\Lambda^n$ instead of $\tensorbundle{p}{q} M$, $V_n M$, and $\Lambda^n M$ respectively.

We say that a manifold $M$ is equipped with a \keyterm{Riemannian metric} $g$ if $g$ is a smooth symmetric tensor field of the type $(0,2)$ and it defines a positive definite non-degenerate inner product in each tangent space $T_p M$. We will call the pair $(M,g)$ a \keyterm{Riemannian manifold}. On a Riemannian manifold we can use the metric for raising and lowering tensor indices:
\begin{gather*}
a_i=g_{ij}a^j,\qquad a^j=g^{ji}a_i,
\end{gather*}
where $g^{ji}$ are components of the inverse matrix to $g_{ij}$, i.e. $g^{ji}g_{ik}=\delta^j_k$. The determinant of the matrix $g_{ij}$ we will denote by $\det g$ or by $\abs{g}$.

For integration over a manifold $M$ we need to define a volume form.

\begin{definition}
An $n$-dimensional manifold $M$ is \keyterm{orientable}, if there exists a constant sign $n$-form $\alpha\in\Lambda^n$.
\end{definition}

In case of an orientable Riemannian $n$-dimensional manifold $(M,g)$ there is a special volume form generated by the metric. It is called the \keyterm{Riemannian volume form} and is given by
\begin{gather*}
\dvol=\eps(\pd_1,\ldots,\pd_n) \abs{g}^{1/2}\, dx^1\wedge\ldots\wedge dx^n,
\end{gather*}
where $\eps(\pd_1,\ldots,\pd_n)$ is equal to $1$ for positively oriented basis $\pd_i$ and $-1$ otherwise. The volume form changes sign whenever a basis changes the orientation; such a geometrical object is called a \keyterm{pseudoform}. We will always work with orientable manifolds and will use only positively oriented bases without specifying it explicitly. In this case $\eps(\pd_1,\ldots,\pd_n)=1$ and we will omit it henceforth.

When computing integrals over a manifold, it is usually necessary to work in several coordinate charts. For using formulas, written in local coordinates, a partition of unity over a manifold is used for splitting an integral into parts, such that each part can be covered by a single chart.

\begin{definition}
Let $(M,g)$ be a Riemannian manifold. The \keyterm{Christoffel symbols} $\Gamma^i_{jk}$ are defined by
\begin{gather*}
\Gamma^i{}_{jk}=\frac{1}{2} g^{il}\left(\pdf{g_{kl}}{x^j}+\pdf{g_{jl}}{x^k}-\pdf{g_{jk}}{x^l}\right).
\end{gather*}
The \keyterm{covariant derivative} (also called the \keyterm{connection, compatible with the metric}, or \keyterm{Levi-Civita connection}) is an operator $\nabla:\Cinf(\tensorbundle{s}{q})\to\Cinf(\tensorbundle{s}{q+1})$, defined by
\begin{gather*}
(\nabla T)^{i_1 \ldots i_s}_{j_1 \ldots j_q k}=
\pd_k T^{i_1 \ldots i_s}_{j_1 \ldots j_q}
+\sum_{l=1}^s \Gamma^{i_l}{}_{mk} T^{i_1 \ldots i_{l-1} m i_{l+1} \ldots i_s}_{j_1 \ldots j_q}
-\sum_{l=1}^q \Gamma^{m}{}_{j_{l} k} T^{i_1 \ldots i_s}_{j_1 \ldots j_{l-1} m j_{l+1} \ldots j_q} \,.
\end{gather*}
For a given point $p\in M$ and a vector $\vect{X}\in T_p M$, the covariant derivative in the direction of the vector $\vect{X}$ is
\begin{gather*}
(\nabla_\vect{X} T)^{i_1 \ldots i_s}_{j_1 \ldots j_q}
=X^k (\nabla T)^{i_1 \ldots i_s}_{j_1 \ldots j_q k}  \,.
\end{gather*}
\end{definition}
We will denote covariant derivatives by indices after a semicolon,
\begin{gather*}
T^{i_1 \ldots i_s}_{j_1 \ldots j_q ; k_1\ldots k_r}=(\nabla^r T)^{i_1 \ldots i_s}_{j_1 \ldots j_q k_1\ldots k_r}.
\end{gather*}

\begin{definition}
The components of the \keyterm{Riemann curvature} tensor are defined by
\begin{gather*}
R^i{}_{jkl}
	=\pdf{\Gamma^i{}_{jl}}{x^k}-\pdf{\Gamma^i{}_{jk}}{x^l}
	+\Gamma^i{}_{ml}\Gamma^m{}_{jk}-\Gamma^i{}_{mk}\Gamma^m{}_{jl} \,.
\end{gather*}
\end{definition}
The Riemann curvature tensor has the following properties:
\begin{gather}
\label{eq:CurvSym}
R_{ijkl}=-R_{jikl}=-R_{ijlk}=R_{klij}
\,,\\
\label{eq:CurvBianci}
R^i{}_{jkl}+R^i{}_{klj}+R^i{}_{ljk}=0
\,,\\
\label{eq:CurvBianciDer}
R^{ij}{}_{kl;m}+R^{ij}{}_{lm;k}R^{ij}{}_{mk;l}=0
\,.
\end{gather}
The last two equations are called \keyterm{Bianci identities}. For the commutator of covariant derivatives there holds
\begin{gather}
\label{eq:CovDerCom}
2\, T^{i_1 \ldots i_s}_{j_1 \ldots j_q ; [kp]}=
\sum_{l=1}^s R^{i_l}{}_{mpk} T^{i_1 \ldots i_{l-1} m i_{l+1} \ldots i_s}_{j_1 \ldots j_q}
-\sum_{l=1}^q R^{m}{}_{j_{l} pk} T^{i_1 \ldots i_s}_{j_1 \ldots j_{l-1} m j_{l+1} \ldots j_q} \,.
\end{gather}
Contraction of two indices of the Riemann curvature tensor gives the \keyterm{Ricci curvature} tensor
\begin{gather*}
R^i{}_{jik}=R_{jk} \,,
\end{gather*}
and contraction of all indices gives the \keyterm{scalar curvature}
\begin{gather*}
R^{ij}{}_{ij}=R^k{}_k=R \,.
\end{gather*}

\begin{definition}
A \keyterm{geodesic} on a Riemannian manifold $(M,g)$ is a critical curve $\gamma:[0,1]\to M$ of the \keyterm{length functional} $\ds\ l(\gamma)=\int_0^1 dt\, \ip{\dot\gamma(t),\dot\gamma(t)}.$
\end{definition}
One can show that geodesics satisfy the equation
\begin{gather*}
\nabla_{\dot\gamma(t)}\dot\gamma(t)=0 \,,
\end{gather*}
which in local coordinates takes the form
\begin{gather} \label{eq:Geodesics}
\frac{d^2 \gamma^i}{dt^2} + \Gamma^i{}_{jk} \frac{d \gamma^j}{dt}\frac{d \gamma^k}{dt}=0\,.
\end{gather}
Since we consider smooth manifolds, the theory of ordinary differential equations guarantees that, for any point $p\in M$ and any vector $\vect{v}\in T_p M$, there exists an $\eps>0$,  such that equation~\eqref{eq:Geodesics} has a unique solution $\gamma_\vect{v}(t)$, $t\in(-\eps,\eps)$, satisfying initial conditions $\gamma_\vect{v}(0)=p$ and $\gamma'_\vect{v}(0)=\vect{v}$.
On a smooth compact manifold without boundary this solution exists for any $t\in\R$.
\begin{definition}
Let $(M,g)$ be a smooth compact Riemannian manifold without boundary. For any point $p\in M$ the \keyterm{exponential map} $\exp_p:T_p M\to M$ is defined by
\begin{gather*}
\exp_p(\vect{v})=\gamma_\vect{v}(1),
\end{gather*}
where $\gamma_\vect{v}(t)$ is the solution of the equation~\eqref{eq:Geodesics} with the initial conditions $\gamma_\vect{v}(0)=p$ and $\gamma'_\vect{v}(0)=\vect{v}$.
\end{definition}
For non-compact manifolds the exponential map is well defined only in some neighborhood of $0\in T_p M$.

\section{World Function}

The local geometry of a smooth Riemannian manifold $(M,g)$ can be described by the world function (see, for example,~\cite{Synge}).
The world function is a real-valued non-negative function $\sigma$ defined in a neighborhood of the diagonal of $M\times M$ as follows.
Let $x'\in M$ be a fixed point.
We fix a sufficiently small neighborhood of $x'$ in $M$ so that every point $x$ in this neighborhood can be connected by a single geodesic to the point $x'$, which means that the exponential map is injective within this neighborhood.
Let $r(x')$ be the geodesic radius of such a neighborhood.
The \keyterm{injectivity radius} of the manifold $M$, denoted by $r_{\rm inj}(M)$, is defined as the infimum of $r(x')$ over the whole manifold:
\begin{gather*}
r_{\rm inj}(M)=\inf_{x\in M}r(x).
\end{gather*}
We assume that the injectivity radius is strictly positive, i.e. $r_{\rm inj}>0$.
For smooth compact manifolds this is always the case.
For non-compact manifolds, this will be one of our assumptions (among other assumptions specified below). 
The \keyterm{world function} $\sigma(x,x')$ is defined as half the square of the length of the geodesic
connecting the points $x$ and $x'$.

Let $y=y(s)$ be a parametrization of the geodesic connecting the points $x$ and $x'$ with $s$ being the natural parameter (arc length) along the geodesic. Let $y(0)=x$ and $y(\tau)=x'$. Then $\dot y^\alpha(s)=\dfrac{d y^\alpha}{ds}(s)$ is the unit tangent vector at the point $y(t)$ of the geodesic.
The world function satisfies the equations~\cite{Synge}
\begin{gather*}
\sigma(x,x')=\frac{1}{2}\, g^{\alpha\beta}(x) \pdf{\sigma}{x^\alpha}(x,x')\, \pdf{\sigma}{x^\beta}(x,x') \,,
\end{gather*}
and
\begin{gather} \label{eq:TanToGeod}
\tau \dot y^\alpha(0)=-g^{\alpha\beta}(x) \pdf{\sigma}{x^\beta}(x,x') \,.
\end{gather}
It is easy to see that
\begin{gather} \label{eq:ExpOfSigmaAlpha}
\exp_x(\tau \dot y^\alpha(0))=x' \,.
\end{gather}

If two points $x,x'\in M$ cannot be connected by a single geodesic, it is still possible to define the \keyterm{distance} between $x$ and $x'$ as the infimum of the length functional over all smooth curves connecting $x$ and $x'$, if they belong to the same connected component. The distance between two submanifolds $A,B\subset M$ we define as
\begin{gather*}
\dist(A,B)=\inf_{x\in A,\ y\in B} \dist(x,y)\,.
\end{gather*}

\section{The Laplacian on a Riemannian Manifold}

Now we will describe the construction of the Laplace operator on a Riemannian manifold $(M,g)$. For details see~\cite{Avramidi1999, Berline, Rosenberg}.

We define the Hilbert space of square integrable functions $\Ltwo(M)$ to be the completion of the space $\Cinfcomp(M)$ of smooth functions with compact support in the norm induced by the \Ltwo-inner product
\begin{gather*}
\la f,h \ra=\int_{M} \dvol(x)\, f(x)\overline{h(x)} \,.
\end{gather*}

The metric $g$ defines the inner product in every tangent space $T_p M$ and induces in a natural way the inner product in the dual space $T^*_p M$ which we will also denote by $g$. It enables one to define the $\Ltwo$-inner product in the space $\Cinfcomp(T^* M)$ of smooth sections of the cotangent bundle $T^* M$ with compact support  by
\begin{gather*}
\la \alpha,\beta \ra=\int_{M} \dvol(x)\, g(\alpha,\beta)
\end{gather*}
and consider a completion of $\Cinfcomp(T^* M)$ in the norm induced by this inner product. We will denote this completion by $\Ltwo(T^* M)$.

The covariant derivative on the space of smooth functions is the mapping $\nabla:\Cinf(M)\to \Cinf(T^* M)$. Since $\Cinfcomp(M)\subset\Cinf(M)$ and $\Cinfcomp(M)$ is dense in $\Ltwo(M)$, the operator $\nabla$ is densely defined. Thus, we can define the adjoint operator $\nabla^*$ by
\begin{gather*}
\la \alpha, \nabla f \ra = \la \nabla^* \alpha, f \ra,
\end{gather*}
where $\alpha\in \Cinfcomp(T^*M)$ and $f\in \Cinfcomp(M)$.

Now we define the \keyterm{Laplace operator} on $\Cinfcomp(M)$ by
\begin{gather*}
 \Delta=-\nabla^* \nabla.
\end{gather*}
It is a symmetric operator, i.e. $\la\Delta f,h\ra=\la f, \Delta h\ra$ for any $f,h\in\Cinfcomp(M)$. Moreover, it is essentially self-adjoint, which means that there is a unique self-adjoint extension of $\Delta$. To simplify notation, we will denote this extension by the same symbol.
In local coordinates the expression for the Laplacian is
\begin{gather*}
\Delta f = \abs{g}^{-1/2} \pd_j \left( g^{ij} \abs{g}^{1/2} \pd_i f \right).
\end{gather*}
Notice, that the leading symbol of the Laplacian is negative definite:
\begin{gather*} 
\sigma_L (\Delta; x,\xi)=-g^{ij}(x)\xi_i \xi_j<0
\end{gather*}
for any $x\in M$ and any $\xi\neq 0$. Therefore, the leading symbol of the operator
\begin{gather*}
F=-\Delta
\end{gather*}
is positive definite. 

\begin{definition}
We will call a differential operator of order $m$ over an $n$-dimensio\-nal manifold $M$ without boundary \keyterm{elliptic}, if there exist constants $C_1,C_2>0$, such that for any $x\in M$ and any $\xi\in\R^n$ its leading symbol $\sigma_L(x,\xi)$ satisfies the inequality
\begin{gather*}
C_1 \abs{\xi}^m \leqslant \abs{\sigma_L (x,\xi)} \leqslant C_2 \abs{\xi}^m,
\end{gather*}
where $\abs{\xi}^2=\xi_1^2+\ldots+\xi_n^2$.
\end{definition}

Let $H$ be a Hilbert space and $A:H\to H$ be a linear operator on it. Let $\lambda\in\CC$. If the inverse of the operator $A-\lambda$ exists, we will call it the \keyterm{resolvent} and denote it by $R_\lambda=(A-\lambda)^{-1}$. There are the following possibilities~\cite{Kreyszig}:
\begin{enumerate}
\item
$R_\lambda$ exists, is bounded, and is densely defined in $H$.
Then $\lambda$ is called a \keyterm{regular value} of the operator $A$.
The set of all regular values is called the \keyterm{resolvent set} and is denoted by $\rho(A)$.
\item
$R_\lambda$ does not exist. That is, the kernel of the operator $A-\lambda$ is non-trivial, $E_\lambda=\Ker(A-\lambda)\neq\set{0}$.
Then $\lambda$ is called an \keyterm{eigenvalue} of the operator $A$.
The vector space $E_\lambda$ is called the \keyterm{eigenspace} corresponding to the eigenvalue $\lambda$, and any non-zero vector $\phi\in E_\lambda$ is called an \keyterm{eigenvector}, corresponding to the eigenvalue $\lambda$.
If the space $E_\lambda$ is finite-dimensional, then $d(\lambda)=\dim E_\lambda$ is called the \keyterm{multiplicity} of $\lambda$.
The set of all eigenvalues is called the \keyterm{point} (or \keyterm{discrete}) \keyterm{spectrum} and is denoted by $\sigma_p(A)$.
\item
$R_\lambda$ exists and is densely defined in $H$, but is unbounded.
The set of all such $\lambda$ is called the \keyterm{continuous spectrum} and is denoted by $\sigma_c(A)$.
\item
$R_\lambda$ exists, but is not densely defined in $H$.
The set of all such $\lambda$ is called the \keyterm{residual spectrum} and is denoted by $\sigma_r(A)$.
\end{enumerate}
The \keyterm{spectrum} $\sigma(A)$ of the operator $A$ is the union of discrete, continuous, and residual spectra,
\begin{gather*}
\sigma(A)=\sigma_p(A) \cup \sigma_c(A) \cup \sigma_r(A) \,.
\end{gather*}
Clearly, the spectrum is the complement of the resolvent set, $\sigma(A)=\CC\setminus\rho(A)$.
If $H$ is a finite dimensional space, then all points of the spectrum of the operator $A$ are eigenvalues, that is, $\sigma(A)=\sigma_p(A)$, but it may not be the case if $H$ is infinite-dimensional.

Both operators $\Delta$ and $F$ are elliptic. There is the following well-known theorem about the spectrum of an elliptic self-adjoint differential operator with a positive definite leading symbol~\cite{Gilkey}.

\begin{theorem}
Let $(M,g)$ be a compact Riemannian manifold without boundary and $F:\Cinf(M)\to\Cinf(M)$ be an elliptic self-adjoint differential operator with a positive definite leading symbol. Then:
\begin{enumerate}
\item The spectrum of the operator $F$ is real, discrete, and bounded from below.
\item All eigenspaces of the operator $F$ are finite-dimensional.
\item Eigenfunctions of the operator $F$ are smooth and form an orthonormal basis in $\Ltwo(M)$.
\end{enumerate}
\end{theorem}

More generally, a second-order elliptic differential operator acting on smooth sections of a vector bundle is called a \keyterm{Laplace type operator}, if it has a scalar leading symbol.

\section{Heat Kernel and Spectral Functions} \label{se:SpectralFunctions}

In this section we will describe the heat kernel and other spectral functions on a Riemannian manifold $(M,g)$. For details see~\cite{Avramidi1999, Berline, Kirsten, Rosenberg, Vassilevich}.

The \keyterm{heat kernel} $U(t;x,x')$ of a Laplace type operator $F$ on an $n$-dimensional Riemannian manifold without boundary $M$ is the integral kernel of the heat semigroup operator
\begin{gather*}
\e{-t F} : \Ltwo(M)\to\Ltwo(M) \,,
\end{gather*}
and can be defined also as the fundamental solution of the heat equation
\begin{gather*}
(\pd_t+F) U(t;x,x')=0
\end{gather*}
 for $t>0$, with the initial condition
\begin{gather*}
U(0^+;x,x')=\delta(x,x')\,,
\end{gather*}
where $\delta(x,x')$ is the Dirac distribution.
Near the diagonal set of $M\times M$ the heat kernel has the form~\cite{Avramidi1991, Avramidi1999}
\begin{gather} \label{eq:HeatKernel}
U(t;x,x')=(4\pi t)^{-n/2} \varDelta^{1/2}(x,x') \e{-\dfrac{\sig{x,x'}}{2t}} \Omega(t;x,x')\,,
\end{gather}
where $\varDelta(x,x')$ is the Van Vleck-Morette determinant,
\begin{gather} \label{eq:VFMdet}
\varDelta(x,x')
	=\abs{g(x)}^{-1/2}
	\det\!\left( -\nabla_\mu \nabla_{\nu'}' \sigma(x,x') \right)
	\abs{g(x')}^{-1/2}
	\,,
\end{gather}
and the function $\Omega(t;x,x')$ has the asymptotic expansion as $t\to 0^+$
\begin{gather} \label{eq:OmegaAE}
\Omega(t;x,x')\sim\sum_{k=0}^\infty \frac{(-t)^k}{k!}\ \omega_k(x,x')\,.
\end{gather}
(Asymptotic equivalence will be defined in Section~\ref{se:LaplaceMethod}.) 

For a compact manifold $M$ the heat semigroup operator $\exp(-t F)$ is a bounded trace-class operator on the Hilbert space $\Ltwo(M)$, with the trace
\begin{gather*}
\Theta(t)=\Trace \exp(-t F) \,,
\end{gather*}
given by
\begin{gather*}
\Theta(t)=\int_{M} \dvol(x)\, U(t;x,x)=\sum_{m=1}^\infty e^{-t\lambda_m},
\end{gather*}
where each eigenvalue $\lambda_m$, $m\in\N$, is repeated as many times as its multiplicity. This trace is also called the \keyterm{trace of the heat kernel} and is one of the spectral functions of the operator $F$. Some other spectral functions may be expressed in terms of the trace of the heat kernel by means of integral transforms.

The \keyterm{distribution} (or \keyterm{counting}) \keyterm{function} $N(\lambda)$ is defined as the number of eigenvalues below $\lambda\in\R$
\begin{gather*}
N(\lambda)=\abs{\set{\lambda_m\ :\ \lambda_m \leqslant \lambda}} \,,
\end{gather*}
and is given by
\begin{gather*}
N(\lambda)
=\int_{\eps-i\infty}^{\eps+i\infty} \dfrac{dt}{2\pi i}\cdot\dfrac{e^{t\lambda}}{t}\, \Theta(t)
=\sum_{m=1}^\infty \theta(\lambda-\lambda_m)
\end{gather*}
for $\lambda\neq\lambda_m$, $m\in\N$, where $\eps$ is a positive constant and $\theta$ is a Heaviside distribution. For $\lambda=\lambda_m$ we have
\begin{gather*}
N(\lambda_m)=\lim_{\lambda\to\lambda_m^+} N(\lambda).
\end{gather*}

The \keyterm{density function} $\rho(\lambda)$ is defined as the derivative of the distribution function
\begin{gather*}
\rho(\lambda)=\dfrac{d}{d\lambda}N(\lambda)
\end{gather*}
and is given by
\begin{gather*}
\rho(\lambda)
=\int_{\eps-i\infty}^{\eps+i\infty} \dfrac{dt}{2\pi i} e^{t\lambda} \,\Theta(t)
=\sum_{m=1}^\infty \delta(\lambda-\lambda_m) \,.
\end{gather*}

Let $\lambda\in\CC$ be such that $\RE\lambda<\lambda_1$. Then the operator $F-\lambda$ is positive. Let $s\in\CC$ be such that $\RE s> \dfrac{n}{2}$ (this condition is necessary for convergence of the following integral and sum). The \keyterm{generalized $\zeta$-function} $\zeta(s,\lambda)$ is defined as the trace of a complex power of the operator $F-\lambda$
\begin{gather*}
\zeta(s,\lambda)=\Trace (F-\lambda)^s
\end{gather*}
and is given by
\begin{gather*}
\zeta(s,\lambda)
=\dfrac{1}{\Gamma(s)}\int_0^\infty dt\, t^{s-1} e^{t\lambda} \,\Theta(t)
=\sum_{m=1}^\infty \frac{1}{(\lambda_n-\lambda)^s} \,.
\end{gather*}
The generalized $\zeta$-function is an analytical function of $s$ for $\RE s>\dfrac{n}{2}$ and can be continued analytically to a meromorphic function in $\CC$. In particular, it is analytic at $s=0$ and, thus, it is possible to define the \keyterm{functional determinant} of the operator $F-\lambda$ by
\begin{gather*}
\Det(F-\lambda)=\e{-\pd_s\zeta\Big|_{s=0}} .
\end{gather*}

These spectral functions are very useful in studying the spectrum of the operator $F$ and, if known exactly, they determine the spectrum. This is not valid for their asymptotic expansions and there are examples of operators with the same asymptotic series of spectral functions, but different spectra~\cite{Estrada}. In principle, all these functions are equivalent to each other, but the heat kernel is more convenient for practical purposes, since it is a smooth function, while the distribution and density functions are singular.

\section{Laplace Method} \label{se:LaplaceMethod}

In this work we will develop a generalized Laplace method for computing asymptotics of integrals over manifolds. Here we present necessary information for both one- and multi-dimensional Laplace methods. For proofs of cited theorems and lemmas see, for example,~\cite{AvramidiNotes, Fedoryuk}.

\begin{definition}
Functions $f(t)$ and $g(t)$ are \keyterm{asymptotically equivalent} as $t\to 0^+$, denoted by
$f\sim g$, if 
\begin{gather*}
\ds\lim_{t\to 0^+}\dfrac{f(t)-g(t)}{t^N}=0
\end{gather*}
for any $N>0$.
\end{definition}

\begin{definition}
Let $S\in\Cinf([a,b])$. A point $\xi\in (a,b)$ is called a \keyterm{non-degenerate critical point} of the function $S$, if $S'(\xi)=0$ and $S''(\xi)\neq 0$.
\end{definition}

\begin{theorem} \label{th:LaplaceM1d}
Let $\phi,S\in \Cinf([a,b])$. Let the function $S$ attain its global minimum on $[a,b]$ at an interior point $\xi\in(a,b)$. Let the point $\xi$ be the only point where the global minimum is attained, that is, $S(x)>S(\xi)$ for any $x\neq\xi$, and assume that $\xi$ is a non-degenerate critical point.
Then there exists an asymptotic expansion as $t\to 0^+$
\begin{gather*}
\int_a^b dx\, \phi(x) \e{-\frac{1}{t}S(x)} \sim
\e{-\frac{1}{t} S(\xi)} \sum_{k=0}^\infty c_k t^{1/2+k}.
\end{gather*}
The coefficients $c_k$ depend polynomially on derivatives of the function $\phi$ at the point $\xi$,  on derivatives of the function $S$ of order higher than two at the point $\xi$, and on $[S''(\xi)]^{-1/2}$.
\end{theorem}

For computing the coefficients $c_k$ in Theorem~\ref{th:LaplaceM1d}, the following result for standard Gaussian integrals is useful.

\begin{lemma} \label{le:GaussianIntegral1d}
Let $a>0$. Then for any $k\in\Z,\ k\geqslant 0$, there holds
\begin{align*}
\allowdisplaybreaks
&\int_{-\infty}^\infty dx\, e^{-a x^2} x^{2k+1}=0 \,,\\
&\int_{-\infty}^\infty dx\, e^{-a x^2} x^{2k}
=\dfrac{\sqrt{\pi}(2k)!}{k!\,2^{2k}}\,  a^{-1/2-k} \,.
\end{align*}
\end{lemma}

A theorem similar to Theorem~\ref{th:LaplaceM1d} holds in the multi-dimensional case.

\begin{definition}
Let $\Omega\subset\R^n$ be a bounded open set. Let $S\in\Cinf(\Omega)$. A point $\xi\in\Omega$ is called a \keyterm{non-degenerate critical point} of the function $S$, if $\,\nabla S(\xi)=0$ and the Hessian matrix of the function $S$ at the point $\xi$, that is, the matrix with entries $(\pd_i \pd_j S)(\xi)$, is either strictly positive or strictly negative definite. 
\end{definition}

\begin{theorem} \label{th:LaplaceMmd}
Let $\Omega\subset\R^n$ be a bounded open set. Let $\phi,S\in \Cinf(\Omega)$. Let the function $S$ attain its global minimum on $\Omega$ at a point $\xi$. Let the point $\xi$ be the only point where the global minimum is attained and assume that $\xi$ is a non-degenerate critical point. Then there exists an asymptotic expansion as $t\to 0^+$
\begin{gather*}
\int_\Omega dx\, \phi(x) \e{-\frac{1}{t}S(x)} \sim
\e{-\frac{1}{t} S(\xi)} \sum_{k=0}^\infty c_k t^{n/2+k}.
\end{gather*}
The coefficients $c_k$ depend polynomially on derivatives of the function $\phi$ at the point $\xi$,  on derivatives of the function $S$ of order higher than two at the point $\xi$, and on the inverse of the Hessian matrix $Q=\big[(\pd_i\pd_j S)(\xi)\big]^{-1}$.
\end{theorem}

\begin{lemma} \label{le:GaussianIntegralmd}
Let $A$ be a strictly positive symmetric matrix and $G=A^{-1}$. Then for any $k\in\Z$, $k\geqslant 0$, there holds
\begin{align*}
\allowdisplaybreaks
&\int_{\R^n} dx\, \e{-x^T A x} x^{i_1}\ldots x^{i_{2k+1}}=0 \,,\\
&\int_{\R^n} dx\, \e{-x^T A x} x^{i_1}\ldots x^{i_{2k}}
=\dfrac{\pi^{n/2}(2k)!}{k!\,2^{2k}}\, (\det G)^{1/2}\, G^{(i_1 i_2}\ldots G^{i_{2k-1} i_{2k})} \,.
\end{align*}
\end{lemma}

\section{Spectral Asymptotics of Riemannian Submanifolds} \label{se:Framework}

In this section we will introduce deformed spectral functions, following~\cite{Donnelly, Gilkey, Oersted}, and formulate the problem considered in this work.

Let $(M,g)$ be a smooth compact Riemannian manifold without boundary, $\Phi:M\to M$ be an isometric mapping on it, and $\Sigma$ be the fixed point set of the mapping $\Phi$, that is,
\begin{gather} \label{eq:Sigma}
\Sigma=\set{x\in M : \Phi(x)=x}.
\end{gather}
The mapping $\Phi$ naturally defines the operator $\Phi:\Cinf(M)\to\Cinf(M)$ by
\begin{gather*}
\left(\Phi f\right)(x)=f(\Phi(x)) \,.
\end{gather*}
The composition of the operator $\Phi$ and the heat semigroup operator defines an operator
\begin{gather} \label{eq:Uphi_t}
U_\Phi(t)=\Phi\e{t\Delta}
\end{gather}
whose kernel is
\begin{gather} \label{eq:Uphi_txxp}
U_\Phi(t;x,x')=U(t;\Phi(x),x'),
\end{gather}
which we will call the deformed heat kernel.

The \keyterm{deformed heat trace} is
\begin{gather} \label{eq:TrUphi_t}
  \Theta_\Phi(t)=\Trace U_\Phi(t)=\int_{M}\dvol(x)\, U(t;\Phi(x),x) \,.
\end{gather}
It can be used to define other deformed spectral functions of the Laplace operator on the manifold $M$.
Namely, the \keyterm{deformed distribution function}
\begin{gather*}
N_\Phi(\lambda)
=\int_{\eps-i\infty}^{\eps+i\infty} \dfrac{dt}{2\pi i}\cdot\dfrac{e^{t\lambda}}{t}\, \Theta_\Phi(t) \,,
\end{gather*}
the \keyterm{deformed density function}
\begin{gather*}
\rho_\Phi(\lambda)
=\int_{\eps-i\infty}^{\eps+i\infty} \dfrac{dt}{2\pi i} e^{t\lambda} \,\Theta_\Phi(t) \,,
\end{gather*}
the \keyterm{deformed $\zeta$-function}
\begin{gather*}
\zeta_\Phi(s,\lambda)
=\dfrac{1}{\Gamma(s)}\int_0^\infty dt\, t^{s-1} e^{t\lambda} \,\Theta_\Phi(t) \,,
\end{gather*}
and the \keyterm{deformed functional determinant} of the operator $F-\lambda$
\begin{gather*}
\Det_\Phi(F-\lambda)=\e{-\pd_s\zeta_\Phi\Big|_{s=0}} \,,
\end{gather*}
where we have used the same notation as in Section~\ref{se:SpectralFunctions}.

These deformed spectral functions were considered in~\cite{Oersted}.
 The following theorem is known~\cite{Gilkey} for the asymptotic expansion of the deformed heat trace.

\begin{theorem} \label{th:AE-Isometry}
Let $M$ be a smooth compact Riemannian manifold without boundary. Let $\Phi:M\to M$ be an isometry  on $M$. Let $\Sigma$ be the fixed point set of $\Phi$. Then:
\begin{enumerate}
\item[1)] $\Sigma$ is a disjoined union of connected submanifolds $\Sigma_i$, $i=1,\ldots,s$;
\item[2)] as $t\to 0^+$
\begin{gather*}
\Theta_\Phi(t) \sim
\sum_{i=1}^s  \sum_{k=0}^\infty t^{k-m_i/2} \int_{\Sigma_i} \dvol(y)\, a^{(i)}_k(y) \,,
\end{gather*}
where $m_i=\dim{\Sigma_i}$ and $a_k^{(i)}(y)$ are scalar invariants on $\Sigma_i$.
\end{enumerate}
\end{theorem}

Donnelly~\cite{Donnelly} computed the coefficients $a_0$ and $a_1$ of this asymptotic expansion. Asymptotic expansions of this kind are used for proofs of the Atiyah-Singer-Lefschetz formulas for compact group actions and index theorems.

Our goal is to consider a more general problem. Instead of the isometry we will consider a smooth mapping $\Phi$ (not necessarily an isometry) with the fixed point set $\Sigma$, satisfying the following assumptions:
\begin{enumerate}

\item The fixed point set $\Sigma$ consists of finitely many disjoint connected components $\Sigma_i$, $i=1,\ldots,s$.

\item Each connected component $\Sigma_i$ is a smooth compact submanifold of $M$ without boundary.

\item Each connected component $\Sigma_i$ is a ``non-degenerate'' fixed point set. We will formulate this assumption precisely for all considered cases in Section~\ref{se:PlaneCurve} (p.~\pageref{se:PlaneCurve}), Section~\ref{se:PlanePoint} (p.~\pageref{se:PlanePoint}), and Section~\ref{se:CurvedPoint} (p.~\pageref{se:CurvedPoint}).

\end{enumerate}
A theorem analogous to Theorem~\ref{th:AE-Isometry} holds for such mappings as well~\cite{Gilkey}. All these assumptions are automatically satisfied for isometries. 

We will develop a technique for computation of the asymptotic expansion coefficients and compute explicitly the coefficients $a_0$, $a_1$, and $a_2$.

We will say that two subsets $A$ and $B$ of the manifold $M$ are separated if
$\dist(A,B)>0$.
It follows from the above assumptions, that all components $\Sigma_i$ are separated, that is,
\begin{gather*}
\min_{i,j} \dist\left( \Sigma_i,\Sigma_j \right)>0.
\end{gather*}
Since the asymptotic expansion of the deformed heat trace $\Theta_\Phi(t)$ is determined by a neighborhood of $\Sigma$ and all connected components $\Sigma_i$ are isolated, it is possible to consider each component separately. The whole asymptotic expansion is equal to the sum of asymptotic expansions generated by each component.
So, without loss of generality and in order to simplify notation, we will assume that $\Sigma$ consists of only one connected component.

We will also consider non-compact manifolds $M$. In this case our additional assumption will be that the injectivity radius of $M$ is non-zero.

\chapter{Symmetric Tensors} \label{ch:SymmetricTensors}

In this chapter we will derive reduction formulas for contraction of symmetric tensors, that allow us to express multi-dimensional Gaussian integrals in a symmetrization-free form. Our main result in this chapter is Lemma~\ref{le:Contraction}, which gives an explicit expression of this kind for a general case.

\section{Algebra of Symmetric Tensors}

Let $\vee: V_p \times V_q \to V_{p+q}$ be the \keyterm{symmetrized tensor product} of symmetric covariant tensors, defined by
\begin{gather*}
  A \vee B= \Sym (A \otimes B) \,.
\end{gather*}

Let $Q$ be a symmetric $(2,0)$-tensor. Let $\tau: V_{m+2} \to V_{m}$ be the operator of contraction of symmetric covariant tensors with $Q$, defined by
\begin{gather}
  \label{eq:contraction_def}
  (\tau A)_{\nu_1\dots\nu_{m}}=Q^{\alpha\beta} A_{\alpha\beta\nu_1\dots\nu_{m}} \,.
\end{gather}
By definition, $\tau A=0$ for any $A$ from $V_0$ or $V_1$.

A \keyterm{multi-index} $\mi{n}$ of the order $k$ is a $k$-tuple $(n_1,\ldots,n_k)$ of non-negative integers. We will use the following notation:
\begin{align*}
\abs{\mi{n}}&=n_1+n_2+\ldots+ n_k\,,\\
\mi{n}!&=n_1!\, n_2!\, \ldots n_k! \,.
\end{align*}

\section{Contraction Formulas: Theory}

Let $k\in\N$ and $\mi{n}$ be a multi-index of the order $k$, such that $\abs{\mi{n}}=2N$.
Let $A^{(i)}\in V_{n_{i}}$, $i=1,\ldots,k$, be symmetric tensors. We will be interested in explicit formulas for the full contraction of the symmetrized product of these tensors with the tensor $Q$, that is, in scalar expressions of the form
\begin{gather} \label{eq:FullContraction}
\tau^N \left(A^{(1)}\vee \ldots \vee A^{(k)}\right).
\end{gather}
By ``explicit formulas'' we mean expressions without symmetrization of indices.
We will develop a technique for obtaining such formulas and then apply it to some particular cases that will be needed in the following chapters.

By definition, we have
\begin{gather} \label{eq:2.3}
  \tau^N \left(A^{(1)}\vee \ldots \vee A^{(k)}\right)=\frac{1}{(2N)!} \sum_{\phi\in S_{2N}} T(\phi),
\end{gather}
where
\begin{gather*}
   T(\phi)= \big(\underbrace{Q \otimes \ldots \otimes Q}_{N}\big)^{\nu_1 \ldots \nu_{2N}} 
  \left( A^{(1)} \otimes \ldots \otimes A^{(k)} \right)_{\nu_{\phi(1)}\ldots\nu_{\phi(2N)}}.
\end{gather*}
We need to combine all similar terms in the sum~\eqref{eq:2.3}, that arise because of the symmetry of the tensors $A^{(i)}$ and $Q$, as well as the possibility to permute factors in the tensor product of $Q$'s. 

With each term $T(\phi)$ we will associate a graph $\Gamma=H(\phi)$. Every graph $\Gamma$ has $k$ distinct vertices of degrees $n_i$, corresponding to the tensors $A^{(i)}$, and $N$ edges, corresponding to the factors of the tensor $Q$. An edge between vertices $i$ and $j$ means that an index of the tensor $A^{(i)}$ is contracted with an index of the tensor $A^{(j)}$ by one of the factors $Q$. Let $\mf{G}(\mi{n})$ be the set of all such graphs. For every permutation $\phi\in S_{2N}$ the term $T(\phi)$ has a unique graph $\Gamma=H(\phi)$ associated with it. Thus, the mapping
\begin{gather*}
  H:S_{2N} \to \mf{G}(\mi{n})
\end{gather*}
is well defined.
From the definition of the mapping $H$ it follows:
\begin{enumerate}
\item Given a graph $\Gamma \in \mf{G}(\mi{n})$ it is possible to construct many terms $T(\phi)$ or, in other words, to find many permutations $\phi\in S_{2N}$, such that $H(\phi)=\Gamma$.
\item For any graph $\Gamma \in \mf{G}(\mi{n})$ all permutations $\phi \in H^{-1}(\Gamma)$ generate similar terms $T(\phi)$ in the sum~\eqref{eq:2.3}.
\item For any two different graphs $\Gamma_1, \Gamma_2 \in \mf{G}(\mi{n})$ and any two permutations $\phi_1 \in H^{-1}(\Gamma_1)$, $\phi_2\in H^{-1}(\Gamma_2)$ terms $T(\phi_1)$ and $T(\phi_2)$ are not similar.
\end{enumerate}
So, the mapping $H$ is surjective, but not injective. Therefore, we have a partition of $S_{2N}$ into equivalence classes
\begin{gather*}
S_{2N}=\bigcup_{\Gamma\in\mf{G}(\mi{n})} H^{-1}(\Gamma) \,.
\end{gather*}
All permutations $\phi$ from the same class $H^{-1}(\Gamma)$ generate similar terms $T(\phi)$.

Our goal is now to list all elements of $\mf{G}(\mi{n})$ and to count how many similar terms $T(\phi)$ in the sum~\eqref{eq:2.3} correspond to each graph. That is, for any graph $\Gamma\in\mf{G}(\mi{n})$ we want to find the number of elements in the equivalence class $H^{-1}(\Gamma)$. We will denote it by
\begin{gather*}
\tld{c}(\Gamma)=\abs{H^{-1}(\Gamma)} .
\end{gather*}

We will identify a graph $\Gamma\in\mf{G}(\mi{n})$ with its adjacency matrix $\Gamma=(\gamma_{ij})$, which is a square matrix of order $k$. By definition, its entry $\gamma_{ij}$ is equal to the number of edges between the vertices $i$ and $j$, that is, to the number of indices of the tensor $A^{(i)}$ contracted with indices of the tensor $A^{(j)}$ by factors of $Q$. It is easy to see, that every matrix $\Gamma\in\mf{G}(\mi{n})$ has the following properties:
\begin{enumerate}
\item For any $i$ and $j$ the entry $\gamma_{ij}$ is a non-negative integer from $0$ to $2N$.
\item The matrix $\Gamma$ is symmetric, that is, for any $i$ and $j$ there holds $\gamma_{ij}=\gamma_{ji}$.
\item For any $i$ the diagonal coefficient $\gamma_{ii}$ is even.
\item For any $i$ the sum of all entries of the $i$-th row, as well as the sum of all entries of the $i$-th column, is equal to the rank of the tensor $A^{(i)}$:
\begin{gather*}
  \sum_{j=1}^k \gamma_{ij}=\sum_{j=1}^k \gamma_{ji}=n_i.
\end{gather*}
\end{enumerate}
Conversely, every matrix of the order $k$ having all these properties is an element of the set $\mf{G}(\mi{n})$. 

Let $\Gamma\in\mf{G}(\mi{n})$. For computing $\tld{c}(\Gamma)$ we select an arbitrary element $\phi_\Gamma \in H^{-1}(\Gamma)$ and the corresponding term $T(\phi)$. We can generate new terms associated with $\Gamma$ by using symmetry of the tensors $Q$ and $A^{(i)}$, and permuting factors in the tensor product of $Q$'s. It gives us $2^N N!\, n_1!\ldots n_k!$ terms. However, some terms will differ just by names of dummy indices. Below we compute how many of these generated terms are the same (up to names of dummy indices).
\begin{enumerate}

\item If $\gamma_{ii}\geqslant 2$, then simultaneous permutation of two indices of the contracted pair in the tensors $Q$ and $A^{(i)}$ does not change the term. For example:
\begin{gather*}
  Q^{\alpha\beta}A^{(i)}_{\alpha\beta\ldots}
  =
  Q^{\beta\alpha}A^{(i)}_{\beta\alpha\ldots} \,.
\end{gather*}
It reduces the number of terms by a factor of $2^{\gamma_{ii}/2}$. Note, that $\dfrac{\gamma_{ii}}{2}$ is the number of self-loops at the vertex $i$, and $\dfrac{1}{2}\tr\Gamma$ is the number of self-loops in the whole graph $\Gamma$.

\item If $\gamma_{ii}\geqslant 4$, then simultaneous permutation of $\dfrac{\gamma_{ii}}{2}$ contracting factors of $Q$ and $\dfrac{\gamma_{ii}}{2}$ contracted pairs of indices of the tensor $A^{(i)}$ does not change the term. For example:
\begin{gather*}
  Q^{\alpha\beta}Q^{\gamma\delta}A^{(i)}_{\alpha\beta\gamma\delta\ldots}
  =
  Q^{\gamma\delta}Q^{\alpha\beta}A^{(i)}_{\gamma\delta\alpha\beta\ldots} \,.
\end{gather*}
It reduces the number of terms by a factor of $\left(\dfrac{\gamma_{ii}}{2}\right)!$.

\item If $\gamma_{ij}\geqslant 2$ for $i<j$, then simultaneous permutation of $\gamma_{ij}$ indices of the tensors $A^{(i)}$ and $A^{(j)}$ contracted with the tensor $Q$, and factors of $Q$, does not change the term. For example:
\begin{gather*}
  Q^{\alpha\gamma}Q^{\beta\delta}A^{(i)}_{\alpha\beta\ldots}A^{(j)}_{\gamma\delta\ldots}
  =
  Q^{\beta\delta}Q^{\alpha\gamma}A^{(i)}_{\beta\alpha\ldots}A^{(j)}_{\delta\gamma\ldots} \,.
\end{gather*}
It reduces the number of terms by a factor of $\gamma_{ij}!$.

\end{enumerate}
Taking all this into account, the number of terms associated with $\Gamma$ is
\begin{gather*}
\tld{c}(\Gamma)
=\frac{2^N N!\, n_1!\ldots n_k!}{\ds \prod_{1\leqslant i<j\leqslant k} \gamma_{ij}! \prod_{1\leqslant i \leqslant k} \left(\frac{\gamma_{ii}}{2}\right)!\, 2^{(\tr\Gamma)/2} } \,.
\end{gather*}
Note, that these coefficients have the following property:
\begin{gather*}
\sum_{\Gamma\in\mf{G}(\mi{n})} \tld{c}(\Gamma)=(2N)! \,.
\end{gather*}
So, it is natural to define weights or probabilities of graphs $\Gamma$ by
\begin{gather*}
  c(\Gamma)=\frac{\tld{c}(\Gamma)}{(2N)!}
  =\frac{\mi{n}!}{(\abs{\mi{n}}-1)!!} \cdot \frac{1}{\ds \prod_{1\leqslant i<j\leqslant k} \gamma_{ij}! \prod_{1\leqslant i \leqslant k} \gamma_{ii}!!} \,,
\end{gather*}
and they have the property
\begin{gather*}
\sum_{\Gamma\in\mf{G}(\mi{n})} c(\Gamma)=1 \,.
\end{gather*}

Thus, we have proved the following lemma.

\begin{lemma} \label{le:Contraction}
Let $k, N\in \N$. Let $\mi{n}$ be a multi-index of the order $k$, such that $\abs{\mi{n}}=2N$.
Let $A^{(i)}\in V_{n_{i}}$, $i=1,\ldots,k$, be symmetric $(0,n_i)$-tensors. Then
\begin{gather*}
\tau^N \left(A^{(1)}\vee \ldots \vee A^{(k)}\right)
=\sum_{\Gamma\in\mf{G}(\mi{n})} c(\Gamma) T(\phi_\Gamma) ,
\end{gather*}
where $\phi_\Gamma$ is an arbitrary permutation from the set $H^{-1}(\Gamma)$ and the coefficients $c(\Gamma)$ are given by
\begin{gather*}
  c(\Gamma)
  =\frac{\mi{n}!}{(\abs{\mi{n}}-1)!!} \cdot \frac{1}{\ds \prod_{1\leqslant i<j\leqslant k} \gamma_{ij}! \prod_{1\leqslant i \leqslant k} \gamma_{ii}!!} \,.
\end{gather*}
\end{lemma}

\section{Contraction Formulas: Particular Cases}

Now we will apply Lemma~\ref{le:Contraction} to a few particular cases. Matrices in the following corollaries are listed in the ``decreasing row-wise lexicographical order''.

\begin{corollary} \label{co:Contraction33}
Let $A, B \in V_3$. Then
\begin{align*}
  \tau^3 (A\vee B)&=
  \frac{1}{5}\, \Big( 3\, T_{(1)} + 2\, T_{(2)} \Big)  \,,
\end{align*}
where $T_{(i)}$ are terms associated with the adjacency matrices $\Gamma_{(i)}$ listed below:
\begin{align*}
\allowdisplaybreaks
  \Gamma_{(1)}&=
  \begin{pmatrix}
    2 & 1\\
    1 & 2
  \end{pmatrix},
  &
  \Gamma_{(2)}&=
  \begin{pmatrix}
    0 & 3\\
    3 & 0
  \end{pmatrix}.
\end{align*}
\end{corollary}

\begin{corollary} \label{co:Contraction35}
Let $A\in V_3$, $B\in V_5$. Then
\begin{align*}
  \tau^4 (A\vee B)&=
  \frac{1}{7}\, \Big( 3\, T_{(1)} + 4\, T_{(2)} \Big)  \,,
\end{align*}
where $T_{(i)}$ are terms associated with the adjacency matrices $\Gamma_{(i)}$ listed below:
\begin{align*}
\allowdisplaybreaks
  \Gamma_{(1)}&=
  \begin{pmatrix}
    2 & 1\\
    1 & 4
  \end{pmatrix},
  &
  \Gamma_{(2)}&=
  \begin{pmatrix}
    0 & 3\\
    3 & 2
  \end{pmatrix}.
\end{align*}
\end{corollary}

\begin{corollary} \label{co:Contraction44}
Let $A\in V_4$, $B\in V_4$. Then
\begin{align*}
  \tau^4 (A\vee B)&=
  \frac{1}{35}\, \Big( 3\, T_{(1)} +24\, T_{(2)} + 8\, T_{(3)} \Big)  \,,
\end{align*}
where $T_{(i)}$ are terms associated with the adjacency matrices $\Gamma_{(i)}$ listed below:
\begin{align*}
\allowdisplaybreaks
  \Gamma_{(1)}&=
  \begin{pmatrix}
    4 & 0\\
    0 & 4
  \end{pmatrix},
  &
  \Gamma_{(2)}&=
  \begin{pmatrix}
    2 & 2\\
    2 & 2
  \end{pmatrix},
  &
  \Gamma_{(3)}&=
  \begin{pmatrix}
    0 & 4\\
    4 & 0
  \end{pmatrix}.
\end{align*}
\end{corollary}

\begin{corollary} \label{co:Contraction334}
Let $A, B \in V_3$, $C\in V_4$. Then
\begin{align*}
  \tau^5 (A\vee B\vee C)=
  \frac{1}{105}\, \Big(
  3\, T_{(1)}
  &+12\, T_{(2)}
  +12\, T_{(3)}
  +8\, T_{(4)}
  +2\, T_{(5)}\\
  &+24\, T_{(6)}
  +12\, T_{(7)}
  +24\, T_{(8)}
  +8\, T_{(9)} \Big)  \,,
\end{align*}
where $T_{(i)}$ are terms associated with the adjacency matrices $\Gamma_{(i)}$ listed below:
\begin{align*}
\allowdisplaybreaks
  \Gamma_{(1)}&=
  \begin{pmatrix}
    2 & 1 & 0\\
    1 & 2 & 0\\
    0 & 0 & 4
  \end{pmatrix},
  &
  \Gamma_{(2)}&=
  \begin{pmatrix}
    2 & 1 & 0\\
    1 & 0 & 2\\
    0 & 2 & 2
  \end{pmatrix},
  &
  \Gamma_{(3)}&=
  \begin{pmatrix}
    2 & 0 & 1\\
    0 & 2 & 1\\
    1 & 1 & 2
  \end{pmatrix},
  \\[2ex]
  \Gamma_{(4)}&=
  \begin{pmatrix}
    2 & 0 & 1\\
    0 & 0 & 3\\
    1 & 3 & 0
  \end{pmatrix},
  &
  \Gamma_{(5)}&=
  \begin{pmatrix}
    0 & 3 & 0\\
    3 & 0 & 0\\
    0 & 0 & 4
  \end{pmatrix},
  &
  \Gamma_{(6)}&=
  \begin{pmatrix}
    0 & 2 & 1\\
    2 & 0 & 1\\
    1 & 1 & 2
  \end{pmatrix},
  \\[2ex]
  \Gamma_{(7)}&=
  \begin{pmatrix}
    0 & 1 & 2\\
    1 & 2 & 0\\
    2 & 0 & 2
  \end{pmatrix},
  &
  \Gamma_{(8)}&=
  \begin{pmatrix}
    0 & 1 & 2\\
    1 & 0 & 2\\
    2 & 2 & 0
  \end{pmatrix},
  &
  \Gamma_{(9)}&=
  \begin{pmatrix}
    0 & 0 & 3\\
    0 & 2 & 1\\
    3 & 1 & 0
  \end{pmatrix}.
\end{align*}
\end{corollary}

\begin{corollary} \label{co:Contraction3333}
Let $A, B, C, D \in V_3$. Then
\begin{multline} \label{eq:Contraction3333}
  \tau^6 (A\vee B\vee C\vee D)=
  \frac{1}{1155}\, \Big(
   9\, T_{(1)}
  +6\, T_{(2)}
  +18\, T_{(3)}
  +18\, T_{(4)}
  +36\, T_{(5)}
  \\
\begin{aligned}
  &
  +18\, T_{(6)}
  +18\, T_{(7)}
  +9\, T_{(8)}
  +18\, T_{(9)}
  +18\, T_{(10)}
  +36\, T_{(11)}
  +6\, T_{(12)}
  \\&
  +9\, T_{(13)}
  +18\, T_{(14)}
  +18\, T_{(15)}
  +6\, T_{(16)}
  +36\, T_{(17)}
  +18\, T_{(18)}
  +6\, T_{(19)}
  \\&
  +4\, T_{(20)}
  +36\, T_{(21)}
  +18\, T_{(22)}
  +36\, T_{(23)}
  +18\, T_{(24)}
  +36\, T_{(25)}
  +36\, T_{(26)}
  \\&
  +18\, T_{(27)}
  +36\, T_{(28)}
  +36\, T_{(29)}
  +18\, T_{(30)}
  +36\, T_{(31)}
  +36\, T_{(32)}
  +144\, T_{(33)}
  \\&
  +36\, T_{(34)}
  +18\, T_{(35)}
  +36\, T_{(36)}
  +36\, T_{(37)}
  +6\, T_{(38)}
  +4\, T_{(39)}
  +18\, T_{(40)}
  \\&
  +36\, T_{(41)}
  +36\, T_{(42)}
  +36\, T_{(43)}
  +18\, T_{(44)}
  +36\, T_{(45)}
  +6\, T_{(46)}
  +4\, T_{(47)}
  \Big)  \,,
\end{aligned}
\end{multline}
where $T_{(i)}$ are terms associated with the adjacency matrices $\Gamma_{(i)}$ listed below:
\begin{align*}
\allowdisplaybreaks
  \Gamma_{(1)}&=
  \begin{pmatrix}
    2 & 1 & 0 & 0\\
    1 & 2 & 0 & 0\\
    0 & 0 & 2 & 1\\
    0 & 0 & 1 & 2
  \end{pmatrix},
  &
  \Gamma_{(2)}&=
  \begin{pmatrix}
    2 & 1 & 0 & 0\\
    1 & 2 & 0 & 0\\
    0 & 0 & 0 & 3\\
    0 & 0 & 3 & 0
  \end{pmatrix},
  &
  \Gamma_{(3)}&=
  \begin{pmatrix}
    2 & 1 & 0 & 0\\
    1 & 0 & 2 & 0\\
    0 & 2 & 0 & 1\\
    0 & 0 & 1 & 2
  \end{pmatrix},
  \displaybreak[0]\\[2ex]
  \Gamma_{(4)}&=
  \begin{pmatrix}
    2 & 1 & 0 & 0\\
    1 & 0 & 1 & 1\\
    0 & 1 & 2 & 0\\
    0 & 1 & 0 & 2
  \end{pmatrix},
  &
  \Gamma_{(5)}&=
  \begin{pmatrix}
    2 & 1 & 0 & 0\\
    1 & 0 & 1 & 1\\
    0 & 1 & 0 & 2\\
    0 & 1 & 2 & 0
  \end{pmatrix},
  &
  \Gamma_{(6)}&=
  \begin{pmatrix}
    2 & 1 & 0 & 0\\
    1 & 0 & 0 & 2\\
    0 & 0 & 2 & 1\\
    0 & 2 & 1 & 0
  \end{pmatrix},
  \displaybreak[0]\\[2ex]
  \Gamma_{(7)}&=
  \begin{pmatrix}
    2 & 0 & 1 & 0\\
    0 & 2 & 1 & 0\\
    1 & 1 & 0 & 1\\
    0 & 0 & 1 & 2
  \end{pmatrix},
  &
  \Gamma_{(8)}&=
  \begin{pmatrix}
    2 & 0 & 1 & 0\\
    0 & 2 & 0 & 1\\
    1 & 0 & 2 & 0\\
    0 & 1 & 0 & 2
  \end{pmatrix},
  &
  \Gamma_{(9)}&=
  \begin{pmatrix}
    2 & 0 & 1 & 0\\
    0 & 2 & 0 & 1\\
    1 & 0 & 0 & 2\\
    0 & 1 & 2 & 0
  \end{pmatrix},
  \displaybreak[0]\\[2ex]
  \Gamma_{(10)}&=
  \begin{pmatrix}
    2 & 0 & 1 & 0\\
    0 & 0 & 2 & 1\\
    1 & 2 & 0 & 0\\
    0 & 1 & 0 & 2
  \end{pmatrix},
  &
  \Gamma_{(11)}&=
  \begin{pmatrix}
    2 & 0 & 1 & 0\\
    0 & 0 & 1 & 2\\
    1 & 1 & 0 & 1\\
    0 & 2 & 1 & 0
  \end{pmatrix},
  &
  \Gamma_{(12)}&=
  \begin{pmatrix}
    2 & 0 & 1 & 0\\
    0 & 0 & 0 & 3\\
    1 & 0 & 2 & 0\\
    0 & 3 & 0 & 0
  \end{pmatrix},
  \displaybreak[0]\\[2ex]
  \Gamma_{(13)}&=
  \begin{pmatrix}
    2 & 0 & 0 & 1\\
    0 & 2 & 1 & 0\\
    0 & 1 & 2 & 0\\
    1 & 0 & 0 & 2
  \end{pmatrix},
  &
  \Gamma_{(14)}&=
  \begin{pmatrix}
    2 & 0 & 0 & 1\\
    0 & 2 & 1 & 0\\
    0 & 1 & 0 & 2\\
    1 & 0 & 2 & 0
  \end{pmatrix},
  &
  \Gamma_{(15)}&=
  \begin{pmatrix}
    2 & 0 & 0 & 1\\
    0 & 2 & 0 & 1\\
    0 & 0 & 2 & 1\\
    1 & 1 & 1 & 0
  \end{pmatrix},
  \displaybreak[0]\\[2ex]
  \Gamma_{(16)}&=
  \begin{pmatrix}
    2 & 0 & 0 & 1\\
    0 & 0 & 3 & 0\\
    0 & 3 & 0 & 0\\
    1 & 0 & 0 & 2
  \end{pmatrix},
  &
  \Gamma_{(17)}&=
  \begin{pmatrix}
    2 & 0 & 0 & 1\\
    0 & 0 & 2 & 1\\
    0 & 2 & 0 & 1\\
    1 & 1 & 1 & 0
  \end{pmatrix},
  &
  \Gamma_{(18)}&=
  \begin{pmatrix}
    2 & 0 & 0 & 1\\
    0 & 0 & 1 & 2\\
    0 & 1 & 2 & 0\\
    1 & 2 & 0 & 0
  \end{pmatrix},
  \displaybreak[0]\\[2ex]
  \Gamma_{(19)}&=
  \begin{pmatrix}
    0 & 3 & 0 & 0\\
    3 & 0 & 0 & 0\\
    0 & 0 & 2 & 1\\
    0 & 0 & 1 & 2
  \end{pmatrix},
  &
  \Gamma_{(20)}&=
  \begin{pmatrix}
    0 & 3 & 0 & 0\\
    3 & 0 & 0 & 0\\
    0 & 0 & 0 & 3\\
    0 & 0 & 3 & 0
  \end{pmatrix},
  &
  \Gamma_{(21)}&=
  \begin{pmatrix}
    0 & 2 & 1 & 0\\
    2 & 0 & 1 & 0\\
    1 & 1 & 0 & 1\\
    0 & 0 & 1 & 2
  \end{pmatrix},
  \displaybreak[0]\\[2ex]
  \Gamma_{(22)}&=
  \begin{pmatrix}
    0 & 2 & 1 & 0\\
    2 & 0 & 0 & 1\\
    1 & 0 & 2 & 0\\
    0 & 1 & 0 & 2
  \end{pmatrix},
  &
  \Gamma_{(23)}&=
  \begin{pmatrix}
    0 & 2 & 1 & 0\\
    2 & 0 & 0 & 1\\
    1 & 0 & 0 & 2\\
    0 & 1 & 2 & 0
  \end{pmatrix},
  &
  \Gamma_{(24)}&=
  \begin{pmatrix}
    0 & 2 & 0 & 1\\
    2 & 0 & 1 & 0\\
    0 & 1 & 2 & 0\\
    1 & 0 & 0 & 2
  \end{pmatrix},
  \displaybreak[0]\\[2ex]
  \Gamma_{(25)}&=
  \begin{pmatrix}
    0 & 2 & 0 & 1\\
    2 & 0 & 1 & 0\\
    0 & 1 & 0 & 2\\
    1 & 0 & 2 & 0
  \end{pmatrix},
  &
  \Gamma_{(26)}&=
  \begin{pmatrix}
    0 & 2 & 0 & 1\\
    2 & 0 & 0 & 1\\
    0 & 0 & 2 & 1\\
    1 & 1 & 1 & 0
  \end{pmatrix},
  &
  \Gamma_{(27)}&=
  \begin{pmatrix}
    0 & 1 & 2 & 0\\
    1 & 2 & 0 & 0\\
    2 & 0 & 0 & 1\\
    0 & 0 & 1 & 2
  \end{pmatrix},
  \displaybreak[0]\\[2ex]
  \Gamma_{(28)}&=
  \begin{pmatrix}
    0 & 1 & 2 & 0\\
    1 & 0 & 1 & 1\\
    2 & 1 & 0 & 0\\
    0 & 1 & 0 & 2
  \end{pmatrix},
  &
  \Gamma_{(29)}&=
  \begin{pmatrix}
    0 & 1 & 2 & 0\\
    1 & 0 & 0 & 2\\
    2 & 0 & 0 & 1\\
    0 & 2 & 1 & 0
  \end{pmatrix},
  &
  \Gamma_{(30)}&=
  \begin{pmatrix}
    0 & 1 & 1 & 1\\
    1 & 2 & 0 & 0\\
    1 & 0 & 2 & 0\\
    1 & 0 & 0 & 2
  \end{pmatrix},
  \displaybreak[0]\\[2ex]
  \Gamma_{(31)}&=
  \begin{pmatrix}
    0 & 1 & 1 & 1\\
    1 & 2 & 0 & 0\\
    1 & 0 & 0 & 2\\
    1 & 0 & 2 & 0
  \end{pmatrix},
  &
  \Gamma_{(32)}&=
  \begin{pmatrix}
    0 & 1 & 1 & 1\\
    1 & 0 & 2 & 0\\
    1 & 2 & 0 & 0\\
    1 & 0 & 0 & 2
  \end{pmatrix},
  &
  \Gamma_{(33)}&=
  \begin{pmatrix}
    0 & 1 & 1 & 1\\
    1 & 0 & 1 & 1\\
    1 & 1 & 0 & 1\\
    1 & 1 & 1 & 0
  \end{pmatrix},
  \displaybreak[0]\\[2ex]
  \Gamma_{(34)}&=
  \begin{pmatrix}
    0 & 1 & 1 & 1\\
    1 & 0 & 0 & 2\\
    1 & 0 & 2 & 0\\
    1 & 2 & 0 & 0
  \end{pmatrix},
  &
  \Gamma_{(35)}&=
  \begin{pmatrix}
    0 & 1 & 0 & 2\\
    1 & 2 & 0 & 0\\
    0 & 0 & 2 & 1\\
    2 & 0 & 1 & 0
  \end{pmatrix},
  &
  \Gamma_{(36)}&=
  \begin{pmatrix}
    0 & 1 & 0 & 2\\
    1 & 0 & 2 & 0\\
    0 & 2 & 0 & 1\\
    2 & 0 & 1 & 0
  \end{pmatrix},
  \displaybreak[0]\\[2ex]
  \Gamma_{(37)}&=
  \begin{pmatrix}
    0 & 1 & 0 & 2\\
    1 & 0 & 1 & 1\\
    0 & 1 & 2 & 0\\
    2 & 1 & 0 & 0
  \end{pmatrix},
  &
  \Gamma_{(38)}&=
  \begin{pmatrix}
    0 & 0 & 3 & 0\\
    0 & 2 & 0 & 1\\
    3 & 0 & 0 & 0\\
    0 & 1 & 0 & 2
  \end{pmatrix},
  &
  \Gamma_{(39)}&=
  \begin{pmatrix}
    0 & 0 & 3 & 0\\
    0 & 0 & 0 & 3\\
    3 & 0 & 0 & 0\\
    0 & 3 & 0 & 0
  \end{pmatrix},
  \displaybreak[0]\\[2ex]
  \Gamma_{(40)}&=
  \begin{pmatrix}
    0 & 0 & 2 & 1\\
    0 & 2 & 1 & 0\\
    2 & 1 & 0 & 0\\
    1 & 0 & 0 & 2
  \end{pmatrix},
  &
  \Gamma_{(41)}&=
  \begin{pmatrix}
    0 & 0 & 2 & 1\\
    0 & 2 & 0 & 1\\
    2 & 0 & 0 & 1\\
    1 & 1 & 1 & 0
  \end{pmatrix},
  &
  \Gamma_{(42)}&=
  \begin{pmatrix}
    0 & 0 & 2 & 1\\
    0 & 0 & 1 & 2\\
    2 & 1 & 0 & 0\\
    1 & 2 & 0 & 0
  \end{pmatrix},
  \displaybreak[0]\\[2ex]
  \Gamma_{(43)}&=
  \begin{pmatrix}
    0 & 0 & 1 & 2\\
    0 & 2 & 1 & 0\\
    1 & 1 & 0 & 1\\
    2 & 0 & 1 & 0
  \end{pmatrix},
  &
  \Gamma_{(44)}&=
  \begin{pmatrix}
    0 & 0 & 1 & 2\\
    0 & 2 & 0 & 1\\
    1 & 0 & 2 & 0\\
    2 & 1 & 0 & 0
  \end{pmatrix},
  &
  \Gamma_{(45)}&=
  \begin{pmatrix}
    0 & 0 & 1 & 2\\
    0 & 0 & 2 & 1\\
    1 & 2 & 0 & 0\\
    2 & 1 & 0 & 0
  \end{pmatrix},
  \displaybreak[0]\\[2ex]
  \Gamma_{(46)}&=
  \begin{pmatrix}
    0 & 0 & 0 & 3\\
    0 & 2 & 1 & 0\\
    0 & 1 & 2 & 0\\
    3 & 0 & 0 & 0
  \end{pmatrix},
  &
  \Gamma_{(47)}&=
  \begin{pmatrix}
    0 & 0 & 0 & 3\\
    0 & 0 & 3 & 0\\
    0 & 3 & 0 & 0\\
    3 & 0 & 0 & 0
  \end{pmatrix}.
\end{align*}
\end{corollary}

We will also need formulas for contractions of equal tensors.

\begin{corollary} \label{co:Contraction3x4}
Let $A\in V_3$. Then
\begin{multline} \label{eq:Contraction3x4}
  \tau^6 (A\vee A\vee A\vee A)
  =
  \frac{1}{385}\, \Big(
  9\, T_{(1)}
  +12\, T_{(2)}
  +72\, T_{(3)}
  +24\, T_{(4)}\\
  +144\, T_{(5)}
  +4\, T_{(20)}
  +72\, T_{(23)}
  +48\, T_{(33)} \Big)  \,,
\end{multline}
where $T_{(i)}$ are defined as in Corollary~\ref{co:Contraction3333}. 
\end{corollary}

\begin{proof}
When some or all tensors in the contraction~\eqref{eq:Contraction3333} are equal, some adjacency matrices generate similar terms. Namely, if tensors $A^{(i)}$ and $A^{(j)}$ in~\eqref{eq:FullContraction} are equal, and matrices $\Gamma_{(1)}$ and $\Gamma_{(2)}$ can be obtained one from another by exchanging $i$-th and $j$-th rows, and $i$-th and $j$-th columns, then $\Gamma_{(1)}$ and $\Gamma_{(2)}$ generate similar terms. 
By direct comparison, the matrices $\Gamma_{(i)}$ from Corollary~\ref{co:Contraction3333} split into eight groups, such that all matrices in each group generate similar terms. Numbers $i$ for groups of equivalent matrices $\Gamma_{(i)}$ are listed below:
\begin{center}
\begin{tabular}{ll}
Group 1: & 1, 8, 13; \\
Group 2: & 2, 12, 16, 19, 38, 46; \\
Group 3: & 3, 6, 9, 10, 14, 18, 22, 24, 27, 35, 40, 44; \\
Group 4: & 4, 7, 15, 30; \\
Group 5: & 5, 11, 17, 21, 26, 28, 31, 32, 34, 37, 41, 43; \\
Group 6: & 20, 39, 47; \\
Group 7: & 23, 25, 29, 36, 42, 45; \\
Group 8: & 33. 
\end{tabular}
\end{center}
By adding the coefficients of the terms $T_{(i)}$ in the formula~\eqref{eq:Contraction3333} separately for each group listed above, we obtain the coefficients in~\eqref{eq:Contraction3x4}.
\end{proof}

In the same way one can obtain the following corollary.
\begin{corollary} \label{co:Contraction3x24}
 Let $A\in V_3$, $B\in V_4$. Then
\begin{multline*}
  \tau^5 (A\vee A\vee B)\\
=
  \frac{1}{105}\, \Big(
  3\, T_{(1)}
  +24\, T_{(2)}
  +12\, T_{(3)}
  +16\, T_{(4)}
  +2\, T_{(5)}
  +24\, T_{(6)}
  +24\, T_{(8)}
  \Big)  \,,
\end{multline*}
where $T_{(i)}$ are defined as in Corollary~\ref{co:Contraction334}. 
\end{corollary}

\chapter{Spectral Geometry of Submanifolds: Flat Manifold} \label{ch:FlatManifolds}
				
In this chapter we consider the plane, $M=\R^2$. 
We provide new proofs for asymptotic expansion Theorems~\ref{th:Expansion-plane-curve} and~\ref{th:Expansion-plane-point}. These proofs will give us a way to compute the asymptotic expansion coefficients. The coefficients $a_1$, $a_2$ in Lemma~\ref{le:Coefficients-plane-curve}, and the coefficients $A_1$, $A_2$ in Lemma~\ref{le:Coefficients-plane-point} are computed explicitly for the first time.

The injectivity radius of the plane is infinite, $r_{\rm inj}(\R^2)=\infty$, and the world function is defined globally, $\sigma:\R^2\times\R^2\to\R$.
We will use the Euclidean metric in $\R^2$ and local coordinates $x^\alpha$.
In Cartesian coordinates the world function is given by
\begin{gather*}
\sig{x,x'} = \dfrac{1}{2}\delta_{\alpha\beta}(x'^\alpha-x^\alpha)(x'^\beta-x^\beta) \,,
\end{gather*}
and the Laplacian by
\begin{gather*}
\Delta=\delta^{\alpha\beta}\pdf{}{x^\alpha}\pdf{}{x^\beta}\,.
\end{gather*}
The heat kernel $U(t;x,x')$ has the well-known form
\begin{gather*}
U(t;x,x')=(4\pi t)^{-1}\e{-\dfrac{\sig{x,x'}}{2t}}\,.
\end{gather*}

Let $\Phi:\R^2\to\R^2$ be a smooth mapping and let $\Sigma$, $U_\Phi(t)$, and $U_\Phi(t;x,x')$ be defined as in Section~\ref{se:Framework} by~\eqref{eq:Sigma}, \eqref{eq:Uphi_t}, \eqref{eq:Uphi_txxp}.
Let $S:\R^2\to\R$ be a real-valued function defined by
\begin{gather*} 
S(x)=\sig{\Phi(x),x}.
\end{gather*}
Since the world function $\sigma$ is defined globally, the function $S$ is well-defined. Obviously, $S(x)\geqslant 0$ for any $x\in\R^2$ and $S(x)=0$ on $\Sigma$.

Since $\dim\R^2=2$, there are only two possible cases: $\dim\Sigma=1$, that is, $\Sigma$ is a curve; and $\dim\Sigma=0$, that is, $\Sigma$ is a point. We will consider these two cases separately.

\section{One-Dimensional Submanifolds} \label{se:PlaneCurve}

Let $\dim\Sigma=1$. Let $\vect{h}$ be a smooth unit normal vector field along the curve $\Sigma$ and $\nabla_\vect{h}$ be the corresponding normal derivative. We define functions $s_k:\Sigma\to\R$ by 
\begin{gather} \label{eq:s_k}
s_k=\nabla_\vect{h}^k S \Big|_\Sigma \,,\  k\in\N\cup\set{0}.
\end{gather}
Since $S\Big|_\Sigma=0$ and $\nabla_\vect{h}S\Big|_\Sigma=0$, it follows that $s_0$ and $s_1$ vanish.

\begin{theorem} \label{th:Expansion-plane-curve}
Let $\Phi$ be a smooth mapping $\Phi:\R^2\to\R^2$.
Let $\Sigma$ be its fixed point set: $\Sigma=\set{x\in\R^2 : \Phi(x)=x}$.
Let $S$ be a real-valued function $S:\R^2\to\R$ defined by $S(x)=\sig{\Phi(x),x}$.
Let $\Phi$ satisfy the assumptions:
\begin{enumerate}
\item[\rm ($\Phi$.1)] $\Sigma$ is a smooth compact connected one-dimensional submanifold of $\,\R^2$ without boundary.
\item[\rm ($\Phi$.2)] There exist constants $p>0,\ R>0,\ C>0$, such that for any $x\in \R^2$,
if $\sig{x,0}>\dfrac{R^2}{2}$, then
\begin{gather} \label{eq:Phi.2}
S(x)\geqslant C\left[\sig{x,0}\right]^p .
\end{gather}
\item[\rm ($\Phi$.3)] $s_2=\nabla_\vect{h}^2 S\Big|_\Sigma>0$ at all points of $\Sigma$, where $s_k$ are defined as in~\eqref{eq:s_k}.
\end{enumerate}
Then there exists the asymptotic expansion as $t\to 0^+$
\begin{gather*}
\Theta_\Phi(t) \sim \sum_{k=0}^{\infty}t^{(2k-1)/2} A_k\,,
\end{gather*}
where the coefficients $A_{k}$ are locally computable in the form
\begin{gather*} 
A_k=\int_\Sigma \dvol(y)\, a_k(y),
\end{gather*}
and $a_k$ are scalars on $\Sigma$ depending polynomially on $s_2^{-1/2}$ and $s_k,\ k\geqslant 3$.
\end{theorem}

\begin{remark}
Assumptions {\rm ($\Phi$.1)--($\Phi$.3)} mean that the fixed point set $\Sigma$ is a smooth closed curve without self-intersections, the function $S(x)$ grows sufficiently fast at infinity and it is non-degenerate near $\Sigma$.
\end{remark}

\begin{proof}[Proof of Theorem~\ref{th:Expansion-plane-curve}]
Let $x^\alpha$ be Cartesian coordinates in $\R^2$. We have
\begin{gather*} 
  \Theta_\Phi(t)=(4\pi t)^{-1}\int_{\R^2} dx\, \e{-\frac{1}{2t}S(x)} \,,
\end{gather*}
where $dx=dx^1 dx^2$ is the standard Lebesgue measure.
The proof of the theorem will be based on a generalized Laplace method and involves the following steps:
\begin{enumerate}
\item Reducing the integration in the deformed heat trace $\Theta_\Phi(t)$ to an integral over a suitable neighborhood $M_\eps$ of $\Sigma$.
\item Choosing a special tangent--normal coordinate system in $M_\eps$.
\item Scaling the normal variable by $t^{1/2}$.
\item Expanding the integrand in the power series in $t^{1/2}$.
\item Computing integrals over the normal variable.
\end{enumerate}

\subsection{Tangent and Normal Coordinates}

Let the fixed point set $\Sigma$ be described by the equations
\begin{gather*}
\Sigma\ :\ x^\alpha=f^\alpha(y),
\end{gather*}
where $y$ is a parameter on $\Sigma$. 
The vector field $\vect{e}=\pdf{x}{y}$ is tangent to $\Sigma$. 
In a sufficiently small neighborhood of $\Sigma$ we can introduce local coordinates for $\R^2$ in the following way.
We construct a family of geodesics $x=x(y,z)$ in $\R^2$, such that at each point of $\Sigma$ the tangent vector to the geodesic is equal to $\vect{h}$, the smooth unit vector field, normal to $\Sigma$. We choose the parameter $z$ in a natural way, so that it is the signed distance to $\Sigma$ along the geodesic. That is, we have the differential equation
\begin{gather} \label{eq:GeodesicsDE}
\pdf[2]{x^\alpha}{z}=0
\end{gather}
with the initial conditions
\begin{align}
  \label{eq:GeodesicsIC1}
  x^\alpha\Big|_{z=0}&=f^\alpha(y) \,, \\[1ex]
  \label{eq:GeodesicsIC2}
  \pdf{x^\alpha}{z}\Big|_{z=0}&=h^\alpha(y) \,.
\end{align}

We restrict ourselves to a tubular neighborhood $M_\eps$
of $\Sigma$ with $\abs{z}<\eps$ for sufficiently small $\eps$, so that the change of coordinates $(x^1,x^2)\to(y,z)$ is not degenerate. We will specify the exact condition later.
The solution of the equations~\eqref{eq:GeodesicsDE}--\eqref{eq:GeodesicsIC2} is given by
\begin{gather*}
x^\alpha(y,z)=f^\alpha(y)+z h^\alpha(y).
\end{gather*}
At every point in $M_\eps$ we have:
\begin{gather*}
\pdf{x^\alpha}{y}=e^\alpha+z\pdf{h^\alpha}{y},\qquad
\pdf{x^\alpha}{z}=h^\alpha.
\end{gather*}
Due to the choice of $z$, we also have:
\begin{gather*}
\delta_{\alpha\beta}\pdf{x^\alpha}{z}\pdf{x^\beta}{z}=1,\qquad
\delta_{\alpha\beta}\pdf{x^\alpha}{y}\pdf{x^\beta}{z}=0,
\end{gather*}
i.e. the vectors $\pdf{x}{y}$ and $\pdf{x}{z}$ form an orthogonal system.

Finally, the volume form in the coordinates $(y,z)$ is
\begin{align*}
\dvol(y,z)&=\left[g(y,z)\right]^{1/2} dy\wedge dz=\eps_{\alpha\beta}\pdf{x^\alpha}{y}\pdf{x^\beta}{z} dy\wedge dz\\
&=\left[\gamma(y)\right]^{1/2} \, [1+z\, \kappa(y)] \, dy\wedge dz,
\end{align*}
where $\eps_{\alpha\beta}$ is the standard Levi-Civita symbol with $\eps_{12}=1$, $g=\det g_{\alpha\beta}$, $\gamma$ is the induced metric on $\Sigma$ and
\begin{gather*} 
\kappa=\left[\gamma(y)\right]^{-1/2} \, \eps_{\alpha\beta}\pdf{h^\alpha}{y} h^\beta\,.
\end{gather*}
Now we see, that it is enough to take $\eps$ less than $\ds \inf_{y\in\Sigma} \abs{\kappa(y)}^{-1}$. Since $\Sigma$ is compact, this infimum is positive.

\subsection{Reduction of the Area of Integration}

Since $S(x)$ is a non-negative function, the asymptotic expansion of the deformed heat trace $\Theta_\Phi(t)$ will be determined by a neighborhood of points where $S(x)=0$, i.e. by $M_\eps$:

\begin{lemma} \label{le:ReductionOfTheAreaOfIntegration}
There holds as $t\to 0^+$
\begin{gather*}
\Theta_\Phi(t)\sim(4\pi t)^{-1} \int_{M_\eps}dx\, \e{-\frac{1}{2t}S(x)}.
\end{gather*}
\end{lemma}
\begin{proof} Let $p,\ R,\ C$ be such constants that the estimation~\eqref{eq:Phi.2} holds.
Let $B_R$ be the ball with radius $R$ centered at the origin: $B_R=\set{x : \sig{x,0}\leqslant \dfrac{R^2}{2}}$.
We split $\Theta_\Phi(t)$ into the sum of three integrals:
\begin{gather*}
\Theta_\Phi(t)=I_1(t)+I_2(t)+I_3(t)\,,
\end{gather*}
where
\begin{align*}
I_1(t)=& (4\pi t)^{-1} \int_{M_\eps}dx\, \e{-\frac{1}{2t}S(x)},\\
I_2(t)=& (4\pi t)^{-1} \int_{B_R\setminus M_\eps}dx\, \e{-\frac{1}{2t}S(x)},\\
I_3(t)=& (4\pi t)^{-1} \int_{\R^2\setminus B_R}dx\, \e{-\frac{1}{2t}S(x)}.
\end{align*}
We will show that $I_2(t)$ and $I_3(t)$ are exponentially small as $t\to 0^+$, i.e. they are asymptotically equivalent to zero and the asymptotics of $\Theta_\Phi(t)$ is determined by $I_1(t)$ only.

Since $B_R\setminus M_\eps$ is compact and $S(x)>0$ on this set, it is separated from zero: $S(x)\geqslant\delta>0$ for all $x\in B_R\setminus M_\eps$. Hence, as $t\to 0^+$, we have
\begin{gather*}
I_2(t) \leqslant
(4\pi t)^{-1}\, \pi R^2 \e {-\frac{\delta}{2t}}\sim 0 \,.
\end{gather*}

For $I_3(t)$ we use~\eqref{eq:Phi.2} to estimate the integrand at infinity:
\begin{align*}
I_3(t)\leqslant & (4\pi t)^{-1} \int_{\R^2\setminus B_R} dx\, \e{-\frac{1}{2t}C \left[\sig{x,0}\right]^p}\\
= & (2t)^{-1} \int_{R}^{\infty}r dr \e{-\dfrac{C r^{2p}}{4t}}\\
= & (2t)^{-1} \dfrac{1}{2p} \left(\dfrac{4t}{C}\right)^{1/p} \Gamma\left(\frac{1}{p},\frac{CR^{2p}}{4t}\right),
\end{align*}
where $\Gamma(x,\alpha)$ is the incomplete $\Gamma$-function:
\begin{gather*}
\Gamma(x,\alpha)=\int_\alpha^\infty e^{-z} z^{x-1}\, dz.
\end{gather*}
It is well known (see, e.g.~\cite{Erdelyi}) that as $\alpha\to\infty$
\begin{gather} \label{eq:IncompleteGammaAsymptotics}
\Gamma(x,\alpha)=\alpha^{x-1} e^{-\alpha} \left(1+\bigO{\dfrac{1}{\alpha}}\right).
\end{gather}
Hence, we have as $t\to 0^+$
\begin{gather*}
I_3(t) \sim 0 \,.
\end{gather*}
Therefore, $\Theta_\Phi(t)\sim I_1(t)$.
\end{proof}

So, after the change of coordinates $x^\alpha=x^\alpha(y,z)$ in $M_\eps$ we get:
\begin{gather*}
\Theta_\Phi(t)\sim
(4\pi t)^{-1} \int_{\Sigma} \dvol(y)\, \int_{-\eps}^{\eps}dz\, \e{-\frac{1}{2t}S(x(y,z))} (1+z\, \kappa(y)) \,.
\end{gather*}

\subsection{Taylor Series Expansions}

By scaling the variable $z\to \sqrt{t}z$ we obtain
\begin{gather*} 
\Theta_\Phi(t)\sim
\dfrac{1}{4\pi \sqrt{t}} \int_{\Sigma}\dvol(y) \int_{-\eps/\sqrt{t}}^{\eps/\sqrt{t}} dz\,
\e{-\frac{1}{2t}S\left(x\left(y,\sqrt{t}z\right)\right)} \left(1+\sqrt{t}z\kappa(y)\right).
\end{gather*}
Now we expand the integrand in the last formula in the power series in $\sqrt{t}$. Below we compute the series expansion for $S\left(x\left(y,\sqrt{t}z\right)\right)$.

We define
\begin{gather*}
\Psi^\alpha(x)=\sigma^\alpha(\Phi(x),x)=\Phi^\alpha(x)-x^\alpha.
\end{gather*}
Then
\begin{gather*}
S\left(x\left(y,\sqrt{t}z\right)\right)=\dfrac{1}{2}\Psi^\alpha \left(x\left(y,\sqrt{t}z\right)\right) \Psi_\alpha \left(x\left(y,\sqrt{t}z\right)\right).
\end{gather*}
Since $z=0$ is the fixed point set of the mapping $\Phi$, we have
\begin{gather*}
\Psi^\alpha(f(y))=0 \,.
\end{gather*}
Therefore,
\begin{gather*}
\Psi^\alpha \left(x\left(y,\sqrt{t}z\right)\right)
\sim \sum_{k=1}^{\infty} \dfrac{t^{k/2}}{k!} \Psi^\alpha_{(k)}(y) \,z^k,
\end{gather*}
where
\begin{align*}
\Psi^\alpha_{(k)}&=\pdf[k]{}{z}\Psi^\alpha \Big|_{z=0}\\[1ex]
&=\dfrac{\pd^k \Psi^\alpha}{\pd x^{\nu_1}\dots\pd x^{\nu_k}} \Big|_{z=0} h^{\nu_1}\dots h^{\nu_k}.
\end{align*}
Now we obtain
\begin{gather} \label{eq:S-Expansion}
S\left(x\left(y,\sqrt{t}z\right)\right)
\sim\sum_{k=2}^{\infty} \dfrac{t^{k/2}}{k!} s_k(y)\, z^k,
\end{gather}
where
\begin{gather*}
s_k=\frac{1}{2} \sum_{l=1}^{k-1} \binom{k}{l} \Psi^\alpha_{(l)} \Psi_{\alpha(k-l)}.
\end{gather*}
We consider a non-degenerate case (see the assumption ($\Phi$.3)), when
\begin{gather*}
\Psi^\alpha_{(1)}=\pdf{\Psi^\alpha}{z} \Big|_{z=0} \neq 0
\end{gather*}
for any point of $\Sigma$. In this case
\begin{gather} \label{eq:s2-plane-curve}
s_2=\Psi^\alpha_{(1)} \Psi_{\alpha(1)}>0
\end{gather}
 everywhere on $\Sigma$. 

\begin{remark}
It is worth mentioning, that if $s_2=0$ somewhere on $\Sigma$, that is, the first non-vanishing term of the Taylor series~\eqref{eq:S-Expansion} of $S(x(y,z))$ is of degree higher than 2, then the asymptotics of $\Theta_\Phi(t)$ as $t\to 0^+$ will be substantially different.
Note also, that since $\Sigma$ is compact, this restriction means that ${s_2\geqslant\delta>0}$ for some constant $\delta$.
\end{remark}

The next several coefficients are:
\begin{align}
s_3&=3\, \Psi^\alpha_{(1)} \Psi_{\alpha(2)} \,,\\
s_4&=4\, \Psi^\alpha_{(1)} \Psi_{\alpha(3)} +3\, \Psi^\alpha_{(2)} \Psi_{\alpha(2)} \,,\\
s_5&=5\, \Psi^\alpha_{(1)} \Psi_{\alpha(4)} +10\, \Psi^\alpha_{(2)} \Psi_{\alpha(3)} \,,\\
s_6&=6\, \Psi^\alpha_{(1)} \Psi_{\alpha(5)} +15\, \Psi^\alpha_{(2)} \Psi_{\alpha(4)} +10\, \Psi^\alpha_{(3)} \Psi_{\alpha(3)} \,.  \label{eq:s6-plane-curve}
\end{align}

Now we rewrite the integral for $\Theta_\Phi(t)$ in the form
\begin{gather} \label{eq:3.9}
\Theta_\Phi(t) \sim \dfrac{1}{4\pi \sqrt{t}} \int_{\Sigma}\dvol(y) \int_{-\eps/\sqrt{t}}^{\eps/\sqrt{t}} dz\, \e{-\frac{s_2}{4}\, z^2} B(t,y,z),
\end{gather}
where
\begin{gather} \label{eq:B(tyz)}
B(t,y,z)=\e{-\frac{1}{2} \sum_{k=3}^{\infty} \frac{t^{k/2-1}}{k!} s_k\, z^k} \left( 1+t^{1/2}z\kappa\right).
\end{gather}
Our next step is to expand $B$ in the power series in $t^{1/2}$
\begin{gather} \label{eq:B-Expansion}
B \sim \sum_{k=0}^\infty b_k t^k + t^{1/2}\sum_{k=0}^\infty \tld{b}_k t^k
\end{gather}
and compute integrals over $z$.
It is easy to see that in this Taylor series:
\begin{enumerate}
\item All coefficients $b_k$ and $\tld{b}_k$ are polynomials in $z$, $\kappa$, and $s_k$, $k\geqslant 3$.
\item Integer powers of $t$ have coefficients $b_k$ with even powers of $z$ only.
\item Half-integer powers of $t$ have coefficients $\tld{b}_k$ with odd powers of $z$ only.
\end{enumerate}
All $\tld{b}_k$ will disappear after integration over $z$, since they are odd functions, so we do not need to keep track of them. For computing integrals involving $b_k$ we will use the following result.

\begin{lemma} \label{le:GIAsymptotics1d}
Let $a>0,\ k\in\Z,\ k\geqslant 0$. Then there holds as $t\to 0^+$
\begin{gather*}
\int_{-\eps/\sqrt{t}}^{\eps/\sqrt{t}} dz\, e^{-az^2}z^{2k} \sim
\int_{-\infty}^{\infty} dz\, e^{-az^2}z^{2k} \,.
\end{gather*}
\end{lemma}
\begin{proof}
We have as $t\to 0^+$
\begin{gather*}
\int_{-\infty}^{\infty} dz\, e^{-az^2}z^{2k} - \int_{-\eps/\sqrt{t}}^{\eps/\sqrt{t}} dz\, e^{-az^2}z^{2k} 
=a^{-(2k+1)/2} \Gamma\left(\dfrac{2k+1}{2},\dfrac{\eps^2 a}{t}\right) \sim 0 \,.
\end{gather*}
Here we used the asymptotic expansion~\eqref{eq:IncompleteGammaAsymptotics}.
\end{proof}

Therefore, by substituting the expansion~\eqref{eq:B-Expansion} into~\eqref{eq:3.9} and using Lemma~\ref{le:GIAsymptotics1d}, we obtain
\begin{gather*}
\Theta_\Phi(t) \sim \sum_{k=0}^{\infty} t^{(2k-1)/2} A_k \,,
\end{gather*}
where
\begin{gather*}
A_k=\int_{\Sigma}\dvol(y)\, a_k
\end{gather*}
and
\begin{gather} \label{eq:a_k-plane-curve}
a_k=\dfrac{1}{4\pi} \int_{-\infty}^\infty dz\, \e{-\frac{s_2}{4}\, z^2} b_k \,.
\end{gather}
Due to Lemma~\ref{le:GaussianIntegral1d}, integration over $z$ will bring inverse powers of $s_2$.
This completes the proof of Theorem~\ref{th:Expansion-plane-curve}.
\end{proof}

Now we will compute the first several coefficients $a_k$. By expanding the function $B(t,y,z)$, given by~\eqref{eq:B(tyz)}, in the power series in $t^{1/2}$, we obtain
\begin{align*}
b_0&=1 \,,\\
b_1&=\left( -\frac{1}{12} s_3 \kappa -\frac{1}{48} s_4 \right) z^4 +\frac{1}{288} s_3^2 z^6 \,,\\
b_2&=
  \left( -\frac{1}{240} s_5 \kappa -\frac{1}{1440} s_6 \right) z^6
  +\left( \frac{1}{576} s_3 s_4 \kappa +\frac{1}{2880} s_3 s_5 +\frac{1}{4608} s_4^2 \right) z^8\\
  &+\left( -\frac{1}{10368} s_3^3 \kappa -\frac{1}{13824} s_3^2 s_4 \right) z^{10}
  +\frac{1}{497664} s_3^4 z^{12} \,.
\end{align*}
Therefore, by using formula~\eqref{eq:a_k-plane-curve}, Lemma~\ref{le:GaussianIntegral1d}, and performing straightforward computations, we get the following lemma.

\begin{lemma} \label{le:Coefficients-plane-curve}
Let conditions of Theorem~\ref{th:Expansion-plane-curve} be satisfied. Then
\begin{align*}
a_0&=(4\pi s_2)^{-1/2} \,,\\
a_1&=(4\pi s_2)^{-1/2} \left[
    \left( -s_3 \kappa - \frac{1}{4} s_4 \right) s_2^{-2} + \frac{5}{12} s_3^2 s_2^{-3} \right] ,\\
a_2&=(4\pi s_2)^{-1/2} \left[
  \left( -\frac{1}{2} s_5 \kappa - \frac{1}{12} s_6 \right) s_2^{-3} \right.
  +\left( \frac{35}{12} s_3 s_4 \kappa + \frac{7}{12} s_3 s_5 + \frac{35}{96} s_4^2 \right) s_2^{-4}\\
 &+\left.\left( -\frac{35}{12} s_3^3 \kappa -\frac{35}{16} s_3^2 s_4 \right) s_2^{-5}
  +\frac{385}{288} s_3^4 s_2^{-6} \right],
\end{align*}
where coefficients $s_k$ are given by~\eqref{eq:s2-plane-curve}--\eqref{eq:s6-plane-curve}.
\end{lemma}

\section{Zero-Dimensional Submanifolds} \label{se:PlanePoint}

In this section we consider the case when the fixed point set $\Sigma$ is just a point, that is, $\dim\Sigma=0$.

\begin{theorem} \label{th:Expansion-plane-point}
Let $\Phi$ be a smooth mapping $\Phi:\R^2\to\R^2$.
Let $\Sigma$ be its fixed point set: $\Sigma=\set{x\in\R^2 : \Phi(x)=x}$.
Let $S$ be a real-valued function $S:\R^2\to\R$ defined by $S(x)=\sig{\Phi(x),x}$.
Let $\Phi$ satisfy the assumptions:
\begin{enumerate}
\item[\rm ($\Phi$.1)] $\Sigma$ is a one-point set, $\Sigma=\set{\xi_0}$.
\item[\rm ($\Phi$.2)] There exist constants $p>0,\ R>0,\ C>0$ such that for any $x\in \R^2$,
if $\sig{x,0}>\dfrac{R^2}{2}$, then
\begin{gather*} 
S(x)\geqslant C\left[\sig{x,0}\right]^p .
\end{gather*}
\item[\rm ($\Phi$.3)] $\xi_0$ is a non-degenerate critical point of the function $S$.
\end{enumerate}
Then there exists the asymptotic expansion as $t\to 0^+$
\begin{gather} \label{eq:Expansion-plane-point}
\Theta_\Phi(t) \sim \sum_{k=0}^{\infty}t^k A_k\,,
\end{gather}
where scalars $A_k$ depend polynomially on derivatives of the function $S$ and the inverse of its Hessian matrix, evaluated at $\xi_0$.
\end{theorem}

Before the proof of the theorem, let us note that, given the fact that coefficients $A_k$ are polynomials in derivatives of the function $S$, it is possible to write down their general form by using dimensional considerations~\cite{Gilkey}. The deformed heat trace $\Theta_\Phi(t)$ is dimensionless, so the expansion on the right hand side of~\eqref{eq:Expansion-plane-point} must also be dimensionless. From the heat equation, $t$ has the dimension $L^2$, so the dimension of $A_k$ is $L^{-2k}$. The function $S$ has the dimension of $L^2$ ($L$~--- dimension of length), since it is equal to the world function. Its $m$-th derivative with respect to the space variables has the dimension $L^{2-m}$. In particular, its second derivative is dimensionless, and all derivatives of order higher than two have dimensions of negative powers of $L$.

In order to construct scalars $A_k$, it is necessary to contract all indices of derivatives with indices of the inverse of the Hessian matrix of the function $S$, which we will denote by $Q$. The order of this contraction can be different, therefore, the coefficients $A_k$ are linear combinations of invariants of the form
\begin{gather*} 
T(\mi{m},\phi)=
\big(\underbrace{Q\otimes\ldots\otimes Q}_{M}\big)^{\nu_1\ldots\nu_{2M}}
\big(\nabla^{m_1}S\otimes\ldots\otimes\nabla^{m_l}S\big)_{\nu_{\phi(1)}\ldots\nu_{\phi(2M)}} \,,
\end{gather*}
where $\mi{m}=(m_1,\ldots,m_l)$ is a multi-index, $\abs{\mi{m}}=m_1+\ldots+m_l=2M$, and $\phi\in S_{2M}$ is a permutation of numbers $1,\ldots,2M$. Note, that $m_i\geqslant 3$ for any $i$, since $A_k$ depends polynomially only on derivatives of the function $S$ of order higher than two (besides, the function $S$ itself and its first derivative vanish on $\Sigma$).

Let us introduce the length, the order, and the dimension of the term $T$ (we omit arguments):
\begin{align*}
\len T&=l \,, \\
\ord T&=2M \,, \\
\dim T&=2M-2l \,.
\end{align*}
In other words, the length is the number of derivatives of the function $S$ in the term $T$, the order is the sum of orders of all derivatives, and the dimension is the (negative) space-dimension of $T$. Clearly,
\begin{gather*}
\dim T=\ord T-2\len T \,.
\end{gather*}
Now we are ready to write down the general formula for the coefficients $A_k$:
\begin{gather*}
A_{k}=
  \sum_{l=1}^{2k}\
  \sum_{\clearfrac{\mi{m}=(m_1,\ldots,m_l)}{\abs{\mi{m}}=2(k+l)}}\
  \sum_{\phi\in S_{2(k+l)}}
   c(\mi{m},\phi)\, T(\mi{m},\phi) \,,
\end{gather*}
where $c(\mi{m},\phi)$ are some numerical coefficients. Note, that this sum is finite for any $k$.

\begin{proof}[Proof of Theorem~\ref{th:Expansion-plane-point}]
The idea of this proof is very similar to the one of Theorem~\ref{th:Expansion-plane-curve}, but in this case we need to use the multi-dimensional version of the Laplace method.
In the same way, due to Lemma~\ref{le:ReductionOfTheAreaOfIntegration}, we can reduce the area of integration to a small neighborhood $M_\eps$ of $\Sigma$ and introduce there a system of Cartesian coordinates $z^\alpha$ with the origin at $\xi_0$. Then we expand the function $S$ in the Taylor series in $z$, scale the variables $z^\alpha\to\sqrt{t}z^\alpha$ and obtain
\begin{gather} \label{eq:3.14}
\Theta_\Phi(t) \sim (4\pi)^{-1} \int_{\R^2} d^2 z\, \e{-\frac{1}{4}\, S_{\alpha\beta} z^\alpha z^\beta} B(t,z),
\end{gather}
where
\begin{align}
\label{eq:B(tz)}
B(t,z)&=\e{-\frac{1}{2} \sum_{k=3}^\infty \dfrac{t^{k/2-1}}{k!} S_{\nu_1\ldots\nu_k} z^{\nu_1}\ldots z^{\nu_k}},\\
\nonumber
S_{\nu_1\ldots\nu_k}&=\pdf{}{z^{\nu_1}}\cdots\pdf{}{z^{\nu_k}} S \Big|_{z=0}\,.
\end{align}
Again, $B$ has the expansion
\begin{gather} \label{eq:B_series-plane-point}
B \sim \sum_{k=0}^\infty b_k t^k + t^{1/2}\sum_{k=0}^\infty \tld{b}_k t^k
\end{gather}
with coefficients $b_k$ and $\tld{b}_k$ being polynomials in $z^\alpha$ and derivatives of $S$ of order higher than two. And there is no need to keep track of $\tld{b}_{k}$ since they do not contribute to the asymptotic expansion of the deformed heat trace $\Theta_\Phi(t)$.

Since $\xi_0$ is a non-degenerate critical point of the function $S$, its Hessian matrix $S_{\alpha\beta}$ is invertible. We will denote the inverse matrix by $Q$. It satisfies the equation
\begin{gather*}
Q^{\mu\nu}S_{\nu\alpha}=\delta^{\mu}_{\alpha}\,.
\end{gather*}
By substituting~\eqref{eq:B_series-plane-point} into the integral~\eqref{eq:3.14} and applying Lemma~\ref{le:GaussianIntegralmd}, we conclude that
\begin{gather*}
\Theta_\Phi(t) \sim \sum_{k=0}^{\infty}t^k A_k \,,
\end{gather*}
where
\begin{gather} \label{eq:A_k-plane-point}
A_k=\dfrac{1}{4\pi} \int_{\R^2} d^2 z\, \e{-\frac{1}{4}\, S_{\alpha\beta} z^\alpha z^\beta} b_k \,.
\end{gather}
Due to Lemma~\ref{le:GaussianIntegralmd}, integration over $z$ will bring factors of $Q$.
This completes the proof of Theorem~\ref{th:Expansion-plane-point}.
\end{proof}

Now we will compute the first several coefficients $A_k$. In order to simplify formulas, we will use the following compact notation for the (symmetrized) derivatives of the function $S$:
\begin{align*}
S_{(k)}&=S_{(\nu_1\ldots\nu_k)} \,,\\
\tld{S}_{(k)}&=\dfrac{1}{k!} \, S_{\nu_1\ldots\nu_k} z^{\nu_1}\ldots z^{\nu_k} \,.
\end{align*}
By expanding the function $B(t,z)$, given by~\eqref{eq:B(tz)}, in the power series in $t^{1/2}$, we obtain
\begin{align*}
b_0&=1 \,,\\
b_1&=
  -\frac{1}{2}\,\tld{S}_{(4)} + \frac{1}{8}\,\tld{S}_{(3)}\tld{S}_{(3)} \,,\\
b_2&=
  -\frac{1}{2}\tld{S}_{(6)}
  +\frac{1}{4}\,\tld{S}_{(5)}\tld{S}_{(3)} +\frac{1}{8}\,\tld{S}_{(4)}\tld{S}_{(4)}
  -\frac{1}{16}\,\tld{S}_{(4)}\tld{S}_{(3)}\tld{S}_{(3)}
  +\frac{1}{384}\,\tld{S}_{(3)}\tld{S}_{(3)}\tld{S}_{(3)}\tld{S}_{(3)}\,.
\end{align*}
In terms of the mapping $\Phi$ and its derivatives, derivatives of the functions $S$ are
\begin{align}
S_{(\alpha\beta)} &=
g_{\alpha\beta}
-2\, g_{\mu(\alpha}\Phi^\mu{}_{\beta)}
+g_{\mu\nu} \Phi^\mu{}_{(\alpha} \Psi^\nu{}_{\beta)}
\,, \label{eq:S2-plane-point} \displaybreak[0]\\
S_{(\alpha\beta\gamma)} &=
-3\, g_{\mu(\alpha} {\Phi^{\mu}}_{\beta\gamma)}
+3\, g_{\mu\nu} {\Phi^{\mu}}_{(\alpha}{\Phi^{\nu}}_{\beta\gamma)}
\,, \displaybreak[0]\\
S_{(\alpha\beta\gamma\delta)} &=
-4\, g_{\mu(\alpha} {\Phi^{\mu}}_{\beta\gamma\delta)}
+4\, g_{\mu\nu} {\Phi^{\mu}}_{(\alpha}{\Phi^{\nu}}_{\beta\gamma\delta)}
+3\, g_{\mu\nu} {\Phi^{\mu}}_{(\alpha\beta}{\Phi^{\nu}}_{\gamma\delta)}
\,, \displaybreak[0]\\
S_{(\alpha\beta\gamma\delta\eps)} &=
-5\, g_{\mu(\alpha} {\Phi^{\mu}}_{\beta\gamma\delta\eps)}
+5\, g_{\mu\nu} {\Phi^{\mu}}_{(\alpha}{\Phi^{\nu}}_{\beta\gamma\delta\eps)}
+10\, g_{\mu\nu} {\Phi^{\mu}}_{(\alpha\beta}{\Phi^{\nu}}_{\gamma\delta\eps)}
\,,\\
S_{(\alpha\beta\gamma\delta\eps\zeta)} &=
-6\, g_{\mu(\alpha} {\Phi^{\mu}}_{\beta\gamma\delta\eps\zeta)}
+6\, g_{\mu\nu} {\Phi^{\mu}}_{(\alpha}{\Phi^{\nu}}_{\beta\gamma\delta\eps\zeta)}
+15\, g_{\mu\nu} {\Phi^{\mu}}_{(\alpha\beta}{\Phi^{\nu}}_{\gamma\delta\eps\zeta)}
\nonumber \displaybreak[0]\\ &
+10\, g_{\mu\nu} {\Phi^{\mu}}_{(\alpha\beta\gamma}{\Phi^{\nu}}_{\delta\eps\zeta)}
\,. \label{eq:S6-plane-point}
\end{align}
These formulas can be obtained as a particular case of~\eqref{eq:S2curvature}--\eqref{eq:S6curvature}, where all derivatives of the world function vanish except for the second order.

Let $\tau:V_{m+2}\to V_{m}$ be the operator of contraction of a symmetric tensor with the matrix $Q$:
\begin{gather*}
(\tau A)_{\nu_1\dots\nu_{m}}=Q^{\alpha\beta} A_{\alpha\beta\nu_1\dots\nu_{m}}\,.
\end{gather*}
By using formula~\eqref{eq:A_k-plane-point}, Lemma~\ref{le:GaussianIntegralmd}, and performing straightforward computations, we get
\begin{align*}
A_0&=\sqrt{\det Q} \,,\\[1ex]
A_1&=\sqrt{\det Q} \left[
  -\frac{1}{4}\, \tau^2 \left( S_{(4)} \right)
  +\frac{5}{12}\, \tau^3 \left( S_{(3)}\vee S_{(3)} \right)
  \right],\\[1ex]
A_2&=\sqrt{\det Q} \left[
  -\frac{1}{12}\, \right. \tau^3 \left( S_{(6)} \right)
  +\frac{7}{12}\, \tau^4 \left( S_{(3)}\vee S_{(5)} \right)
  +\frac{35}{96}\, \tau^4 \left( S_{(4)}\vee S_{(4)} \right) \\[1ex]
  &-\frac{35}{16}\, \tau^5 \left( S_{(3)}\vee S_{(3)}\vee S_{(4)} \right)
  +\left. \frac{385}{288}\, \tau^6 \left( S_{(3)}\vee S_{(3)}\vee S_{(3)}\vee S_{(3)} \right)
  \right].
\end{align*}
Since expressions involving symmetrization are inconvenient for practical computations, we will expand them in a symmetrization-free form. We will use the following notation 
\begin{gather}
S_{[2k]\nu_1\ldots\nu_m}=\left(\tau^k\left(S_{(2k+m)}\right)\right)_{\nu_1\ldots\nu_m}
=Q^{\mu_1\mu_2}\ldots Q^{\mu_{2k-1}\mu_{2k}}S_{\mu_1\mu_2\ldots\mu_{2k-1}\mu_{2k}\nu_1\ldots\nu_m} \,,
\label{eq:contractedS}
\end{gather}
that is, the number in brackets is the number of indices contracted with $Q$. By using corollaries~\ref{co:Contraction33}, \ref{co:Contraction35}, \ref{co:Contraction44}, \ref{co:Contraction3x24}, and \ref{co:Contraction3x4}, we obtain the following lemma.

\begin{lemma} \label{le:Coefficients-plane-point}
Let conditions of Theorem~\ref{th:Expansion-plane-point} be satisfied. Then
\begin{align*}
A_0&=\sqrt{\det Q} \,,
    \displaybreak[0]\\[1ex]
A_1&=\sqrt{\det Q} \left[ 
  -\frac{1}{4} S_{[4]}
  +\frac{1}{4} S_{[2]\nu_1} S_{[2]\nu_2} Q^{\nu_1\nu_2}
  +\frac{1}{6} S_{\nu_1\nu_2\nu_3} S_{\nu_4\nu_5\nu_6}  Q^{\nu_1\nu_4} Q^{\nu_2\nu_5} Q^{\nu_3\nu_6}
  \right],
    \displaybreak[0]\\[1ex]
A_2&=\sqrt{\det Q} \left[
  -\frac{1}{12} \right. S_{[6]}
  +\frac{1}{4} S_{[2]\nu_1} S_{[4]\nu_2} Q^{\nu_1\nu_2}
  +\frac{1}{3} S_{\nu_1\nu_2\nu_3} S_{[2]\nu_4\nu_5\nu_6} Q^{\nu_1\nu_4} Q^{\nu_2\nu_5} Q^{\nu_3\nu_6}
    \displaybreak[0]\\[1ex]
  &+\frac{1}{32} S_{[4]} S_{[4]}
  +\frac{1}{4} S_{[2]\nu_1\nu_2} S_{[2]\nu_3\nu_4} Q^{\nu_1\nu_3} Q^{\nu_2\nu_4}
    \displaybreak[0]\\[1ex]
  &+\frac{1}{12} S_{\nu_1\nu_2\nu_3\nu_4} S_{\nu_5\nu_6\nu_7\nu_8}
    Q^{\nu_1\nu_5} Q^{\nu_2\nu_6} Q^{\nu_3\nu_7} Q^{\nu_4\nu_8}
  -\frac{1}{16} S_{[2]\nu_1} S_{[2]\nu_2} S_{[4]} Q^{\nu_1\nu_2}
    \displaybreak[0]\\[1ex]
  &-\frac{1}{2} S_{[2]\nu_1} S_{\nu_2\nu_3\nu_4} S_{[2]\nu_5\nu_6}
    Q^{\nu_1\nu_2} Q^{\nu_3\nu_5} Q^{\nu_4\nu_6}
  -\frac{1}{4} S_{[2]\nu_1} S_{[2]\nu_2} S_{[2]\nu_3\nu_4}  Q^{\nu_1\nu_3} Q^{\nu_2\nu_4}
    \displaybreak[0]\\[1ex]
  &-\frac{1}{3} S_{[2]\nu_1} S_{\nu_2\nu_3\nu_4} S_{\nu_5\nu_6\nu_7\nu_8}
    Q^{\nu_1\nu_5} Q^{\nu_2\nu_6} Q^{\nu_3\nu_7} Q^{\nu_4\nu_8}
    \displaybreak[0]\\[1ex]
  &-\frac{1}{24} S_{\nu_1\nu_2\nu_3} S_{\nu_4\nu_5\nu_6} S_{[4]}
    Q^{\nu_1\nu_4} Q^{\nu_2\nu_5} Q^{\nu_3\nu_6}
    \displaybreak[0]\\[1ex]
  &-\frac{1}{2} S_{\nu_1\nu_2\nu_3} S_{\nu_4\nu_5\nu_6} S_{[2]\nu_7\nu_8}
    Q^{\nu_1\nu_4} Q^{\nu_2\nu_5} Q^{\nu_3\nu_7} Q^{\nu_6\nu_8}
    \displaybreak[0]\\[1ex]
  &-\frac{1}{2} S_{\nu_1\nu_2\nu_3} S_{\nu_4\nu_5\nu_6} S_{\nu_7\nu_8\nu_9\nu_{10}}
    Q^{\nu_1\nu_4} Q^{\nu_2\nu_7} Q^{\nu_3\nu_8} Q^{\nu_5\nu_9} Q^{\nu_6\nu_{10}}
    \displaybreak[0]\\[1ex]
  &+\frac{1}{32} S_{[2]\nu_1} S_{[2]\nu_2} S_{[2]\nu_3} S_{[2]\nu_4}
    Q^{\nu_1\nu_2} Q^{\nu_3\nu_4}
    \displaybreak[0]\\[1ex]
  &+\frac{1}{24} S_{[2]\nu_1} S_{[2]\nu_2} S_{\nu_3\nu_4\nu_5} S_{\nu_6\nu_7\nu_8}
    Q^{\nu_1\nu_2} Q^{\nu_3\nu_6} Q^{\nu_4\nu_7} Q^{\nu_5\nu_8}
    \displaybreak[0]\\[1ex]
  &+\frac{1}{4} S_{[2]\nu_1} S_{\nu_2\nu_3\nu_4} S_{\nu_5\nu_6\nu_7} S_{[2]\nu_8}
    Q^{\nu_1\nu_2} Q^{\nu_3\nu_5} Q^{\nu_4\nu_6} Q^{\nu_7\nu_8}
    \displaybreak[0]\\[1ex]
  &+\frac{1}{12} S_{[2]\nu_1} S_{\nu_2\nu_3\nu_4} S_{[2]\nu_5} S_{[2]\nu_6}
    Q^{\nu_1\nu_2} Q^{\nu_3\nu_5} Q^{\nu_4\nu_6}
    \displaybreak[0]\\[1ex]
  &+\frac{1}{2} S_{[2]\nu_1} S_{\nu_2\nu_3\nu_4} S_{\nu_5\nu_6\nu_7} S_{\nu_8\nu_9\nu_{10}}
    Q^{\nu_1\nu_2} Q^{\nu_3\nu_5} Q^{\nu_4\nu_8} Q^{\nu_6\nu_9} Q^{\nu_7\nu_{10}}
    \displaybreak[0]\\[1ex]
  &+\frac{1}{72} S_{\nu_1\nu_2\nu_3} S_{\nu_4\nu_5\nu_6} S_{\nu_7\nu_8\nu_9} S_{\nu_{10}\nu_{11}\nu_{12}}
    Q^{\nu_1\nu_4} Q^{\nu_2\nu_5} Q^{\nu_3\nu_6} Q^{\nu_7\nu_{10}} Q^{\nu_8\nu_{11}} Q^{\nu_9\nu_{12}}
    \displaybreak[0]\\[1ex]
  &+\frac{1}{4} S_{\nu_1\nu_2\nu_3} S_{\nu_4\nu_5\nu_6} S_{\nu_7\nu_8\nu_9} S_{\nu_{10}\nu_{11}\nu_{12}}
    Q^{\nu_1\nu_4} Q^{\nu_2\nu_5} Q^{\nu_3\nu_7} Q^{\nu_6\nu_{10}} Q^{\nu_8\nu_{11}} Q^{\nu_9\nu_{12}}
    \displaybreak[0]\\[1ex]
  &\left.
  +\frac{1}{6} S_{\nu_1\nu_2\nu_3} S_{\nu_4\nu_5\nu_6} S_{\nu_7\nu_8\nu_9} S_{\nu_{10}\nu_{11}\nu_{12}}
    Q^{\nu_1\nu_4} Q^{\nu_2\nu_7} Q^{\nu_3\nu_{10}} Q^{\nu_5\nu_8} Q^{\nu_6\nu_{11}} Q^{\nu_9\nu_{12}}
  \right],
\end{align*}
where tensors $S_{[2k]\nu_1\ldots\nu_m}$ are given by~\eqref{eq:contractedS}, \eqref{eq:S2-plane-point}--\eqref{eq:S6-plane-point}, and $Q$ is the inverse matrix to $S_{\alpha\beta}$.
\end{lemma}

\chapter{Spectral Geometry of Submanifolds: Curved Manifold} \label{ch:CurvedManifolds}

In this chapter we will consider a curved manifold $M$, and a smooth mapping $\Phi:M\to M$ on it, such that its fixed point set $\Sigma$ is zero-dimensional, that is, it is just a single point.
We provide a new proof for asymptotic expansion Theorem~\ref{th:Expansion-curved-point}. This proof will give us a way to compute the asymptotic expansion coefficients. The coefficients $A_1$, $A_2$ in Lemma~\ref{le:Coefficients-curved-point}, and the coefficient $A_2$ in Lemma~\ref{le:Coefficients-curved-point-isometry} are computed explicitly for the first time.

\section{Main Theorem} \label{se:CurvedPoint}

\begin{theorem} \label{th:Expansion-curved-point}
Let $(M,g)$ be an $n$-dimensional compact smooth Riemannian manifold without boundary.
Let $\Phi$ be a smooth mapping $\Phi:M\to M$.
Let $\Sigma$ be its fixed point set: $\Sigma=\set{x\in M : \Phi(x)=x}$.
Let $\Phi$ satisfy the assumptions:
\begin{enumerate}
\item[\rm ($\Phi$.1)] $\Sigma$ is a one-point set, $\Sigma=\set{x_0}$.
\item[\rm ($\Phi$.2)] $\det(I-d\Phi)(x_0)\neq 0$.
\end{enumerate}
Then there exists the asymptotic expansion as $t\to 0^+$
\begin{gather} \label{eq:Expansion-curved-point}
\Theta_\Phi(t) \sim \sum_{k=0}^{\infty}t^k A_k\,,
\end{gather}
where the coefficients $A_k$ are scalar invariants depending polynomially on covariant derivatives of the curvature of the metric $g$, symmetrized covariant derivatives of the differential $d\Phi$ of the mapping $\Phi$, and the matrix
\begin{gather} \label{eq:Q-curved-point}
Q= \left\{ (I-d\Phi)^T g (I-d\Phi) \right\}^{-1} ,
\end{gather}
evaluated at $x_0$.
\end{theorem}

\begin{proof}
The idea of this proof is the same as for a flat manifold. We will choose a convenient system of coordinates, reduce the integral over the whole manifold $M$ to the integral over a neighborhood of the fixed point set $\Sigma$, expand integrand in Taylor series at $x_0$, and compute the asymptotics of the integral by reducing it to the standard Gaussian integrals, which are known exactly.

In this proof we will make several references to some computations in the upcoming sections. However, it does not affect rigor, since those computations do not use this theorem.

Let $S$ be a real-valued function, defined in some neighborhood of $\Sigma$ by $S(x)=\sig{x,\Phi(x)}$. Obviously, $x_0$ is a point of minimum of the function $S$, so
$S(x_0)=0$ and $S_{;\alpha}(x_0)=0$. We will show below (Section~\ref{se:SymDerOfS},~\eqref{eq:S_alphabeta}) that
\begin{gather*}
S_{;\alpha\beta} (x_0)=
	\big\{
	\delta^\mu{}_\alpha
	-(d\Phi(x_0))^\mu{}_\alpha
	\big\}
 	g_{\mu\nu}(x_0)
	\big\{
	\delta^\nu{}_\beta
	-(d\Phi(x_0))^\nu{}_\beta
	\big\}
	\,,
\end{gather*}
or, in the matrix notation,
\begin{gather} \label{eq:S2Matrix}
\nabla\nabla S (x_0)=\big\{I-d\Phi(x_0)\big\}^T g(x_0)\big\{I-d\Phi(x_0)\big\} \,.
\end{gather}
Therefore, the matrix $Q$ defined by~\eqref{eq:Q-curved-point} is nothing but $\big\{\nabla\nabla S (x_0)\big\}^{-1}$.

From now on for simplicity we will omit a semicolon in covariant derivatives of scalar functions, for example, $S_\alpha=S_{;\alpha}$, $S_{\alpha\beta}=S_{;\alpha\beta}$, etc.

Let $M_\eps$ be an open geodesic ball with center at the point $x_0$ and small enough radius $\eps$, so that the functions $\sigma$ and $S$ are well-defined in $M_\eps$. Since $M\setminus M_\eps$ is compact, there exists a positive constant $\delta$, such that for any point $x\in M\setminus M_\eps$ we have
\begin{gather*}
\dist(x,\Phi(x))>\delta \,.
\end{gather*}
From the last inequality and the standard off-diagonal estimates for the heat kernel in~\cite{Grigoryan} it follows that as $t\to 0^+$
\begin{gather*}
\int_{M\setminus M_\eps} \dvol(x)\, U(t;x,\Phi(x)) \sim 0 \,.
\end{gather*}
Therefore,
\begin{gather*}
\Theta_\Phi(t) \sim \int_{M_\eps} \dvol(x)\, U(t;x,\Phi(x))
\end{gather*}
and according to~\eqref{eq:HeatKernel} we have
\begin{gather} \label{eq:curvI}
\Theta_\Phi(t) \sim \int_{M_\eps}
	\dvol(x) (4\pi t)^{-n/2} \varDelta^{1/2}(x,\Phi(x)) \e{-\dfrac{S(x)}{2t}} \Omega(t;x,\Phi(x))
	\,.
\end{gather}

Now we will choose a convenient coordinate system in $M_\eps$ and expand the integrand of~\eqref{eq:curvI} in the covariant Taylor series at $x_0$.

\subsection{Normal Coordinates}

Let $x$ be a point in $M_\eps$.
Due to~\eqref{eq:TanToGeod}, the components of the tangent vector $\vect{Z}=z^\alpha \pdf{}{x^\alpha}$ to the geodesic connecting the points $x_0$ and $x$ at the point $x_0$ in the direction of the point $x$ and with the norm equal to the length of this geodesic are given by
\begin{gather*}
z^\alpha(x)=-g^{\alpha\beta}(x_0) \pdf{\sigma}{x_0^\beta}(x_0,x) \,.
\end{gather*}
According to~\eqref{eq:ExpOfSigmaAlpha}, we have
\begin{gather*}
x=\exp_{x_0}(\vect{Z}) \,.
\end{gather*}
One can use this vector to define new coordinates in $M_\eps$ by assigning to a point $x$ with  coordinates $x^\alpha$ new coordinates $z^\alpha$ equal to the components of the vector $\vect{Z}$ in the coordinates $x^\alpha$.
Let $\phi:M_\eps\to\R^n$ be such a coordinate mapping.
Since $M_\eps$ is a geodesic ball with radius $\eps$, the image of $M_\eps$ is
\begin{gather*}
\phi(M_\eps)=\set{z\in\R^n : \norm{z}<\eps} \,.
\end{gather*}

The differentials are related by
\begin{gather*}
dz^\alpha =
	\pdf{z^\alpha}{x^\beta} dx^\beta
	=-g^{\alpha\gamma}(x_0) \frac{\pd^2 \sigma}{\pd x^\beta \pd x_0^\gamma}(x_0,x) dx^\beta
	\,.
\end{gather*}
Thus, for the measures we obtain
\begin{align*}
dz &=\abs{g(x_0)}^{-1} \det\left(-\frac{\pd^2 \sigma}{\pd x^\beta \pd x_0^\gamma}(x_0,x) \right) dx
	\\
	&=\abs{g(x_0)}^{-1/2} \varDelta(x_0,x)\, \abs{g(x)}^{1/2} dx
	\\
	&=\abs{g(x_0)}^{-1/2} \varDelta(x_0,x)\, \dvol(x)
	\,.
\end{align*}
So, the Riemannian volume element has in new coordinates the form
\begin{gather*}
\dvol(x)=\varDelta^{-1}(x_0,x)\, \abs{g(x_0)}^{1/2} dz \,.
\end{gather*}

After the change of variables $x\to z$ integral \eqref{eq:curvI} takes the form
\begin{multline} \label{eq:curvIM'}
\Theta_\Phi(t) \sim
	(4\pi t)^{-n/2}
	\int_{\norm{z}<\eps}  dz\,
	\abs{g(x_0)}^{1/2}
	\e{-\dfrac{S(x(z))}{2t}}
	\\ \times
	\varDelta^{-1}\big(x_0,x(z)\big)
	\varDelta^{1/2}\big(x(z),\Phi(x(z))\big)
	\Omega\big(t;x(z),\Phi(x(z))\big)
	\,.
\end{multline}

\subsection{Taylor Series Expansion}

Let $f$ be a smooth function on the manifold $M$. Then in the neighborhood $M_\eps$ of the fixed point $x_0$ this function can be expanded in covariant Taylor series~\cite{Avramidi2000}
\begin{gather} \label{eq:TaylorSeries}
f(x(z))=\sum_{k=0}^\infty \frac{1}{k!}\,
	\big(\nabla^k f\big)_{(\mu_1\ldots\mu_k)}(x_0)\, z^{\mu_1}\ldots z^{\mu_k}
	\,,
\end{gather}
where $z^\mu$ are the coordinates introduced in the previous section. Since this expansion involves only covariant derivatives of the function $f$ and the coordinates $z^\mu$ are covariant derivatives of the world function, all parts of this expression are written in an invariant form. Note also, that all covariant derivatives of the function $f$ are symmetrized, which makes~\eqref{eq:TaylorSeries} very convenient for practical computations. We will use the notation
\begin{align*}
f_{(k)} &= \Sym\left(\nabla^k f\right) (x_0)
	\,,\\ 
\tld{f}_{(k)} &= \frac{1}{k!} \left(f_{(k)}\right)_{\mu_1\ldots\mu_k} z^{\mu_1}\ldots z^{\mu_k}
	\,.
\end{align*}

We will expand each factor of the integrand in~\eqref{eq:curvIM'} in the Taylor series at $x_0$ separately and then multiply the expansions. Note, that we will need to compute symmetrized covariant derivatives of functions of the form $f(x,\Phi(x))$ at the point $x_0$, which has the property $\Phi(x_0)=x_0$.

The function $\Omega$ has asymptotic expansion~\eqref{eq:OmegaAE} in powers of $t$. Let
\begin{gather*}
a_m(x(z))=\omega_m\big(x(z),\Phi(x(z))\big)=\sum_{k=0}^\infty \tld{a}_{m(k)} \,,
\end{gather*}
that is, the first index stands for the coefficient of the power series in $t$, and the second one for the coefficient of the covariant Taylor series in $z$.
It is convenient to represent the Van Fleck-Morette determinant in the form~\cite{Avramidi2000}
\begin{gather} \label{eq:VFMzeta}
\varDelta=\exp (2\zeta) \,.
\end{gather}
The function $\zeta$ has a known expansion in a Taylor series with the first two coefficients being equal to zero~\cite{Avramidi2000},
\begin{gather*}
\zeta(x_0,x(z))=\sum_{k=2} \tld{\zeta}_{(k)} \,.
\end{gather*}
Let
\begin{gather*}
\eta(x(z))=\zeta\big(x(z),\Phi(x(z))\big)
\end{gather*}
and
\begin{gather*}
\eta(x(z))=\sum_{k=0}^\infty \tld{\eta}_{(k)} \,.
\end{gather*}
We will show in Section~\ref{se:VFMdet} that $\tld{\eta}_{(0)}=\tld{\eta}_{(1)}=0$.
The Taylor series of the function $S$ also starts from the second term,
\begin{gather*}
S(x(z))=\sum_{k=2}^\infty \tld{S}_{(k)} \,.
\end{gather*}

\subsection{Scaling of Coordinates}

Now we rescale the coordinates
\begin{gather*}
z^\alpha \to \sqrt{t} z^\alpha \,.
\end{gather*}
Then all terms in covariant Taylor series are transformed as
\begin{gather*}
\tld{f}_{(k)} \to t^{k/2} \tld{f}_{(k)} \,,
\end{gather*}
and, hence, they become asymptotic expansions in $t^{1/2}$.
By substituting expansions from the previous section and~\eqref{eq:OmegaAE} into~\eqref{eq:curvIM'}, we obtain
\begin{gather} \label{eq:4.5}
\Theta_\Phi(t) \sim
	(4\pi)^{-n/2}
	\int_{\norm{z}<\eps/\sqrt{t}}  dz\,
	\abs{g(x_0)}^{1/2}
	\e{-\dfrac{1}{2}\,\tld{S}_{(2)}}
	B(t)
	\,,
\end{gather}
where
\begin{gather*}
B(t)=
	\e{-\sum_{k=1}^\infty t^{k/2}\, h_k}
	\sum_{m=0}^\infty \frac{(-t)^m}{m!}
	\sum_{k=0}^\infty t^{k/2}\, \tld{a}_{m(k)}
	\,,
\end{gather*}
and
\begin{gather} \label{eq:h_k}
h_k = \dfrac{1}{2}\,\tld{S}_{(k+2)} + 2\,\tld{\zeta}_{(k)} - \tld{\eta}_{(k)} \,.
\end{gather}
Note, that terms $\tld{S}_{(k)}$ and $\tld{\eta}_{(k)}$ depend on the mapping $\Phi$, while $\tld{\zeta}_{(k)}$ does not.
Let $b_k$ be the coefficients of a power series in $t^{1/2}$ for the function $B$,
\begin{gather} \label{eq:4.9}
B(t)=\sum_{k=0}^\infty t^{k/2}\, b_k \,.
\end{gather}
It is easy to see, that the coefficients $b_k$ are polynomials in $z$. Moreover, they are polynomials in $h_k$ and $\tld{a}_{m(k)}$. From~\cite{Avramidi2000} and our computations in Section~\ref{se:CurvePointComputation} it follows, that $h_k$ and $\tld{a}_{m(k)}$ can be expressed in terms of covariant derivatives of the curvature of the metric $g$ and symmetrized covariant derivatives of the differential of the mapping $\Phi$. Note, also, that half-integer (integer) powers of $t$ correspond to polynomials of odd (even) degree in $z$. By substituting~\eqref{eq:4.9} into~\eqref{eq:4.5}, we obtain
\begin{gather*}
\Theta_\Phi(t) \sim
	\sum_{k=0}^\infty
	t^{k/2}\, 
	(4\pi)^{-n/2}
	\int_{\norm{z}<\eps/\sqrt{t}}  dz\,
	\abs{g(x_0)}^{1/2}
	\e{-\dfrac{1}{2}\,\tld{S}_{(2)}}
	b_k
	\,.
\end{gather*}
In the same way, as in Lemma~\ref{le:GIAsymptotics1d}, we expand the integration over the set $\norm{z}<\eps/\sqrt{t}$ to the integration over the whole $\R^n$, i.e. one can show that
\begin{gather*}
\int_{\norm{z}<\eps/\sqrt{t}} dz\, \e{-\dfrac{1}{2}\,\tld{S}_{(2)}} z^{i_1}\ldots z^{i_k}
\sim
\int_{\R^n} dz\, \e{-\dfrac{1}{2}\,\tld{S}_{(2)}} z^{i_1}\ldots z^{i_k}
\,,
\end{gather*}
for both odd and even $k$. After this operation all terms with half-integer powers of $t$ disappear, since they contain integrals of odd polynomials over a symmetric set. Therefore, we have
\begin{gather*}
\Theta_\Phi(t) \sim
	\sum_{k=0}^\infty
	t^k\, 
	(4\pi)^{-n/2}
	\int_{\R^n}  dz\,
	\abs{g(x_0)}^{1/2}
	\e{-\dfrac{1}{2}\,\tld{S}_{(2)}}
	b_{2k}
	\,.
\end{gather*}
By comparing this expression with~\eqref{eq:Expansion-curved-point}, we finally obtain
\begin{gather} \label{eq:A_k-curved-point}
A_k=(4\pi)^{-n/2}\, 
	\int_{\R^n}  dz\,
	\abs{g(x_0)}^{1/2}
	\e{-\dfrac{1}{2}\,\tld{S}_{(2)}}
	b_{2k}
	\,.
\end{gather}
In this expression:
\begin{enumerate}
\item $\tld{S}_{(2)}=\dfrac{1}{2} S_{\alpha\beta}z^\alpha z^\beta$, where $S_{\alpha\beta}$ is the positive-definite quadratic form, defined by~\eqref{eq:S2Matrix}.
\item $b_{2k}$ are polynomials in $z$, defined by~\eqref{eq:4.9}.
\item All integrals in~\eqref{eq:A_k-curved-point} are standard Gaussian integrals which can be computed by using Lemma~~\ref{le:GaussianIntegralmd}.
\item $b_{2k}$ are polynomials in covariant derivatives of the curvature of the metric $g$ and symmetrized covariant derivatives of the differential $d\Phi$ of the mapping~$\Phi$. Due to Lemma~~\ref{le:GaussianIntegralmd}, Gaussian integrals in~\eqref{eq:A_k-curved-point} depend on the inverse of the matrix $S_{\alpha\beta}$, that is, the matrix $Q$.
\end{enumerate}
This completes the proof of the theorem.
\end{proof}

\section{Computation of the Coefficients $A_0, A_1,$ and $A_2$} \label{se:CurvePointComputation}

In the remaining part of the chapter we will develop an algorithm for computation of the coefficients $A_k$ in terms of covariant derivatives of the differential $d\Phi$ and the curvature of the metric $g$. We will compute coefficients $A_0, A_1$, and $A_2$ explicitly. In the next section we describe tools for construction of Taylor series of two-point functions.

\subsection{Coincidence Limits of Two-Point Functions}

Let $f:M\times M\to R$ be a two-point function defined on a manifold $M$. We will denote by $[f]$ its coincidence limit, that is,
\begin{gather*}
[f](x)=\lim_{x'\to x} f(x,x') \,.
\end{gather*}
Obviously, covariant derivatives of the function $f$ at different points commute with each other. We will use primes for indices of covariant derivatives at the point $x'$. For coincidence limits of mixed derivatives of the function $f$ there holds~\cite{Synge}
\begin{gather*}
[f_{\ldots; \alpha'}]=[f_{\ldots}]_{;\alpha}-[f_{\ldots;\alpha}] \,,
\end{gather*}
where dots stand for a fixed set of some indices and may include other covariant derivatives. By using this property repeatedly, we obtain
\begin{align}
\label{eq:CoincDer2}
[f_{\ldots;\alpha'\beta'}] &=
	[f_{\ldots}]_{;\alpha\beta}
	-[f_{\ldots;\alpha}]_{;\beta}
	-[f_{\ldots;\beta}]_{;\alpha}
	+[f_{\ldots;\beta\alpha}]
	\,,\\
\label{eq:CoincDer3}
[f_{\ldots;\alpha'\beta'\gamma'}] &=
	[f_{\ldots}]_{;\alpha\beta\gamma}
	-[f_{\ldots;\gamma}]_{;\alpha\beta}
	-[f_{\ldots;\beta}]_{;\alpha\gamma}
	-[f_{\ldots;\alpha}]_{;\beta\gamma}
	\\ \nonumber &
	+[f_{\ldots;\gamma\beta}]_{;\alpha}
	+[f_{\ldots;\gamma\alpha}]_{;\beta}
	+[f_{\ldots;\beta\alpha}]_{;\gamma}
	-[f_{\ldots;\gamma\beta\alpha}]
	\,.
\end{align}
One can easily obtain the following general formula of this kind for
$[f_{\ldots;\mu_1'\ldots\mu_m'}]$.
Let $\N_m=\set{1,\ldots,m}$ and $\PP(\N_m)$ be its power set.
For any $I\in \PP(\N_m)$, such that $I=\set{i_1,\ldots,i_k}$ and $i_1<\ldots<i_k$, let
$\underrightarrow{\mu(I)}$ be the increasing sequence of indices $\mu_{i_1}\ldots\mu_{i_k}$ and
$\underleftarrow{\mu(I)}$ be the decreasing sequence of indices $\mu_{i_k}\ldots\mu_{i_1}$. Then
\begin{gather*}
[f_{\ldots;\mu_1'\ldots\mu_m'}] =
	\sum_{I\in \PP(\N_m)}\
	(-1)^{\abs{I}}\,
	[f_{\ldots;\underleftarrow{\mu(I)}}]_{;\underrightarrow{\mu(J)}} \,,
\end{gather*}
where $J=\N_m\setminus I$.

\subsection{Derivatives of the Deformed Diagonal Functions} \label{se:DerDefDiag}

Let $f$ be a two-point function $f:M\times M\to\R$. In this section we will develop an algorithm for computing symmetrized covariant derivatives of the function $g(x)=f(x,\Phi(x))$. We will agree, that indices with primes denote covariant differentiation of the function $f$ with respect to the second argument and indices without primes~--- with respect to the first argument. Obviously, derivatives with respect to different arguments commute with each other.

Let $\tptb{p}{q}{r}{s}$ be the set of two-point tensors of the type $(p,q)$ at a point $x$ and of the type $(r,s)$ at a point $x'$, that is, the set of linear mappings
\begin{multline*}
T:
\underbrace{T_x M\times\ldots\times T_x M}_q
\times
\underbrace{T^*_x M\times\ldots\times T^*_x M}_p
\\
\times
\underbrace{T_{x'} M\times\ldots\times T_{x'} M}_{s}
\times
\underbrace{T^*_{x'} M\times\ldots\times T^*_{x'} M}_{r}
\to\R \,.
\end{multline*}
If some of the numbers $p,q,r,s$ are equal to zero, we will omit them. The external tensor product of two-point tensors
$A\in \tptb{p}{q}{r}{s}$
and
$B\in \tptb{t}{u}{v}{w}$
is the two-point tensor
$A \boxtimes B \in \tptb{p+t}{q+u}{r+v}{s+w}$,
defined by
\begin{gather*}
(A \boxtimes B)
	^{\alpha_1\ldots \alpha_p
	  \mu_1\ldots \mu_t
	  \gamma'_1\ldots \gamma'_r
	  \xi'_1\ldots \xi'_v}
	_{\beta_1\ldots \beta_q
	  \nu_1\ldots \nu_u
	  \delta'_1\ldots \delta'_s
	  \pi'_1\ldots \pi'_w}
=
A
	^{\alpha_1\ldots \alpha_p
	  \gamma'_1\ldots \gamma'_r}
	_{\beta_1\ldots \beta_q
	  \delta'_1\ldots \delta'_s}
B
	^{\mu_1\ldots \mu_t
	  \xi'_1\ldots \xi'_v}
	_{\nu_1\ldots \nu_u
	  \pi'_1\ldots \pi'_w}
\,.
\end{gather*}
We define a mapping $\bullet:\tptb{}{p}{}{q}\times \tptb{}{r}{q}{}\to\tptb{}{p+r}{}{}$ by
\begin{gather*}
(A \bullet B)_{\alpha_1\ldots \alpha_p \mu_1\ldots \mu_r}
=A_{\alpha_1\ldots \alpha_p \beta'_1\ldots \beta'_q} B^{\beta'_1\ldots \beta'_q}_{\mu_1\ldots \mu_r} \,.
\end{gather*}

It is easy to see, that the $n$-th symmetrized derivative of the function $g$ is a linear combination of terms of the form
\begin{gather*}
T(k,m;j_1,\ldots,j_k)=
\Sym \Big(
\left( \nabla'^k \nabla^m f \right)
\bullet
\left(
\left( \nabla^{j_1} d\Phi \right) \boxtimes \ldots \boxtimes \left( \nabla^{j_k} d\Phi \right)
\right)
\Big)
\,,
\end{gather*}
where numbers $k,m,j_1,\ldots,j_k$ are non-negative integers, such that
\begin{gather} \label{eq:Sder-kmjn}
k+m+j_1+\ldots+j_k=n \,.
\end{gather}
A simple computation shows that there exist $2^n$ such terms.

Let numbers $c(k,m;j_1,\ldots,j_k)$ be the coefficients of this linear combination, that is,
\begin{gather*}
\Sym(\nabla^n g)=\sum_{k=0}^n\ \sum_{m=0}^{n-k}\
	\sum_{\clearfrac{\mi{j}=(j_1,\ldots,j_k)}{\abs{\mi{j}}=n-m-k}}\
	c(k,m;j_1,\ldots,j_k)\ T(k,m;j_1,\ldots,j_k) \,.
\end{gather*}
We will derive a recursive formula for these coefficients. Note, that the term
$T(k,m;j_1,\ldots,j_k)$ of the $n$-th derivative can be obtained from the following terms of the $(n-1)$-th derivative:
\begin{enumerate}
\item $T(k,m-1;j_1,\ldots,j_k)$, by differentiating $f$ with respect to the first argument, if $m\geqslant 1$.
\item $T(k-1,m;j_1,\ldots,j_{k-1})$, by differentiating $f$ with respect to the second argument, if $k\geqslant 1$ and $j_k=0$.
\item $T(k,m;j_1,\ldots,j_i-1,\ldots,j_k)$, by differentiating $d\Phi$, if $j_i\geqslant 1$.
\end{enumerate}
Let $c(k,m;j_1,\ldots,j_k)=0$, if one or more arguments are negative. The above rules give us
\begin{align*}
c(k,m;j_1,\ldots,j_k) &
	=c(k,m-1;j_1,\ldots,j_k)
	+\sum_{i=1}^k c(k,m;j_1,\ldots,j_i-1,\ldots,j_k) \,,
\end{align*}
if $j_k>0$, and
\begin{align*}
c(k,m;j_1,\ldots,j_k) &
	=c(k,m-1;j_1,\ldots,j_k)
	+\sum_{i=1}^k c(k,m;j_1,\ldots,j_i-1,\ldots,j_k)
	\\&
	+c(k-1,m;j_1,\ldots,j_{k-1}) \,,
\end{align*}
if $j_k=0$.
The obvious initial condition for this recursion is
\begin{gather*}
c(0,0;)=1.
\end{gather*}
One can also get a non-recursive formula for the coefficients $c(k,m;j_1,\ldots,j_k)$.

\begin{lemma}
Let $n\in\N$ and $k,m,j_1,\ldots,j_k$ be non-negative integers, satisfying~\eqref{eq:Sder-kmjn}. Then
\begin{gather*}
c(k,m;j_1,\ldots,j_k)=
	\frac{n!}{m! j_1! \ldots j_k! (1+j_k) (2+j_k+j_{k-1}) \ldots (k+j_k+\ldots+j_1)}
	\,.
\end{gather*}
\end{lemma}

\begin{proof}
It is not hard to check that the non-symmetrized covariant derivative $g_{\nu_1\ldots\nu_n}$ is the sum of all terms of the form
\begin{gather*}
\Big(
\left( \nabla'^k \nabla^m f \right)
\bullet
\left(
\left( \nabla^{j_1} d\Phi \right) \boxtimes \ldots \boxtimes \left( \nabla^{j_k} d\Phi \right)
\right)
\Big)
_{\nu_{\phi(1)}\ldots\nu_{\phi(n)}}
\,,
\end{gather*}
for all possible values of $k,m,j_1,\ldots,j_k$ and all permutations $\phi$ of numbers $1,\ldots,n$ satisfying certain conditions, that we will determine below. Symmetrization of indices $\nu_1\ldots\nu_n$ makes all terms with fixed $k,m,j_1,\ldots,j_k$ the same. Hence, the coefficient $c(k,m;j_1,\ldots,j_k)$ is the number of all possible permutations $\phi$.

Let $J_i$ be defined for $i=1,\ldots,k$ by
\begin{gather*}
J_i=m+i+j_1+\ldots+j_{i-1} \,.
\end{gather*}
Note that indices $\nu_{\phi(1)}\ldots\nu_{\phi(m)}$ correspond to the derivative of $f$, index $\nu_{\phi(J_i)}$ corresponds to the $i$-th factor $d\Phi$, and indices $\nu_{\phi(J_i+1)}\ldots\nu_{\phi(J_i+j_i)}$ correspond to the derivative of this factor. This observation leads us to the following conditions on the permutation $\phi$:
\begin{enumerate}
\item $\phi(1)<\phi(2)<\ldots<\phi(m)\,$;
\item $\phi(J_i)<\phi(J_i+1)<\ldots<\phi(J_i+j_i)\,$;
\item $\phi(J_1)<\phi(J_2)<\ldots<\phi(J_k)\,$.
\end{enumerate}
In words, indices on every factor and the first indices of factors $d\Phi$ must be in the ``increasing alphabetical order.'' So, there are $\binom{n}{m}$ ways to choose indices corresponding to the derivative of $f$. After this step $\phi(J_1)$ is determined uniquely and there are $\binom{n-m-1}{j_1}$ ways to choose indices corresponding to the derivative of the first factor $d\Phi$. By continuing this reasoning, we get
\begin{multline*}
c(k,m;j_1,\ldots,j_k)=
	\binom{n}{m}
	\binom{n-m-1}{j_1}
	\\[2ex]  \times
	\binom{n-m-2-j_1}{j_2}
	\ldots
	\binom{n-m-(k-1)-j_1-\ldots-j_{k-2}}{j_{k-1}}
	\,.
\end{multline*}
After straightforward transformations, we get the desired formula.
\end{proof}

\subsection{Symmetrized Derivatives of the Function $S$} \label{se:SymDerOfS}

Now we will apply the technique developed in Section~\ref{se:DerDefDiag} to the computation of symmetrized covariant derivatives of the function $S(x)=\sigma(x,\Phi(x))$ up to the sixth order.
Since we are interested in values of the derivatives on the fixed point set $\Sigma$, all appearing derivatives of the world function will be actually equal to their coincidence limits. We will take into account the following properties of the world function~\cite{Avramidi2000}:
\begin{enumerate}
\item The world function is symmetric in its arguments.
\item The first two derivatives with respect to the same argument commute.
\item $[\sigma]=0$.
\item $[\sigma_\alpha]=0$.
\item $[\sigma_{\alpha\beta}]=-[\sigma_{\alpha\beta'}]=g_{\alpha\beta}$.
\item $[\sigma_{\alpha\beta\gamma}]=[\sigma_{\alpha\beta\gamma'}]=0$.
\item $[\sigma_{(\mu_1\ldots\mu_k)}]=0$, $k\geqslant 3$.
\item $[\sigma_{\alpha(\mu_1\ldots\mu_k)}]=[\sigma_{\alpha'(\mu_1\ldots\mu_k)}]=0$, $k\geqslant 2$.
\end{enumerate}

The recursive algorithm, described above, for computing the coef\-ficients $c(k,m;j_1,\ldots,j_k)$ was realized as the Python script listed in Appen\-dix~\ref{ap:Script}.  This script also generates preliminary expressions for derivatives of the function $S$ in the \LaTeX\ format.

For the derivatives of the differential $d\Phi$ of the mapping $\Phi$, evaluated at the point $x_0$, the notation
\begin{gather*}
\Psi^{\mu}{}_{\nu_1\nu_2\ldots\nu_k}=(d\Phi)^{\mu'}{}_{\nu_1;\nu_2\ldots\nu_k}(x_0)
\end{gather*}
will be used. The first six symmetrized derivatives of $S$ are
\begin{align*}
S&=0
\,,\\
S_{\alpha}&=0
\,,\\
S_{(\alpha\beta)}&=g_{\alpha\beta}-2\, \Psi_{(\alpha\beta)}+\Psi_{\mu(\alpha} \Psi^{\mu}{}_{\beta)}
\,,\\
S_{(\alpha\beta\gamma)}&=
-3\, \Psi_{(\alpha\beta\gamma)}
+3\, \Psi_{\mu(\alpha} \Psi^{\mu}{}_{\beta\gamma)}
\,,\\
S_{(\alpha\beta\gamma\delta)}&=
-4\, \Psi_{(\alpha\beta\gamma\delta)}
+4\, \Psi_{\mu(\alpha} \Psi^{\mu}{}_{\beta\gamma\delta)}
+3\, \Psi_{\mu(\alpha\beta} \Psi^{\mu}{}_{\gamma\delta)}
+6\, [\sigma_{\mu'\nu'(\alpha\beta}] \Psi^{\mu}{}_{\gamma} \Psi^{\nu}{}_{\delta)}
\,,\\
S_{(\alpha\beta\gamma\delta\eps)}&=
-5\, \Psi_{(\alpha\beta\gamma\delta\eps)}
+5\, \Psi_{\mu(\alpha} \Psi^{\mu}{}_{\beta\gamma\delta\eps)}
+10\, \Psi_{\mu(\alpha\beta} \Psi^{\mu}{}_{\gamma\delta\eps)}
\displaybreak[0]\\ &
+30\, [\sigma_{\mu'\nu'(\alpha\beta}] \Psi^{\mu}{}_{\gamma} \Psi^{\nu}{}_{\delta\eps)}
\displaybreak[0]\\ &
+[\sigma_{\mu'\nu'\xi'(\alpha}] \Big\{
	5\, \Psi^{\mu}{}_{\beta} \Psi^{\nu}{}_{\gamma} \Psi^{\xi}{}_{\delta\eps)}
	+25\, \Psi^{\mu}{}_{\beta} \Psi^{\nu}{}_{\gamma\delta} \Psi^{\xi}{}_{\eps)}
	\Big\}
\displaybreak[0]\\ &
+10\, [\sigma_{\mu'\nu'(\alpha\beta\gamma}] \Psi^{\mu}{}_{\delta} \Psi^{\nu}{}_{\eps)}
+10\, [\sigma_{\mu'\nu'\xi'(\alpha\beta}] \Psi^{\mu}{}_{\gamma} \Psi^{\nu}{}_{\delta} \Psi^{\xi}{}_{\eps)}
\,,\\
S_{(\alpha\beta\gamma\delta\eps\zeta)}&=
-6\, \Psi_{(\alpha\beta\gamma\delta\eps\zeta)}
+6\, \Psi_{\mu(\alpha} \Psi^{\mu}{}_{\beta\gamma\delta\eps\zeta)}
+15\, \Psi_{\mu(\alpha\beta} \Psi^{\mu}{}_{\gamma\delta\eps\zeta)}
+10\, \Psi_{\mu(\alpha\beta\gamma} \Psi^{\mu}{}_{\delta\eps\zeta)}
\displaybreak[0]\\ &
+[\sigma_{\mu'\nu'(\alpha\beta}] \Big\{
	60\, \Psi^{\mu}{}_{\gamma} \Psi^{\nu}{}_{\delta\eps\zeta)}
	+45\, \Psi^{\mu}{}_{\gamma\delta} \Psi^{\nu}{}_{\eps\zeta)}
	\Big\}
\displaybreak[0]\\ &
+[\sigma_{\mu'\nu'\xi'(\alpha}] \Big\{
	6\, \Psi^{\mu}{}_{\beta} \Psi^{\nu}{}_{\gamma} \Psi^{\xi}{}_{\delta\eps\zeta)}
	+42\, \Psi^{\mu}{}_{\beta} \Psi^{\nu}{}_{\gamma\delta} \Psi^{\xi}{}_{\eps\zeta)}
\displaybreak[0]\\ &
	+54\, \Psi^{\mu}{}_{\beta} \Psi^{\nu}{}_{\gamma\delta\eps} \Psi^{\xi}{}_{\zeta)}
	+48\, \Psi^{\mu}{}_{\beta\gamma} \Psi^{\nu}{}_{\delta\eps} \Psi^{\xi}{}_{\zeta)}
	\Big\}
\displaybreak[0]\\ &
+[\sigma_{\mu'\nu'\xi'\pi'}] \Big\{
	3\, \Psi^{\mu}{}_{(\alpha} \Psi^{\nu}{}_{\beta} \Psi^{\xi}{}_{\gamma\delta} \Psi^{\pi}{}_{\eps\zeta)}
	+9\, \Psi^{\mu}{}_{(\alpha} \Psi^{\nu}{}_{\beta\gamma} \Psi^{\xi}{}_{\delta} \Psi^{\pi}{}_{\eps\zeta)}
\displaybreak[0]\\ &
	+18\, \Psi^{\mu}{}_{(\alpha} \Psi^{\nu}{}_{\beta\gamma} \Psi^{\xi}{}_{\delta\eps} \Psi^{\pi}{}_{\zeta)}
	+15\, \Psi^{\mu}{}_{(\alpha\beta} \Psi^{\nu}{}_{\gamma\delta} \Psi^{\xi}{}_{\eps} \Psi^{\pi}{}_{\zeta)}
	\Big\}
\displaybreak[0]\\ &
+60\, [\sigma_{\mu'\nu'(\alpha\beta\gamma}] \Psi^{\mu}{}_{\delta} \Psi^{\nu}{}_{\eps\zeta)}
\displaybreak[0]\\ &
+[\sigma_{\mu'\nu'\xi'(\alpha\beta}] \Big\{
	15\, \Psi^{\mu}{}_{\gamma} \Psi^{\nu}{}_{\delta} \Psi^{\xi}{}_{\eps\zeta)}
	+75\, \Psi^{\mu}{}_{\gamma} \Psi^{\nu}{}_{\delta\eps} \Psi^{\xi}{}_{\zeta)}
	\Big\}
\displaybreak[0]\\ &
+[\sigma_{\mu'\nu'\xi'\pi'(\alpha}] \Big\{
	6\, \Psi^{\mu}{}_{\beta} \Psi^{\nu}{}_{\gamma} \Psi^{\xi}{}_{\delta} \Psi^{\pi}{}_{\eps\zeta)}
	+12\, \Psi^{\mu}{}_{\beta} \Psi^{\nu}{}_{\gamma} \Psi^{\xi}{}_{\delta\eps} \Psi^{\pi}{}_{\zeta)}
\displaybreak[0]\\ &
	+42\, \Psi^{\mu}{}_{\beta} \Psi^{\nu}{}_{\gamma\delta} \Psi^{\xi}{}_{\eps} \Psi^{\pi}{}_{\zeta)}
	\Big\}
\displaybreak[0]\\ &
+15\, [\sigma_{\mu'\nu'(\alpha\beta\gamma\delta}] \Psi^{\mu}{}_{\eps} \Psi^{\nu}{}_{\zeta)}
+20\, [\sigma_{\mu'\nu'\xi'(\alpha\beta\gamma}] \Psi^{\mu}{}_{\delta} \Psi^{\nu}{}_{\eps} \Psi^{\xi}{}_{\zeta)}
\displaybreak[0]\\ &
+15\, [\sigma_{\mu'\nu'\xi'\pi'(\alpha\beta}] \Psi^{\mu}{}_{\gamma} \Psi^{\nu}{}_{\delta} \Psi^{\xi}{}_{\eps} \Psi^{\pi}{}_{\zeta)}
\,.
\end{align*}
The second derivative can also be written in the form
\begin{align}
S_{\alpha\beta}&= \nonumber
	g_{\alpha\beta}-2\, \Psi_{(\alpha\beta)}+\Psi_{\mu(\alpha} \Psi^{\mu}{}_{\beta)}
	\\ &= \nonumber
	g_{\alpha\beta}-\Psi_{\alpha\beta}-\Psi_{\beta\alpha}+\Psi^\mu{}_\alpha g_{\mu\nu} \Psi^\nu{}_\beta
	\\ &= \label{eq:S_alphabeta}
	\big(
	\delta^\mu{}_\alpha
	-(d\Phi)^\mu{}_\alpha
	\big)
 	g_{\mu\nu}
 	\big(
	\delta^\nu{}_\beta
	-(d\Phi)^\nu{}_\beta
	\big)
	\,.
\end{align}

We need the following coincidence limits for derivatives of the world function:
\begin{gather*}
[\sigma_{\mu'\nu'\alpha\beta}] \,,\
[\sigma_{\mu'\nu'\xi'\alpha}] \,,\
[\sigma_{\mu'\nu'\xi'\pi'}] \,,\\
[\sigma_{\mu'\nu'\xi'\alpha\beta}] \,,\
[\sigma_{\mu'\nu'\xi'\pi'\alpha}] \,,\\
[\sigma_{\mu'\nu'(\alpha\beta\gamma\delta)}] \,,\
[\sigma_{(\mu'\nu'\xi')(\alpha\beta\gamma)}] \,.
\end{gather*}
We will reduce them to the limits which are known~\cite{Avramidi2000}
\begin{align*}
[\sigma_{\alpha'\beta (\mu_1\mu_2)}]&=
	-\frac{1}{3}\, R_{\alpha(\mu_1|\beta|\mu_2)}
	\,,\displaybreak[0]\\[1ex]
[\sigma_{\alpha'\beta (\mu_1\mu_2\mu_3)}]&=
	-\frac{1}{2}\, R_{\alpha(\mu_1|\beta|\mu_2;\mu_3)}
	\,,\displaybreak[0]\\[1ex]
[\sigma_{\alpha'\beta (\mu_1\mu_2\mu_3\mu_4)}]&=
	-\frac{3}{5}\, R_{\alpha(\mu_1|\beta|\mu_2;\mu_3\mu_4)}
	-\frac{7}{15}\, R_{\alpha(\mu_1|\gamma|\mu_2} R^\gamma{}_{\mu_3|\beta|\mu_4)}
	\,.
\end{align*}
For this purpose we will use properties of coincidence limits~\eqref{eq:CoincDer2}, \eqref{eq:CoincDer3}, expression for the commutator of covariant derivatives~\eqref{eq:CovDerCom}, symmetry properties of the curvature tensor~\eqref{eq:CurvSym}, and Bianci identities~\eqref{eq:CurvBianci}, \eqref{eq:CurvBianciDer}. After straightforward but somewhat cumbersome computation we have:
\begin{align*}
[\sigma_{\mu'\nu'\alpha\beta}] &=
	-\frac{2}{3}\, R_{\alpha(\mu|\beta|\nu)}
	\,,\displaybreak[0]\\
[\sigma_{\mu'\nu'\xi'\alpha}] &=
	\frac{2}{3}\, R_{\alpha(\mu|\xi|\nu)}
	\,,\displaybreak[0]\\
[\sigma_{\mu'\nu'\xi'\pi'}] &=
	-\frac{2}{3}\, R_{\pi(\mu|\xi|\nu)}
	\,,\displaybreak[0]\\
[\sigma_{\mu'\nu'\xi'\alpha\beta}] &=
	\frac{1}{4}\, R_{(\mu|\xi|\nu)(\alpha;\beta)}
	-\frac{5}{12}\, R_{(\mu|\alpha|\nu)\beta;\xi}
	-\frac{1}{12}\, R_{(\alpha|\xi|\beta)(\mu;\nu)}
	\,,\displaybreak[0]\\
[\sigma_{\mu'\nu'\xi'\pi'\alpha}] &=
	\frac{5}{6}\, R_{(\mu|\alpha|\nu)(\xi;\pi)}
	-\frac{1}{6}\, R_{(\xi|\alpha|\pi)(\mu;\nu)}
	\,,\displaybreak[0]\\
[\sigma_{\mu'\nu'(\alpha\beta\gamma\delta)}] &=
	-\frac{2}{5}\, R_{\mu(\alpha|\nu|\beta;\gamma\delta)}
	-\frac{8}{15}\, R_{ a (\alpha|\mu|\beta} R^a{}_{\gamma|\nu|\delta)}
	\,,\displaybreak[0]\\
[\sigma_{(\mu'\nu'\xi')(\alpha\beta\gamma)}] &=
	\left(
	-\frac{3}{10}\, R^{(\pi}{}_{(\alpha}{}^\rho{}_{\beta;}{}^{\tau)}{}_{\gamma)}
	-\frac{4}{3}\, R_a{}^{(\pi}{}_{(\alpha}{}^\rho R^{|a|}{}_\beta{}^{\tau)}{}_{\gamma)}
	\right)
	g_{\pi\mu} g_{\rho\nu} g_{\tau\xi}
	\,.
\end{align*}
Detailed computation is given in Appendix~\ref{ap:CLofWF}.

By combining obtained expressions, we finally have
\begin{align}
S_{(\alpha\beta)}&=g_{\alpha\beta}-2\, \Psi_{(\alpha\beta)}+\Psi_{\mu(\alpha} \Psi^{\mu}{}_{\beta)}
\label{eq:S2curvature}
\,,\\
S_{(\alpha\beta\gamma)}&=
-3\, \Psi_{(\alpha\beta\gamma)}
+3\, \Psi_{\mu(\alpha} \Psi^{\mu}{}_{\beta\gamma)}
\,,\\
S_{(\alpha\beta\gamma\delta)}&=
-4\, R_{(\alpha|\mu|\beta|\nu|} \Psi^{\mu}{}_{\gamma} \Psi^{\nu}{}_{\delta)}
+3\, \Psi_{\mu(\alpha\beta} \Psi^{\mu}{}_{\gamma\delta)}
+4\, \Psi_{\mu(\alpha} \Psi^{\mu}{}_{\beta\gamma\delta)}
-4\, \Psi_{(\alpha\beta\gamma\delta)}
\,,\\
S_{(\alpha\beta\gamma\delta\eps)}&=
-5\, R_{\mu(\alpha|\nu|\beta;\gamma}
\Psi^{\mu}{}_{\delta} \Psi^{\nu}{}_{\eps)}
-5\, R_{\mu(\alpha|\nu|\beta;|\xi|}
\Psi^{\mu}{}_{\gamma} \Psi^{\nu}{}_{\delta} \Psi^{\xi}{}_{\eps)}
\nonumber
\displaybreak[0]\\ &
-20\, R_{(\alpha|\mu|\beta|\nu} \Psi^{\mu}{}_{\gamma} \Psi^{\nu}{}_{\delta\eps)}
+\frac{10}{3}\, R_{(\mu|\xi|\nu)(\alpha} \Big\{
	\Psi^{\mu}{}_{\beta} \Psi^{\nu}{}_{\gamma} \Psi^{\xi}{}_{\delta\eps)}
	+5\, \Psi^{\mu}{}_{\beta} \Psi^{\nu}{}_{\gamma\delta} \Psi^{\xi}{}_{\eps)}
	\Big\}
\nonumber
\displaybreak[0]\\ &
+10\, \Psi_{\mu(\alpha\beta} \Psi^{\mu}{}_{\gamma\delta\eps)}
-5\, \Psi_{(\alpha\beta\gamma\delta\eps)}
+5\, \Psi_{\mu(\alpha} \Psi^{\mu}{}_{\beta\gamma\delta\eps)}
\,,\\
S_{(\alpha\beta\gamma\delta\eps\zeta)}&=
-\Big\{
	6\, R_{\mu(\alpha|\nu|\beta;\gamma\delta}
	+8\, R_{ \rho (\alpha|\mu|\beta} R^\rho{}_{\gamma|\nu|\delta}
\Big\} \Psi^{\mu}{}_{\eps} \Psi^{\nu}{}_{\zeta)}
\nonumber
\displaybreak[0]\\ &
-\frac{2}{3}\Big\{
	9\, R_{\mu(\alpha|\nu|\beta;|\xi|\gamma}
	+40\, R_{ \rho \mu(\alpha|\nu|} R^\rho{}_{\beta|\xi|\gamma}
\Big\} \Psi^{\mu}{}_{\delta} \Psi^{\nu}{}_{\eps} \Psi^{\xi}{}_{\zeta)}
\nonumber
\displaybreak[0]\\ &
-\Big\{
	6\, R_{\mu(\alpha|\nu|\beta;|\xi\pi|}
	+8\, R_{ \rho \mu\nu(\alpha} R^\rho{}_{|\xi\pi|\beta}
\Big\} \Psi^{\mu}{}_{\gamma} \Psi^{\nu}{}_{\delta} \Psi^{\xi}{}_{\eps} \Psi^{\pi}{}_{\zeta)}
\nonumber
\displaybreak[0]\\ &
-30\, R_{\mu(\alpha|\nu|\beta;\gamma} \Psi^{\mu}{}_{\delta} \Psi^{\nu}{}_{\eps\zeta)}
\nonumber
\displaybreak[0]\\ &
+\frac{5}{8}\, \Big\{
	6\, R_{(\mu|\xi|\nu)(\alpha;\beta}
	-10\, R_{\mu(\alpha|\nu|\beta;|\xi|}
	-R_{\xi(\alpha|\mu|\beta;|\nu|}
	-R_{\xi(\alpha|\nu|\beta;|\mu|}
\Big\}
\nonumber
\displaybreak[0]\\ &
\times \Big\{
	\Psi^{\mu}{}_{\gamma} \Psi^{\nu}{}_{\delta} \Psi^{\xi}{}_{\eps\zeta)}
	+5\, \Psi^{\mu}{}_{\gamma} \Psi^{\nu}{}_{\delta\eps} \Psi^{\xi}{}_{\zeta)}
	\Big\}
\nonumber
\displaybreak[0]\\ &
+\Big\{
	5\, R_{(\mu}{}^\rho{}_{\nu)(\xi;\pi)}
	-R_{(\xi}{}^\rho{}_{\pi)(\mu;\nu)}
\Big\} g_{\rho(\alpha} \Big\{
	\Psi^{\mu}{}_{\beta} \Psi^{\nu}{}_{\gamma} \Psi^{\xi}{}_{\delta} \Psi^{\pi}{}_{\eps\zeta)}
\nonumber
\displaybreak[0]\\ &
	+2\, \Psi^{\mu}{}_{\beta} \Psi^{\nu}{}_{\gamma} \Psi^{\xi}{}_{\delta\eps} \Psi^{\pi}{}_{\zeta)}
	+7\, \Psi^{\mu}{}_{\beta} \Psi^{\nu}{}_{\gamma\delta} \Psi^{\xi}{}_{\eps} \Psi^{\pi}{}_{\zeta)}
	\Big\}
\nonumber
\displaybreak[0]\\ &
-10\, R_{\mu(\alpha|\nu|\beta} \Big\{
	4\, \Psi^{\mu}{}_{\gamma} \Psi^{\nu}{}_{\delta\eps\zeta)}
	+3\, \Psi^{\mu}{}_{\gamma\delta} \Psi^{\nu}{}_{\eps\zeta)}
	\Big\}
\nonumber
\displaybreak[0]\\ &
+4\, R_{(\mu|\xi|\nu)(\alpha} \Big\{
	\Psi^{\mu}{}_{\beta} \Psi^{\nu}{}_{\gamma} \Psi^{\xi}{}_{\delta\eps\zeta)}
	+7\, \Psi^{\mu}{}_{\beta} \Psi^{\nu}{}_{\gamma\delta} \Psi^{\xi}{}_{\eps\zeta)}
\nonumber
\displaybreak[0]\\ &
	+9\, \Psi^{\mu}{}_{\beta} \Psi^{\nu}{}_{\gamma\delta\eps} \Psi^{\xi}{}_{\zeta)}
	+8\, \Psi^{\mu}{}_{\beta\gamma} \Psi^{\nu}{}_{\delta\eps} \Psi^{\xi}{}_{\zeta)}
	\Big\}
\nonumber
\displaybreak[0]\\ &
-2\, R_{\xi(\mu|\pi|\nu)} \Big\{
	\Psi^{\mu}{}_{(\alpha} \Psi^{\nu}{}_{\beta} \Psi^{\xi}{}_{\gamma\delta} \Psi^{\pi}{}_{\eps\zeta)}
	+3\, \Psi^{\mu}{}_{(\alpha} \Psi^{\nu}{}_{\beta\gamma} \Psi^{\xi}{}_{\delta} \Psi^{\pi}{}_{\eps\zeta)}
\nonumber
\displaybreak[0]\\ &
	+6\, \Psi^{\mu}{}_{(\alpha} \Psi^{\nu}{}_{\beta\gamma} \Psi^{\xi}{}_{\delta\eps} \Psi^{\pi}{}_{\zeta)}
	+5\, \Psi^{\mu}{}_{(\alpha\beta} \Psi^{\nu}{}_{\gamma\delta} \Psi^{\xi}{}_{\eps} \Psi^{\pi}{}_{\zeta)}
	\Big\}
\nonumber
\displaybreak[0]\\ &
-6\, \Psi_{(\alpha\beta\gamma\delta\eps\zeta)}
+6\, \Psi_{\mu(\alpha} \Psi^{\mu}{}_{\beta\gamma\delta\eps\zeta)}
+15\, \Psi_{\mu(\alpha\beta} \Psi^{\mu}{}_{\gamma\delta\eps\zeta)}
+10\, \Psi_{\mu(\alpha\beta\gamma} \Psi^{\mu}{}_{\delta\eps\zeta)}
\label{eq:S6curvature}
\,.
\end{align}

\subsection{The Heat Kernel Coefficients}

By using the general algorithm for computing coincidence limits of the heat kernel coefficients and their symmetrized covariant derivatives, developed in~\cite{Avramidi1999,Avramidi2000}, we obtain
\begin{align*}
[\omega_0] &=
	1
	\,,\displaybreak[0]\\
[\omega_1] &=
	-\frac{1}{6}\, R
	\,,\displaybreak[0]\\
[\omega_{1;\alpha}] &=
	-\frac{1}{12}\, R_{;\alpha}
	\,,\displaybreak[0]\\
[\omega_{1;\alpha\beta}] &=
	-E_{\alpha\beta}
	\,,\displaybreak[0]\\
[\omega_2] &=
	\frac{1}{36}\, R^2
	+\frac{1}{15}\, R_{;\mu}{}^{\mu}
	-\frac{1}{90}\, R_{\mu\nu} R^{\mu\nu}
	+\frac{1}{90}\, R_{\mu\nu\xi\pi} R^{\mu\nu\xi\pi}
	\,,
\end{align*}
where
\begin{gather} \label{eq:E_2}
E_{\alpha\beta}=
	\frac{1}{20}\, R_{;\alpha\beta}
	+\frac{1}{60}\, R_{\alpha\beta;\mu}{}^\mu
	-\frac{1}{45}\, R_{\alpha\mu} R^\mu{}_\beta
	+\frac{1}{90}\, R_{\alpha\mu\nu\xi} R_\beta{}^{\mu\nu\xi}
	+\frac{1}{90}\, R_{\mu\nu} R^\mu{}_\alpha{}^\nu{}_\beta
	\,.
\end{gather}

So, for the coefficients $a_k(x)=\omega_k(x,\Phi(x))$ at the point $x_0$ we have
\begin{align*}
a_0(x_0) &=
	1
	\,,\displaybreak[0]\\
a_1(x_0) &=
	-\frac{1}{6}\, R
	\,,\displaybreak[0]\\
a_{1;\alpha}(x_0) &=
	-\frac{1}{12}\, R_{;\alpha}
	-\frac{1}{12}\, R_{;\mu} \Psi^\mu{}_\alpha
	\,,\displaybreak[0]\\
a_{1;\alpha\beta}(x_0) &=
	-E_{\alpha\beta}
	-\frac{1}{6}\, R_{;\mu(\alpha} \Psi^\mu{}_{\beta)}
	+2\,E_{\mu(\alpha} \Psi^\mu{}_{\beta)}
	-E_{\mu\nu} \Psi^\mu{}_{(\alpha} \Psi^\nu{}_{\beta)}
	-\frac{1}{12}\, R_{;\mu} \Psi^\mu{}_{(\alpha\beta)}
	\,,\displaybreak[0]\\
a_2(x_0) &=
	\frac{1}{36}\, R^2
	+\frac{1}{15}\, R_{;\mu}{}^{\mu}
	-\frac{1}{90}\, R_{\mu\nu} R^{\mu\nu}
	+\frac{1}{90}\, R_{\mu\nu\xi\pi} R^{\mu\nu\xi\pi}
	\,.
\end{align*}

\subsection{The Van Fleck-Morette Determinant} \label{se:VFMdet}

For the function $\zeta$ in representation~\eqref{eq:VFMzeta} we have~\cite{Avramidi2000}
\begin{align*}
[\zeta] &= 0
	\,,\\
[\zeta_\alpha] &= 0
	\,,\\
[\zeta_{\alpha\beta}] &=
	\frac{1}{6}\, R_{\alpha\beta}
	\,,\\
[\zeta_{(\alpha\beta\gamma)}] &=
	\frac{1}{4}\, R_{(\alpha\beta;\gamma)}
	\,,\\
[\zeta_{(\alpha\beta\gamma\delta)}] &=
	\frac{3}{10}\, R_{(\alpha\beta;\gamma\delta)}
	+\frac{1}{15}\, R_{ \mu (\alpha| \nu |\beta} R^\mu{}_\gamma{}^\nu{}_{\delta)}
	\,.
\end{align*}
In the same way, as for the world function, we obtain non-symmetrized derivatives of the function $\zeta$ from these expressions for symmetrized derivatives. After performing straightforward transformations, we have
\begin{align*}
[\zeta_{\alpha\beta\gamma}] &=
	\frac{1}{4}\, R_{(\alpha\beta;\gamma)}
	\,,\\
[\zeta_{\alpha\beta\gamma\delta}] &=
	F_{\alpha\beta\gamma\delta}
	\,,
\end{align*}
where
\begin{gather} \label{eq:F_4}
F_{\alpha\beta\gamma\delta}=
	\frac{3}{10}\, R_{(\alpha\beta;\gamma\delta)}
	+\frac{1}{15}\, R_{ \mu (\alpha| \nu |\beta} R^\mu{}_\gamma{}^\nu{}_{\delta)}
	-\frac{1}{9}\, R_{\mu(\alpha} R^\mu{}_{\beta|\delta|\gamma)}
	\,.
\end{gather}

By using the method developed in Section~\ref{se:DerDefDiag}, for derivatives of the function $\eta(x)=\zeta(x,\Phi(x))$ at the point $x_0$ we obtain
\begin{align*}
\eta(x_0) &= 0
	\,,\\
\eta_\alpha(x_0) &= 0
	\,,\\
\eta_{(\alpha\beta)}(x_0) &=
	\frac{1}{6}\, R_{\alpha\beta}
	-\frac{1}{3}\, R_{\mu(\alpha} \Psi^{\mu}{}_{\beta)}
	+\frac{1}{6}\, R_{\mu\nu} \Psi^{\mu}{}_{\alpha} \Psi^{\nu}{}_{\beta}
	\,,\\
\eta_{(\alpha\beta\gamma)} (x_0) &=
	-\frac{1}{2}\, R_{\mu(\alpha} \Psi^{\mu}{}_{\beta\gamma)}
	+\frac{1}{2}\, R_{\mu\nu} \Psi^{\mu}{}_{(\alpha} \Psi^{\nu}{}_{\beta\gamma)}
	+\frac{1}{4}\, R_{(\alpha\beta;\gamma)}
	+\frac{1}{4}\, R_{(\alpha\beta;|\mu|} \Psi^{\mu}{}_{\gamma)}
	\displaybreak[0]\\ &
	-\frac{1}{2}\, R_{\mu(\alpha;\beta} \Psi^{\mu}{}_{\gamma)}
	+\frac{1}{4}\, R_{\mu\nu;(\alpha} \Psi^{\mu}{}_{\beta} \Psi^{\nu}{}_{\gamma)}
	-\frac{1}{2}\, R_{\mu(\alpha;|\nu|} \Psi^{\mu}{}_{\beta} \Psi^{\nu}{}_{\gamma)}
	\displaybreak[0]\\ &
	+\frac{1}{4}\, R_{\mu\nu;\xi} \Psi^{\mu}{}_{(\alpha} \Psi^{\nu}{}_{\beta} \Psi^{\xi}{}_{\gamma)}
	\,,\\
\eta_{(\alpha\beta\gamma\delta)} (x_0) &=
	-\frac{2}{3}\, R_{\mu(\alpha} \Psi^{\mu}{}_{\beta\gamma\delta)}
	+\frac{1}{6}\, R_{\mu\nu}
	\Big\{
	4\, \Psi^{\mu}{}_{(\alpha} \Psi^{\nu}{}_{\beta\gamma\delta)}
	+3\, \Psi^{\mu}{}_{(\alpha\beta} \Psi^{\nu}{}_{\gamma\delta)}
	\Big\}
	\displaybreak[0]\\ &
	+R_{(\alpha\beta;|\mu|} \Psi^{\mu}{}_{\gamma\delta)}
	+2\, R_{\mu\nu;(\alpha} \Psi^{\mu}{}_{\beta} \Psi^{\nu}{}_{\gamma\delta)}
	+R_{\mu\nu;(\alpha\beta} \Psi^{\mu}{}_{\gamma} \Psi^{\nu}{}_{\delta)}
	\displaybreak[0]\\ &
	-\frac{1}{2}\, R_{(\alpha\beta;|\mu} \Psi^{\mu}{}_{\gamma\delta)}
	-R_{\mu(\alpha;\beta} \Psi^{\mu}{}_{\gamma\delta)}
	-R_{\mu\nu;(\alpha} \Psi^{\mu}{}_{\beta} \Psi^{\nu}{}_{\gamma\delta)}
	\displaybreak[0]\\ &
	-2\, R_{\mu(\alpha;|\nu|} \Psi^{(\mu}{}_{\beta} \Psi^{\nu)}{}_{\gamma\delta)}
	+\frac{3}{2}\, R_{(\mu\nu;\xi)} \Psi^{\mu}{}_{(\alpha} \Psi^{\nu}{}_{\beta} \Psi^{\xi}{}_{\gamma\delta)}
	+R_{(\alpha\beta;\gamma|\mu|} \Psi^{\mu}{}_{\delta)}
	\displaybreak[0]\\ &
	-R_{\mu\nu;(\alpha\beta} \Psi^{\mu}{}_{\gamma} \Psi^{\nu}{}_{\delta)}
	-2\, R_{\mu(\alpha;|\nu|\beta} \Psi^{\mu}{}_{\gamma} \Psi^{\nu}{}_{\delta)}
	+R_{\mu\nu;\xi(\alpha} \Psi^{\mu}{}_{\beta} \Psi^{\nu}{}_{\gamma} \Psi^{\xi}{}_{\delta)}
	\displaybreak[0]\\ &
	+\frac{3}{10}\, R_{(\alpha\beta;\gamma\delta)}
	+\frac{1}{15}\, R_{ \mu (\alpha| \nu |\beta} R^\mu{}_\gamma{}^\nu{}_{\delta)}
	-4\, F_{(\alpha\beta\gamma|\mu|} \Psi^{\mu}{}_{\delta)}
	\displaybreak[0]\\ &
	+6\, F_{\mu\nu(\alpha\beta} \Psi^{\mu}{}_{\gamma} \Psi^{\nu}{}_{\delta)}
	-4\, F_{\mu\nu\xi(\alpha} \Psi^{\mu}{}_{\beta} \Psi^{\nu}{}_{\gamma} \Psi^{\xi}{}_{\delta)}
	\displaybreak[0]\\ &
	+\bigg\{
	\frac{3}{10}\, R_{\mu\nu;\xi\pi}
	+\frac{1}{15}\, R_{ \rho \mu \theta \nu} R^\rho{}_\xi{}^\theta{}_{\pi}
	\bigg\}
	\Psi^{\mu}{}_{(\alpha} \Psi^{\nu}{}_{\beta} \Psi^{\xi}{}_{\gamma} \Psi^{\pi}{}_{\delta)}
	\,.
\end{align*}

\subsection{Coefficients $b_k$}

By direct expansion we get the coefficients $b_0$, $b_2$, and $b_4$, defined by~\eqref{eq:4.9} (all objects are computed at the point $x_0$ and we will omit arguments)
\begin{align*}
b_0 &=
	1
	\,,\displaybreak[0]\\
b_2 &=
	-\tld{a}_{1(0)}
	-h_2
	+\frac{1}{2}\, h_1^2
	\,,\displaybreak[0]\\
b_4 &=
	\frac{1}{2}\, \tld{a}_{2(0)}
	-\tld{a}_{1(2)}
	+h_1 \tld{a}_{1(1)}
	+(h_2-\frac{1}{2}\, h_1^2)\, \tld{a}_{1(0)} 
	\\ &
	-h_4
	+h_1 h_3
	+\frac{1}{2}\, h_2^2
	-\frac{1}{2}\, h_1^2 h_2
	+\frac{1}{24}\, h_1^4 
	\,.
\end{align*}
Due to~\eqref{eq:h_k} and taking into account that $\tld{\zeta}_{(1)}=\tld{\eta}_{(1)}=0$, we have
\begin{align*}
h_1 &= \dfrac{1}{2}\,\tld{S}_{(3)}
	\,,\displaybreak[0]\\
h_2 &= \dfrac{1}{2}\,\tld{S}_{(4)} + 2\,\tld{\zeta}_{(2)} - \tld{\eta}_{(2)}
	\,,\displaybreak[0]\\
h_3 &= \dfrac{1}{2}\,\tld{S}_{(5)} + 2\,\tld{\zeta}_{(3)} - \tld{\eta}_{(3)}
	\,,\displaybreak[0]\\
h_4 &= \dfrac{1}{2}\,\tld{S}_{(6)} + 2\,\tld{\zeta}_{(4)} - \tld{\eta}_{(4)}
	\,.
\end{align*}
Therefore,
\begin{align*}
b_0 &=
	1
	\,,\displaybreak[0]\\
b_2 &=
	-a_1
	- 2\,\tld{\zeta}_{(2)}
	+ \tld{\eta}_{(2)}
	-\dfrac{1}{2}\,\tld{S}_{(4)}
	+\frac{1}{8}\, \left(\tld{S}_{(3)}\right)^2
	\,,\displaybreak[0]\\
b_4 &=
	\frac{1}{2}\, a_2
	-\tld{a}_{1(2)}
	+ 2\, a_1 \tld{\zeta}_{(2)}
	-  a_1 \tld{\eta}_{(2)}
	+\dfrac{1}{2}\, a_1 \tld{S}_{(4)}
	- 2\,\tld{\zeta}_{(4)}
	+ \tld{\eta}_{(4)}
	+\dfrac{1}{2}\, \tld{S}_{(3)} \tld{a}_{1(1)}
	\displaybreak[0]\\ &
	+2\, \left(\tld{\zeta}_{(2)}\right)^2
	+\frac{1}{2}\, \left(\tld{\eta}_{(2)}\right)^2
	-\tld{\zeta}_{(2)} \tld{\eta}_{(2)}
	-\dfrac{1}{2}\,\tld{S}_{(6)}
	+\dfrac{1}{2}\, \tld{S}_{(4)} \tld{\zeta}_{(2)}
	-\dfrac{1}{4}\, \tld{S}_{(4)} \tld{\eta}_{(2)}
	\displaybreak[0]\\ &	
	-\frac{1}{8}\, a_1 \left(\tld{S}_{(3)}\right)^2
	+\tld{S}_{(3)}\tld{\zeta}_{(3)}
	-\frac{1}{2}\,\tld{S}_{(3)}\tld{\eta}_{(3)}
	+\dfrac{1}{4}\,\tld{S}_{(5)}\tld{S}_{(3)}
	+\dfrac{1}{8}\, \left(\tld{S}_{(4)}\right)^2
	\displaybreak[0]\\ &	
	-\frac{1}{4}\, \left(\tld{S}_{(3)}\right)^2 \tld{\zeta}_{(2)}
	+\frac{1}{8}\, \left(\tld{S}_{(3)}\right)^2 \tld{\eta}_{(2)}
	-\dfrac{1}{16}\, \tld{S}_{(4)} \left(\tld{S}_{(3)}\right)^2
	+\frac{1}{384}\, \left(\tld{S}_{(3)}\right)^4 
	\,.
\end{align*}

\subsection{Gaussian Average}

Since integrals in~\eqref{eq:A_k-curved-point} are Gaussian, it is convenient to introduce a Gaussian average for a function $f$, defined on $\R^n$, by
\begin{gather*}
\la f \ra=
	(4\pi)^{-n/2}
	\rho^{-1} 
	\int_{\R^n} dz\,
	\abs{g(x_0)}^{1/2}
	\e{-\dfrac{1}{4}\,S_{\alpha\beta} z^\alpha z^\beta}
	f(z)
	\,,
\end{gather*}
where
\begin{gather*}
\rho=\abs{\det \big(I-d\Phi(x_0)\big)}^{-1} \,.
\end{gather*}
The Gaussian average is normalized so that $\la 1 \ra=1$. Recall that $\dfrac{1}{2} S_{\alpha\beta} z^\alpha z^\beta= \tld{S}_{(2)}$, so \eqref{eq:A_k-curved-point} takes the form
\begin{gather} \label{eq:A_k-curved-point-avg}
A_k=\rho \la b_{2k} \ra
	\,.
\end{gather}
Due to Lemma~\ref{le:GaussianIntegralmd}, for polynomials in $z$ we have
\begin{align*}
\la z^{\alpha_1}\ldots z^{\alpha_{2k+1}} \ra &= 0
\,,\\
\la z^{\alpha_1}\ldots z^{\alpha_{2k}} \ra &=
\dfrac{(2k)!}{k!}\, Q^{(\alpha_1 \alpha_2}\ldots Q^{\alpha_{2k-1} \alpha_{2k})}
\,,
\end{align*}
where $Q$ is the matrix inverse to $S_{\alpha\beta}$ and it is given by~\eqref{eq:Q-curved-point}. 
Let $\tau:V_{m+2}\to V_{m}$ be the operator of contraction of a symmetric tensor with the matrix $Q$:
\begin{gather*}
(\tau A)_{\nu_1\dots\nu_{m}}=Q^{\alpha\beta} A_{\alpha\beta\nu_1\dots\nu_{m}}\,.
\end{gather*}
By using formula~\eqref{eq:A_k-curved-point-avg} and performing straightforward computations, we get \begin{align*}
A_0 &=\rho 
	\,,\displaybreak[0]\\[2ex]
A_1 &=\rho 
	\bigg[
	-a_1
	-2\, \tau \left( \zeta_{(2)} \right)
	+\tau \left( \eta_{(2)} \right)
	-\frac{1}{4}\, \tau^2 \left( S_{(4)} \right)
	+\frac{5}{12}\, \tau^3 \left( S_{(3)}\vee S_{(3)} \right)
	\bigg]
	\,,\displaybreak[0]\\[2ex]
A_2 &=\rho 
	\bigg[
	\frac{1}{2}\, a_2
	-\tau \left( a_{1(2)} \right)
	+2\, a_1 \tau \left( \zeta_{(2)} \right)
	-a_1 \tau \left( \eta_{(2)} \right)
	\displaybreak[0]\\ &
	+\dfrac{1}{4}\, a_1 \tau^2 \left( S_{(4)} \right)
	-\tau^2 \left( \zeta_{(4)} \right)
	+\frac{1}{2}\, \tau^2 \left( \eta_{(4)} \right)
	+\tau^2 \left( S_{(3)} \vee a_{1(1)} \right)
	\displaybreak[0]\\ &
	+6\, \tau^2 \left( \zeta_{(2)} \vee \zeta_{(2)} \right)
	+\frac{3}{2}\, \tau^2 \left( \eta_{(2)} \vee \eta_{(2)} \right)
	-3\, \tau^2 \left( \zeta_{(2)} \vee \eta_{(2)} \right)
	\displaybreak[0]\\ &
	-\dfrac{1}{12}\, \tau^3 \left( S_{(6)} \right)
	+\dfrac{5}{4}\, \tau^3 \left( S_{(4)} \vee \zeta_{(2)} \right)
	-\dfrac{5}{8}\, \tau^3 \left( S_{(4)} \vee \eta_{(2)} \right)
	\displaybreak[0]\\ &	
	-\frac{5}{12}\, a_1 \tau^3 \left( S_{(3)} \vee S_{(3)} \right)
	+\frac{10}{3}\, \tau^3 \left( S_{(3)} \vee \zeta_{(3)} \right)
	-\frac{5}{3}\, \tau^3 \left( S_{(3)} \vee \eta_{(3)} \right)
	\displaybreak[0]\\ &	
	+\dfrac{7}{12}\, \tau^4 \left( S_{(5)} \vee S_{(3)} \right)
	+\dfrac{35}{96}\, \tau^4 \left( S_{(4)} \vee S_{(4)} \right)
	\displaybreak[0]\\ &	
	-\frac{35}{6}\, \tau^4 \left( S_{(3)} \vee S_{(3)} \vee \zeta_{(2)} \right)
	+\frac{35}{12}\, \tau^4 \left( S_{(3)} \vee S_{(3)} \vee \eta_{(2)} \right)
	\displaybreak[0]\\ &	
	-\dfrac{35}{16}\, \tau^5 \left( S_{(4)} \vee S_{(3)} \vee S_{(3)} \right)
	\displaybreak[0]\\ &	
	+\frac{385}{288}\, \tau^6 \left( S_{(3)} \vee S_{(3)} \vee S_{(3)} \vee S_{(3)} \right)
	\bigg]
	\,.
\end{align*}

All terms in these expressions were evaluated above. Their substitution gives expressions for the computed coefficients in terms of the curvature, the differential $d\Phi$ of the mapping $\Phi$, and their covariant derivatives.

\begin{lemma} \label{le:Coefficients-curved-point}
Let conditions of Theorem~\ref{th:Expansion-curved-point} be satisfied. Then
\begin{align*}
A_0 &=\rho 
	\,,\displaybreak[0]\\[2ex]
A_1 &=\rho 
	\bigg\{
	H^{(1)}
	+H^{(1)}_{(\alpha\beta)}
	Q^{\alpha\beta}
	+H^{(1)}_{(\alpha\beta\gamma\delta)}
	Q^{\alpha\beta} Q^{\gamma\delta}
	+H^{(1)}_{(\alpha\beta\gamma\delta\eps\zeta)}
	Q^{\alpha\beta} Q^{\gamma\delta} Q^{\eps\zeta}
	\bigg\}
	\,,\displaybreak[0]\\[2ex]
A_2 &=\rho
	\bigg\{
	H^{(2)}
	+H^{(2)}_{(\alpha\beta)}
	Q^{\alpha\beta}
	+H^{(2)}_{(\alpha\beta\gamma\delta)}
	Q^{\alpha\beta} Q^{\gamma\delta}
	+H^{(2)}_{(\alpha\beta\gamma\delta\eps\zeta)}
	Q^{\alpha\beta} Q^{\gamma\delta} Q^{\eps\zeta}
	\displaybreak[0]\\ &
	+H^{(2)}_{(\alpha\beta\gamma\delta\eps\zeta\eta\theta)}
	Q^{\alpha\beta} Q^{\gamma\delta} Q^{\eps\zeta} Q^{\eta\theta}
	+H^{(2)}_{(\alpha\beta\gamma\delta\eps\zeta\eta\theta\kappa\lambda)}
	Q^{\alpha\beta} Q^{\gamma\delta} Q^{\eps\zeta} Q^{\eta\theta} Q^{\kappa\lambda}
	\displaybreak[0]\\ &	
	+H^{(2)}_{(\alpha\beta\gamma\delta\eps\zeta\eta\theta\kappa\lambda\psi\omega)}
	Q^{\alpha\beta} Q^{\gamma\delta} Q^{\eps\zeta} Q^{\eta\theta} Q^{\kappa\lambda} Q^{\psi\omega}
	\bigg\}
	,
\end{align*}
where
\begin{align*}
H^{(1)} &=
	\frac{1}{6}\, R
	\,,\displaybreak[0]\\[2ex]
H^{(1)}_{\alpha\beta} &=
	-\frac{1}{6}\, R_{\alpha\beta}
	-\frac{1}{3}\, R_{\mu\alpha} \Psi^{\mu}{}_{\beta}
	+\frac{1}{6}\, R_{\mu\nu} \Psi^{\mu}{}_{\alpha} \Psi^{\nu}{}_{\beta}
	\,,\displaybreak[0]\\[2ex]
H^{(1)}_{\alpha\beta\gamma\delta} &=
	R_{\alpha\mu\beta\nu} \Psi^{\mu}{}_{\gamma} \Psi^{\nu}{}_{\delta}
	+\Psi_{\alpha\beta\gamma\delta}
	-\frac{3}{4}\, \Psi_{\mu\alpha\beta} \Psi^{\mu}{}_{\gamma\delta}
	-\Psi_{\mu\alpha} \Psi^{\mu}{}_{\beta\gamma\delta}
	\,,\displaybreak[0]\\[2ex]
H^{(1)}_{\alpha\beta\gamma\delta\eps\zeta} &=
	\frac{15}{4}\, \Psi_{\alpha\beta\gamma} \Psi_{\delta\eps\zeta}
	-\frac{15}{2}\, \Psi_{\mu\alpha} \Psi^{\mu}{}_{\beta\gamma} \Psi_{\delta\eps\zeta}
	+\frac{15}{4}\, \Psi_{\mu\alpha} \Psi^{\mu}{}_{\beta\gamma}\Psi_{\nu\delta} \Psi^{\nu}{}_{\eps\zeta}
	\,,\displaybreak[0]\\[2ex]
H^{(2)} &=
	\frac{1}{72}\, R^2
	+\frac{1}{30}\, R_{;\mu}{}^{\mu}
	-\frac{1}{180}\, R_{\mu\nu} R^{\mu\nu}
	+\frac{1}{180}\, R_{\mu\nu\xi\pi} R^{\mu\nu\xi\pi}
	\,,\displaybreak[0]\\[2ex]
H^{(2)}_{\alpha\beta} &=
	-\frac{1}{36}\, R R_{\alpha\beta}
	+E_{\alpha\beta}
	+\bigg(
	\frac{1}{6}\, R_{;\mu\alpha}
	-\frac{1}{18}\, R R_{\mu\alpha}
	-2\,E_{\mu\alpha}
	\bigg)
	\Psi^\mu{}_{\beta}
	\displaybreak[0]\\ &
	+\bigg(
	\frac{1}{36}\, R R_{\mu\nu}
	+E_{\mu\nu}
	\bigg)
	\Psi^\mu{}_{\alpha} \Psi^\nu{}_{\beta}
	+\frac{1}{12}\, R_{;\mu} \Psi^\mu{}_{\alpha\beta}
	\,,\displaybreak[0]\\[2ex]
H^{(2)}_{\alpha\beta\gamma\delta} &=
	\frac{1}{8}\, R_{\alpha\beta} R_{\gamma\delta}
	-\frac{1}{30}\, R_{ \mu \alpha \nu \beta} R^\mu{}_\gamma{}^\nu{}_{\delta}
	-\frac{3}{20}\, R_{\alpha\beta;\gamma\delta}
	\displaybreak[0]\\ &
	+\bigg(
	\frac{1}{2}\, R_{\alpha\beta;\gamma\mu}
	-2\, F_{\alpha\beta\gamma\mu}
	\bigg)
	\Psi^{\mu}{}_{\delta}
	\displaybreak[0]\\ &
	+\bigg(
	\frac{1}{6}\, R_{\mu\alpha} R_{\nu\beta}
	+\frac{1}{6}\, R R_{\mu\alpha\nu\beta}
	-R_{\mu\alpha;\nu\beta}
	+3\, F_{\mu\nu\alpha\beta}
	\bigg)
	\Psi^{\mu}{}_{\gamma} \Psi^{\nu}{}_{\delta}
	\displaybreak[0]\\ &
	+\bigg(
	\frac{1}{2}\, R_{\mu\nu;\xi\alpha}
	-\frac{1}{6}\, R_{\mu\alpha} R_{\nu\xi}
	-2\, F_{\mu\nu\xi\alpha}
	\bigg)
	\Psi^{\mu}{}_{\beta} \Psi^{\nu}{}_{\gamma} \Psi^{\xi}{}_{\delta}
	\displaybreak[0]\\ &
	+\bigg(
	\frac{1}{24}\, R_{\mu\nu} R_{\xi\pi}
	+\frac{1}{30}\, R_{ \rho \mu \tau \nu} R^\rho{}_\xi{}^\tau{}_{\pi}
	+\frac{3}{20}\, R_{\mu\nu;\xi\pi}
	\bigg)
	\Psi^{\mu}{}_{\alpha} \Psi^{\nu}{}_{\beta} \Psi^{\xi}{}_{\gamma} \Psi^{\pi}{}_{\delta}
	\displaybreak[0]\\ &
	+\frac{1}{4}\, R_{;\alpha} \Psi_{\beta\gamma\delta} 
	+\bigg(
	\frac{1}{4}\, R_{\alpha\beta;\mu}
	-\frac{1}{2}\, R_{\mu\alpha;\beta}
	\bigg)
	\Psi^{\mu}{}_{\gamma\delta}
	-\frac{1}{4}\, R_{;\alpha} \Psi_{\mu\beta} \Psi^{\mu}{}_{\gamma\delta}
	\displaybreak[0]\\ &
	+\frac{1}{4}\, R_{;\mu} \Psi^\mu{}_\alpha \Psi_{\beta\gamma\delta}
	+\bigg(
	\frac{1}{2}\, R_{\mu\nu;\alpha}
	-\frac{1}{2}\, R_{\mu\alpha;\nu}
	-\frac{1}{2}\, R_{\nu\alpha;\mu}
	\bigg)
	\Psi^{\mu}{}_{\beta} \Psi^{\nu}{}_{\gamma\delta}
	\displaybreak[0]\\ &
	-\frac{1}{4}\, R_{;\mu} \Psi^\mu{}_\alpha \Psi_{\nu\beta} \Psi^{\nu}{}_{\gamma\delta}
	+\bigg(
	\frac{1}{4}\, R_{\mu\nu;\xi}
	+\frac{1}{2}\, R_{\xi\mu;\nu}
	\bigg) \Psi^{\mu}{}_{\alpha} \Psi^{\nu}{}_{\beta} \Psi^{\xi}{}_{\gamma\delta}
	\displaybreak[0]\\ &
	+\frac{1}{6}\, R \Psi_{\alpha\beta\gamma\delta}
	-\frac{1}{3}\, R_{\mu\alpha} \Psi^{\mu}{}_{\beta\gamma\delta}
	-R \bigg(
	\frac{1}{6}\, \Psi_{\mu\alpha} \Psi^{\mu}{}_{\beta\gamma\delta}
	+\frac{1}{8}\, \Psi_{\mu\alpha\beta} \Psi^{\mu}{}_{\gamma\delta}
	\bigg)
	\displaybreak[0]\\ &
	+R_{\mu\nu} \bigg(
	\frac{1}{3}\, \Psi^{\mu}{}_{\alpha} \Psi^{\nu}{}_{\beta\gamma\delta}
	+\frac{1}{4}\, \Psi^{\mu}{}_{\alpha\beta} \Psi^{\nu}{}_{\gamma\delta}
	\bigg)
	\,,\displaybreak[0]\\[2ex]
H^{(2)}_{\alpha\beta\gamma\delta\eps\zeta} &=
	\bigg(
	\dfrac{2}{3}\, R_{\xi\alpha\mu\beta} R^\xi{}_{\gamma\nu\delta}
	-\dfrac{5}{12}\, R_{\mu\alpha\nu\beta} R_{\gamma\delta}
	+\dfrac{1}{2}\, R_{\mu\alpha\nu\beta;\gamma\delta}
	\bigg)
	\Psi^{\mu}{}_{\eps} \Psi^{\nu}{}_{\zeta}
	\displaybreak[0]\\ &
	+\bigg(
	\dfrac{20}{9}\, R_{\pi\mu\alpha\nu} R^\pi{}_{\beta\xi\gamma}
	-\dfrac{5}{6}\, R_{\mu\alpha\nu\beta} R_{\xi\gamma}
	+\dfrac{1}{2}\, R_{\mu\alpha\nu\beta;\xi\gamma}
	\bigg)
	\Psi^{\mu}{}_{\delta} \Psi^{\nu}{}_{\eps} \Psi^{\xi}{}_{\zeta}
	\displaybreak[0]\\ &
	+\bigg(
	\dfrac{5}{12}\, R_{\mu\alpha\nu\beta} R_{\xi\pi}
	+\dfrac{2}{3}\, R_{ \rho \mu\alpha\nu} R^\rho{}_{\xi\beta\pi}
	+\dfrac{1}{2}\, R_{\mu\alpha\nu\beta;\xi\pi}
	\bigg)
	\Psi^{\mu}{}_{\gamma} \Psi^{\nu}{}_{\delta} \Psi^{\xi}{}_{\eps} \Psi^{\pi}{}_{\zeta}
	\displaybreak[0]\\ &
	+\frac{7}{4}\, R_{\alpha\beta;\gamma} \bigg(
	\Psi_{\mu\delta} \Psi^{\mu}{}_{\eps\zeta}
	-\Psi_{\delta\eps\zeta}
	\bigg)
	+\bigg(
	\frac{3}{4}\, R_{\alpha\beta;\mu}
	-\frac{3}{2}\, R_{\mu\alpha;\beta}
	\bigg)
	\Psi^{\mu}{}_{\gamma} \Psi_{\delta\eps\zeta}
	\displaybreak[0]\\ &
	+\dfrac{5}{2}\, R_{\mu\alpha\nu\beta;\gamma} \Psi^{\mu}{}_{\delta} \Psi^{\nu}{}_{\eps\zeta}
	+\bigg(
	\frac{3}{2}\, R_{\mu\alpha;\beta}
	-\frac{3}{4}\, R_{\alpha\beta;\mu}
	\bigg)
	\Psi^{\mu}{}_{\gamma} \Psi_{\nu\delta} \Psi^{\nu}{}_{\eps\zeta}
	\displaybreak[0]\\ &
	+\bigg(
	\frac{3}{4}\, R_{\mu\nu;\alpha}
	-\frac{3}{2}\, R_{\alpha\mu;\nu}
	\bigg)
	\Psi^{\mu}{}_{\beta} \Psi^{\nu}{}_{\gamma} \Psi_{\delta\eps\zeta}
	\displaybreak[0]\\ &
	+\bigg(
	\dfrac{5}{4}\, R_{\xi\mu\alpha\nu;\beta}
	+\frac{15}{4}\, R_{\xi\alpha\mu\beta;\nu}
	\bigg)
	\Psi^{\mu}{}_{\gamma} \Psi^{\nu}{}_{\delta} \Psi^{\xi}{}_{\eps\zeta}
	\displaybreak[0]\\ &
	+\bigg(
	\frac{3}{2}\, R_{\alpha\mu;\nu}
	-\frac{3}{4}\, R_{\mu\nu;\alpha}
	\bigg)
	\Psi^{\mu}{}_{\beta} \Psi^{\nu}{}_{\gamma} \Psi_{\xi\delta} \Psi^{\xi}{}_{\eps\zeta}
	\displaybreak[0]\\ &
	+\frac{3}{4}\, R_{\mu\nu;\xi}
	\Psi^{\mu}{}_{\alpha} \Psi^{\nu}{}_{\beta} \Psi^{\xi}{}_{\gamma} \Psi_{\delta\eps\zeta}
	+R_{\mu\alpha\nu\pi;\xi}
	\Psi^{\mu}{}_{\beta} \Psi^{\nu}{}_{\gamma} \Psi^{\xi}{}_{\delta} \Psi^{\pi}{}_{\eps\zeta}
	\displaybreak[0]\\ &
	-\frac{3}{4}\, R_{\mu\nu;\xi} \Psi^{\mu}{}_{\alpha} \Psi^{\nu}{}_{\beta} \Psi^{\xi}{}_{\gamma} \Psi_{\pi\delta} \Psi^{\pi}{}_{\eps\zeta}
	-\dfrac{5}{12}\, R_{\eps\zeta} \Psi_{\alpha\beta\gamma\delta}
	+\frac{5}{8}\, R \Psi_{\alpha\beta\gamma} \Psi_{\delta\eps\zeta}
	\displaybreak[0]\\ &
	+R_{\alpha\beta} \bigg(
	\dfrac{5}{16}\, \Psi_{\mu\gamma\delta} \Psi^{\mu}{}_{\eps\zeta}
	+\dfrac{5}{12}\, \Psi_{\mu\gamma} \Psi^{\mu}{}_{\delta\eps\zeta}
	\bigg)
	\displaybreak[0]\\ &	
	-R_{\mu\alpha} \bigg(
	\dfrac{5}{6}\, \Psi^{\mu}{}_{\beta} \Psi_{\gamma\delta\eps\zeta}
	+\frac{3}{2}\, \Psi^{\mu}{}_{\beta\gamma} \Psi_{\delta\eps\zeta}
	\bigg)
	\displaybreak[0]\\ &	
	+R_{\mu\alpha\nu\beta} \bigg(
	\dfrac{10}{3}\, \Psi^{\mu}{}_{\gamma} \Psi^{\nu}{}_{\delta\eps\zeta}
	+\dfrac{5}{2}\, \Psi^{\mu}{}_{\gamma\delta} \Psi^{\nu}{}_{\eps\zeta}
	\bigg)
	-\frac{5}{4}\, R \Psi_{\mu\alpha} \Psi^{\mu}{}_{\beta\gamma} \Psi_{\delta\eps\zeta}
	\displaybreak[0]\\ &
	+R_{\mu\alpha} \bigg(
	\dfrac{5}{8}\, \Psi^{\mu}{}_{\beta} \Psi_{\nu\gamma\delta} \Psi^{\nu}{}_{\eps\zeta}
	+\dfrac{5}{6}\, \Psi^{\mu}{}_{\beta} \Psi_{\nu\gamma} \Psi^{\nu}{}_{\delta\eps\zeta}
	+\dfrac{3}{2}\, \Psi^{\mu}{}_{\beta\gamma} \Psi_{\pi\delta} \Psi^{\pi}{}_{\eps\zeta}
	\bigg)
	\displaybreak[0]\\ &
	+R_{\mu\nu} \bigg(
	\dfrac{5}{12}\, \Psi^{\mu}{}_{\alpha} \Psi^{\nu}{}_{\beta} \Psi_{\gamma\delta\eps\zeta}
	+\frac{3}{2}\, \Psi^{\mu}{}_{\alpha} \Psi^{\nu}{}_{\beta\gamma} \Psi_{\delta\eps\zeta}
	\bigg)
	\displaybreak[0]\\ &	
	+R_{\mu\xi\nu\alpha} \bigg(
	\dfrac{7}{6}\,  \Psi^{\mu}{}_{\beta} \Psi^{\nu}{}_{\gamma} \Psi^{\xi}{}_{\delta\eps\zeta}
	+\dfrac{3}{2}\, \Psi^{\mu}{}_{\beta} \Psi^{\nu}{}_{\gamma\delta} \Psi^{\xi}{}_{\eps\zeta}
	\bigg)
	\displaybreak[0]\\ &
	+\frac{5}{8}\, R \Psi_{\mu\alpha} \Psi^{\mu}{}_{\beta\gamma} \Psi_{\nu\delta} \Psi^{\nu}{}_{\eps\zeta}
	+\dfrac{1}{4}\, R_{\mu\xi\nu\pi}
	\Psi^{\mu}{}_{\alpha} \Psi^{\nu}{}_{\beta} \Psi^{\xi}{}_{\gamma\delta} \Psi^{\pi}{}_{\eps\zeta}
	\displaybreak[0]\\ &
	-R_{\mu\nu} \bigg(
	\dfrac{5}{16}\, \Psi^{\mu}{}_{\alpha} \Psi^{\nu}{}_{\beta} \Psi_{\xi\gamma\delta} \Psi^{\xi}{}_{\eps\theta}
	+\dfrac{5}{12}\, \Psi^{\mu}{}_{\alpha} \Psi^{\nu}{}_{\beta} \Psi_{\xi\gamma} \Psi^{\xi}{}_{\delta\eps\theta}
	\displaybreak[0]\\ &
	+\frac{3}{2}\, \Psi^{\mu}{}_{\alpha} \Psi^{\nu}{}_{\beta\gamma} \Psi_{\xi\delta} \Psi^{\xi}{}_{\eps\zeta}
	\bigg)
	+\dfrac{1}{2}\, \Psi_{\alpha\beta\gamma\delta\eps\zeta}
	-\dfrac{1}{2}\, \Psi_{\mu\alpha} \Psi^{\mu}{}_{\beta\gamma\delta\eps\zeta}
	\displaybreak[0]\\ &
	-\dfrac{5}{4}\, \Psi_{\mu\alpha\beta} \Psi^{\mu}{}_{\gamma\delta\eps\zeta}
	-\dfrac{5}{6}\, \Psi_{\mu\alpha\beta\gamma} \Psi^{\mu}{}_{\delta\eps\zeta}
	\,,\displaybreak[0]\\[2ex]
H^{(2)}_{\alpha\beta\gamma\delta\eps\zeta\eta\theta} &=
	\dfrac{35}{6}\, R_{\alpha\mu\beta\nu} R_{\eps\xi\zeta\pi} \Psi^{\mu}{}_{\gamma} \Psi^{\nu}{}_{\delta} \Psi^{\xi}{}_{\eta} \Psi^{\pi}{}_{\theta}
	+\dfrac{35}{4}\, R_{\mu\alpha\nu\beta;\gamma} \Psi^{\mu}{}_{\delta} \Psi^{\nu}{}_{\eps} \Psi_{\zeta\eta\theta}
	\displaybreak[0]\\ &
	+\dfrac{35}{4}\, R_{\mu\alpha\nu\beta;\xi} \Psi^{\mu}{}_{\gamma} \Psi^{\nu}{}_{\delta} \Psi^{\xi}{}_{\eps} \Psi_{\zeta\eta\theta}
	-\dfrac{35}{4}\, R_{\mu\alpha\nu\beta;\gamma} \Psi^{\mu}{}_{\delta} \Psi^{\nu}{}_{\eps}  \Psi_{\pi\zeta} \Psi^{\pi}{}_{\eta\theta}
	\displaybreak[0]\\ &
	-\dfrac{35}{4}\, R_{\mu\alpha\nu\beta;\xi} \Psi^{\mu}{}_{\gamma} \Psi^{\nu}{}_{\delta} \Psi^{\xi}{}_{\eps}  \Psi_{\pi\zeta} \Psi^{\pi}{}_{\eta\theta}
	-\frac{35}{4}\, R_{\mu\alpha} \Psi^{\mu}{}_{\beta} \Psi_{\gamma\delta\eps} \Psi_{\zeta\eta\theta}
	\displaybreak[0]\\ &
	+R_{\alpha\beta} \bigg(
	\frac{35}{4}\, \Psi_{\mu\gamma} \Psi^{\mu}{}_{\delta\eps} \Psi_{\zeta\eta\theta}
	-\frac{35}{8}\, \Psi_{\gamma\delta\eps} \Psi_{\zeta\eta\theta}
	\bigg)
	\displaybreak[0]\\ &
	+R_{\mu\alpha\nu\beta} \bigg(
	\dfrac{35}{3}\, R_{\mu\alpha\nu\beta} \Psi^{\mu}{}_{\gamma} \Psi^{\nu}{}_{\delta} \Psi_{\eps\zeta\eta\theta}
	+35\, \Psi^{\mu}{}_{\gamma} \Psi^{\nu}{}_{\delta\eps} \Psi_{\zeta\eta\theta}
	\bigg)
	\displaybreak[0]\\ &	
	-\frac{35}{8}\, R_{\alpha\beta} \Psi_{\mu\gamma} \Psi^{\mu}{}_{\delta\eps} \Psi_{\nu\zeta} \Psi^{\nu}{}_{\eta\theta}
	+\frac{35}{2}\, R_{\mu\alpha} \Psi^{\mu}{}_{\beta} \Psi_{\nu\gamma} \Psi^{\nu}{}_{\delta\eps} \Psi_{\zeta\eta\theta}
	\displaybreak[0]\\ &
	+\frac{35}{8}\, R_{\mu\nu} \Psi^{\mu}{}_{\alpha} \Psi^{\nu}{}_{\beta} \Psi_{\gamma\delta\eps} \Psi_{\zeta\eta\theta}
	-R_{\mu\alpha\nu\beta} \bigg(
	35\, \Psi^{\mu}{}_{\gamma} \Psi^{\nu}{}_{\delta\eps} \Psi_{\pi\zeta} \Psi^{\pi}{}_{\eta\theta}
	\displaybreak[0]\\ &
	+\dfrac{35}{4}\, \Psi^{\mu}{}_{\gamma} \Psi^{\nu}{}_{\delta}  \Psi_{\xi\eps\zeta} \Psi^{\xi}{}_{\eta\theta}
	+\dfrac{35}{3}\, \Psi^{\mu}{}_{\gamma} \Psi^{\nu}{}_{\delta} \Psi_{\xi\eps} \Psi^{\xi}{}_{\zeta\eta\theta}
	\bigg)
	\displaybreak[0]\\ &
	+\frac{35}{4}\, R_{\mu\xi\nu\alpha} \Psi^{\mu}{}_{\beta} \Psi^{\nu}{}_{\gamma} \Psi^{\xi}{}_{\delta\eps} \Psi_{\zeta\eta\theta}
	-\frac{35}{4}\, R_{\mu\alpha} \Psi^{\mu}{}_{\beta} \Psi_{\nu\gamma} \Psi^{\nu}{}_{\delta\eps} \Psi_{\xi\zeta} \Psi^{\xi}{}_{\eta\theta}
	\displaybreak[0]\\ &
	-\frac{35}{4}\, R_{\mu\nu} \Psi^{\mu}{}_{\alpha} \Psi^{\nu}{}_{\beta} \Psi_{\xi\gamma} \Psi^{\xi}{}_{\delta\eps} \Psi_{\zeta\eta\theta}
	-\frac{35}{4}\, R_{\mu\xi\nu\alpha} \Psi^{\mu}{}_{\beta} \Psi^{\nu}{}_{\gamma} \Psi^{\xi}{}_{\delta\eps} \Psi_{\pi\zeta} \Psi^{\pi}{}_{\eta\theta}
	\displaybreak[0]\\ &
	+\frac{35}{8}\, R_{\mu\nu} \Psi^{\mu}{}_{\alpha} \Psi^{\nu}{}_{\beta} \Psi_{\xi\gamma} \Psi^{\xi}{}_{\delta\eps} \Psi_{\pi\zeta} \Psi^{\pi}{}_{\eta\theta}
	+\dfrac{35}{6}\, \Psi_{\alpha\beta\gamma\delta} \Psi_{\eps\zeta\eta\theta}
	\displaybreak[0]\\ &
	+\dfrac{35}{4}\, \Psi_{\alpha\beta\gamma\delta\eps} \Psi_{\zeta\eta\theta}
	-\dfrac{35}{2}\, \Psi_{\mu\alpha\beta} \Psi^{\mu}{}_{\gamma\delta\eps} \Psi_{\zeta\eta\theta}
	-\dfrac{35}{4}\, \Psi_{\mu\alpha} \Psi^{\mu}{}_{\beta\gamma\delta\eps} \Psi_{\zeta\eta\theta}
	\displaybreak[0]\\ &
	-\dfrac{35}{4}\, \Psi_{\alpha\beta\gamma\delta\eps} \Psi_{\pi\zeta} \Psi^{\pi}{}_{\eta\theta}
	-\dfrac{35}{4}\, \Psi_{\mu\alpha\beta} \Psi^{\mu}{}_{\gamma\delta} \Psi_{\eps\zeta\eta\theta}
	-\dfrac{35}{3}\, \Psi_{\mu\alpha} \Psi^{\mu}{}_{\beta\gamma\delta} \Psi_{\eps\zeta\eta\theta}
	\displaybreak[0]\\ &	
	+\dfrac{35}{2}\, \Psi_{\mu\alpha\beta} \Psi^{\mu}{}_{\gamma\delta\eps} \Psi_{\pi\zeta} \Psi^{\pi}{}_{\eta\theta}
	+\dfrac{35}{4}\, \Psi_{\mu\alpha} \Psi^{\mu}{}_{\beta\gamma\delta\eps} \Psi_{\pi\zeta} \Psi^{\pi}{}_{\eta\theta}
	\displaybreak[0]\\ &
	+\frac{105}{32}\, \Psi_{\mu\alpha\beta} \Psi^{\mu}{}_{\gamma\delta} \Psi_{\xi\eps\zeta} \Psi^{\xi}{}_{\eta\theta}
	+\dfrac{35}{6}\, \Psi_{\mu\alpha} \Psi^{\mu}{}_{\beta\gamma\delta} \Psi_{\xi\eps} \Psi^{\xi}{}_{\zeta\eta\theta}
	\displaybreak[0]\\ &
	+\dfrac{35}{4}\, \Psi_{\mu\alpha\beta} \Psi^{\mu}{}_{\gamma\delta} \Psi_{\xi\eps} \Psi^{\xi}{}_{\zeta\eta\theta}
	\,,\displaybreak[0]\\[2ex]
H^{(2)}_{\alpha\beta\gamma\delta\eps\zeta\eta\theta\kappa\lambda} &=
	\dfrac{315}{4}\, 
	\bigg[
	R_{\mu\alpha\nu\beta} \bigg(
	\Psi^{\mu}{}_{\gamma} \Psi^{\nu}{}_{\delta} \Psi_{\eps\zeta\eta} \Psi_{\theta\kappa\lambda}
	-2\, \Psi^{\mu}{}_{\gamma} \Psi^{\nu}{}_{\delta} \Psi_{\xi\eps} \Psi^{\xi}{}_{\zeta\eta} \Psi_{\theta\kappa\lambda}
	\displaybreak[0]\\ &	
	+\Psi^{\mu}{}_{\gamma} \Psi^{\nu}{}_{\delta} \Psi_{\xi\eps} \Psi^{\xi}{}_{\zeta\eta} \Psi_{\pi\theta} \Psi^{\pi}{}_{\kappa\lambda}
	\bigg)
	+\Psi_{\alpha\beta\gamma} \Psi_{\delta\eps\zeta} \Psi_{\eta\theta\kappa\lambda}
	\displaybreak[0]\\ &	
	-\frac{3}{4}\, \Psi_{\alpha\beta\gamma} \Psi_{\delta\eps\zeta} \Psi_{\mu\eta\theta} \Psi^{\mu}{}_{\kappa\lambda}
	-\Psi_{\alpha\beta\gamma} \Psi_{\delta\eps\zeta} \Psi_{\mu\eta} \Psi^{\mu}{}_{\theta\kappa\lambda}
	\displaybreak[0]\\ &
	-2\, \Psi_{\mu\alpha} \Psi^{\mu}{}_{\beta\gamma} \Psi_{\delta\eps\zeta} \Psi_{\eta\theta\kappa\lambda}
	+\frac{3}{2}\, \Psi_{\mu\alpha} \Psi^{\mu}{}_{\beta\gamma} \Psi_{\delta\eps\zeta} \Psi_{\nu\eta\theta} \Psi^{\nu}{}_{\kappa\lambda}
	\displaybreak[0]\\ &	
	+2\, \Psi_{\mu\alpha} \Psi^{\mu}{}_{\beta\gamma} \Psi_{\delta\eps\zeta} \Psi_{\nu\eta} \Psi^{\nu}{}_{\theta\kappa\lambda}
	+\Psi_{\mu\alpha} \Psi^{\mu}{}_{\beta\gamma} \Psi_{\nu\delta} \Psi^{\nu}{}_{\eps\zeta} \Psi_{\eta\theta\kappa\lambda}
	\displaybreak[0]\\ &	
	-\frac{3}{4}\, \Psi_{\mu\alpha} \Psi^{\mu}{}_{\beta\gamma} \Psi_{\nu\delta} \Psi^{\nu}{}_{\eps\zeta} \Psi_{\xi\eta\theta} \Psi^{\xi}{}_{\kappa\lambda}
	-\Psi_{\mu\alpha} \Psi^{\mu}{}_{\beta\gamma} \Psi_{\nu\delta} \Psi^{\nu}{}_{\eps\zeta} \Psi_{\xi\eta} \Psi^{\xi}{}_{\theta\kappa\lambda}
	\bigg]
	\,,\displaybreak[0]\\[2ex]
H^{(2)}_{\alpha\beta\gamma\delta\eps\zeta\eta\theta\kappa\lambda\psi\omega} &=
	\frac{3465}{32}\,
	\bigg[
	\Psi_{\alpha\beta\gamma} \Psi_{\delta\eps\zeta} \Psi_{\eta\theta\kappa} \Psi_{\lambda\psi\omega}
	-4\, \Psi_{\mu\alpha} \Psi^{\mu}{}_{\beta\gamma} \Psi_{\delta\eps\zeta} \Psi_{\eta\theta\kappa} \Psi_{\lambda\psi\omega}
	\displaybreak[0]\\ &	
	+6\, \Psi_{\mu\alpha} \Psi^{\mu}{}_{\beta\gamma} \Psi_{\nu\delta} \Psi^{\nu}{}_{\eps\zeta} \Psi_{\eta\theta\kappa} \Psi_{\lambda\psi\omega}
	\displaybreak[0]\\ &	
	-4\, \Psi_{\mu\alpha} \Psi^{\mu}{}_{\beta\gamma} \Psi_{\nu\delta} \Psi^{\nu}{}_{\eps\zeta} \Psi_{\xi\eta} \Psi^{\xi}{}_{\theta\kappa} \Psi_{\lambda\psi\omega}
	\displaybreak[0]\\ &	
	+\Psi_{\mu\alpha} \Psi^{\mu}{}_{\beta\gamma} \Psi_{\nu\delta} \Psi^{\nu}{}_{\eps\zeta} \Psi_{\xi\eta} \Psi^{\xi}{}_{\theta\kappa} \Psi_{\pi\lambda} \Psi^{\pi}{}_{\psi\omega}
	\bigg]
	\,,\\
\Psi^{\mu}{}_{\nu_1\nu_2\ldots\nu_k}&=(d\Phi)^{\mu'}{}_{\nu_1;\nu_2\ldots\nu_k}(x_0)
	\,,\\
\rho&=\abs{\det \big(I-d\Phi(x_0)\big)}^{-1} \,,
\end{align*}
and $E_{\alpha\beta}$, $F_{\alpha\beta\gamma\delta}$ are given by~\eqref{eq:E_2}, \eqref{eq:F_4} respectively.
\end{lemma}

\section{Isometric Case}

If the mapping $\Phi$ is an isometry, all derivatives of $d\Phi$ vanish (see, for example,~\cite{Kobayashi}) and from Lemma~\ref{le:Coefficients-curved-point} we get the following lemma.

\begin{lemma} \label{le:Coefficients-curved-point-isometry}
Let conditions of Theorem~\ref{th:Expansion-curved-point} be satisfied. Then
\begin{align*}
A_0 &=\rho 
	\,,\displaybreak[0]\\[2ex]
A_1 &=\rho 
	\bigg\{
	\frac{1}{6}\, R
	-\bigg[
	\frac{1}{6}\, R_{\alpha\beta}
	+\frac{1}{3}\, R_{\mu\alpha} \Psi^{\mu}{}_{\beta}
	-\frac{1}{6}\, R_{\mu\nu} \Psi^{\mu}{}_{\alpha} \Psi^{\nu}{}_{\beta}
	\bigg]
	Q^{\alpha\beta}
	\displaybreak[0]\\ &
	+R_{\alpha\mu\beta\nu} \Psi^{\mu}{}_{\gamma} \Psi^{\nu}{}_{\delta}
	Q^{(\alpha\beta} Q^{\gamma\delta)}
	\bigg\}
	\,,\displaybreak[0]\\[2ex]
A_2 &=\rho
	\bigg\{
	\frac{1}{72}\, R^2
	+\frac{1}{30}\, R_{;\mu}{}^{\mu}
	-\frac{1}{180}\, R_{\mu\nu} R^{\mu\nu}
	+\frac{1}{180}\, R_{\mu\nu\xi\pi} R^{\mu\nu\xi\pi}
	\displaybreak[0]\\ &
	+\bigg[
	E_{\alpha\beta}
	-2\,E_{\mu\alpha} \Psi^\mu{}_{\beta}
	+E_{\mu\nu} \Psi^\mu{}_{\alpha} \Psi^\nu{}_{\beta}
	+\frac{1}{6}\, R_{;\mu\alpha} \Psi^\mu{}_{\beta}
	\displaybreak[0]\\ &
	-\frac{1}{36}\, R R_{\alpha\beta}
	-\frac{1}{18}\, R R_{\mu\alpha} \Psi^{\mu}{}_{\beta}
	+\frac{1}{36}\, R R_{\mu\nu} \Psi^{\mu}{}_{\alpha} \Psi^{\nu}{}_{\beta}
	\bigg]
	Q^{\alpha\beta}
	\displaybreak[0]\\ &
	+\bigg[
	-\frac{3}{20}\, R_{\alpha\beta;\gamma\delta}
	+\frac{1}{2}\, R_{\alpha\beta;\gamma\mu} \Psi^{\mu}{}_{\delta}
	-R_{\mu\alpha;\nu\beta} \Psi^{\mu}{}_{\gamma} \Psi^{\nu}{}_{\delta}
	\displaybreak[0]\\ &
	+\frac{1}{2}\, R_{\mu\nu;\xi\alpha} \Psi^{\mu}{}_{\beta} \Psi^{\nu}{}_{\gamma} \Psi^{\xi}{}_{\delta}
	+\frac{3}{20}\, R_{\mu\nu;\xi\pi}
	\Psi^{\mu}{}_{\alpha} \Psi^{\nu}{}_{\beta} \Psi^{\xi}{}_{\gamma} \Psi^{\pi}{}_{\delta}
	\displaybreak[0]\\ &
	+\frac{1}{8}\, R_{\alpha\beta} R_{\gamma\delta}
	+\frac{1}{6}\, R_{\mu\alpha} R_{\xi\gamma} \Psi^{\mu}{}_{\beta} \Psi^{\xi}{}_{\delta}
	-\frac{1}{6}\, R_{\mu\alpha} R_{\xi\pi} \Psi^{\mu}{}_{\beta} \Psi^{\xi}{}_{\gamma} \Psi^{\pi}{}_{\delta}
	\displaybreak[0]\\ &
	+\frac{1}{24}\, R_{\mu\nu} R_{\xi\pi} \Psi^{\mu}{}_{\alpha} \Psi^{\nu}{}_{\beta} \Psi^{\xi}{}_{\gamma} \Psi^{\pi}{}_{\delta}
	+\frac{1}{6}\, R R_{\alpha\mu\beta\nu} \Psi^{\mu}{}_{\gamma} \Psi^{\nu}{}_{\delta}
	\displaybreak[0]\\ &
	-\frac{1}{30}\, R_{ \mu \alpha \nu \beta} R^\mu{}_\gamma{}^\nu{}_{\delta}
	+\frac{1}{30}\, R_{ \rho \mu \tau \nu} R^\rho{}_\xi{}^\tau{}_{\pi}
	\Psi^{\mu}{}_{\alpha} \Psi^{\nu}{}_{\beta} \Psi^{\xi}{}_{\gamma} \Psi^{\pi}{}_{\delta}
	\displaybreak[0]\\ &
	-2\, F_{\alpha\beta\gamma\mu} \Psi^{\mu}{}_{\delta}
	+3\, F_{\mu\nu\alpha\beta} \Psi^{\mu}{}_{\gamma} \Psi^{\nu}{}_{\delta}
	-2\, F_{\mu\nu\xi\alpha} \Psi^{\mu}{}_{\beta} \Psi^{\nu}{}_{\gamma} \Psi^{\xi}{}_{\delta}
	\bigg]
	Q^{(\alpha\beta} Q^{\gamma\delta)}
	\displaybreak[0]\\ &
	+
	\bigg[
	-\dfrac{5}{12}\, R_{\alpha\mu\beta\nu} R_{\eps\zeta} \Psi^{\mu}{}_{\gamma} \Psi^{\nu}{}_{\delta}
	-\dfrac{5}{6}\, R_{\alpha\mu\beta\nu} R_{\xi\eps} \Psi^{\xi}{}_{\zeta} \Psi^{\mu}{}_{\gamma} \Psi^{\nu}{}_{\delta}
	\displaybreak[0]\\ &
	+\dfrac{5}{12}\, R_{\alpha\mu\beta\nu} R_{\xi\pi} \Psi^{\xi}{}_{\eps} \Psi^{\pi}{}_{\zeta} \Psi^{\mu}{}_{\gamma} \Psi^{\nu}{}_{\delta}
	+\dfrac{2}{3}\, R_{\xi\alpha\mu\beta} R^\xi{}_{\gamma\nu\delta} \Psi^{\mu}{}_{\eps} \Psi^{\nu}{}_{\zeta}
	\displaybreak[0]\\ &
	+\dfrac{20}{9}\, R_{\pi\mu\alpha\nu} R^\pi{}_{\beta\xi\gamma} \Psi^{\mu}{}_{\delta} \Psi^{\nu}{}_{\eps} \Psi^{\xi}{}_{\zeta}
	+\dfrac{2}{3}\, R_{ \rho \mu\nu\alpha} R^\rho{}_{\xi\pi\beta} \Psi^{\mu}{}_{\gamma} \Psi^{\nu}{}_{\delta} \Psi^{\xi}{}_{\eps} \Psi^{\pi}{}_{\zeta}
	\displaybreak[0]\\ &
	+\dfrac{1}{2}\, R_{\mu\alpha\nu\beta;\gamma\delta} \Psi^{\mu}{}_{\eps} \Psi^{\nu}{}_{\zeta}
	+\dfrac{1}{2}\, R_{\mu\alpha\nu\beta;\xi\gamma} \Psi^{\mu}{}_{\delta} \Psi^{\nu}{}_{\eps} \Psi^{\xi}{}_{\zeta}
	\displaybreak[0]\\ &
	+\dfrac{1}{2}\, R_{\mu\alpha\nu\beta;\xi\pi} \Psi^{\mu}{}_{\gamma} \Psi^{\nu}{}_{\delta} \Psi^{\xi}{}_{\eps} \Psi^{\pi}{}_{\zeta}
	\bigg]
	Q^{(\alpha\beta} Q^{\gamma\delta} Q^{\eps\zeta)}
	\displaybreak[0]\\ &	
	+
	\dfrac{35}{6}\, R_{\alpha\mu\beta\nu} R_{\eps\xi\zeta\pi} \Psi^{\mu}{}_{\gamma} \Psi^{\nu}{}_{\delta} \Psi^{\xi}{}_{\eta} \Psi^{\pi}{}_{\theta}
	Q^{(\alpha\beta} Q^{\gamma\delta} Q^{\eps\zeta} Q^{\eta\theta)}
	\bigg\}
	,\\
\Psi^{\mu}{}_{\nu}&=(d\Phi)^{\mu'}{}_{\nu}(x_0)
	\,,\\
\rho&=\abs{\det \big(I-d\Phi(x_0)\big)}^{-1} \,,
\end{align*}
and $E_{\alpha\beta}$, $F_{\alpha\beta\gamma\delta}$ are given by~\eqref{eq:E_2}, \eqref{eq:F_4} respectively.
\end{lemma}

\section{Comparison of the Results}

In this section we compare our expression for the coefficient $A_1$ in a particular case when the mapping $\Phi$ is an isometry
\begin{align}
A_1 &=\rho 
	\bigg\{
	\frac{1}{6}\, R
	-\bigg[
	\frac{1}{6}\, R_{\alpha\beta}
	+\frac{1}{3}\, R_{\mu\alpha} \Psi^{\mu}{}_{\beta}
	-\frac{1}{6}\, R_{\mu\nu} \Psi^{\mu}{}_{\alpha} \Psi^{\nu}{}_{\beta}
	\bigg]
	Q^{\alpha\beta}
	\nonumber
	\displaybreak[0]\\ &
	+R_{\alpha\mu\beta\nu} \Psi^{\mu}{}_{\gamma} \Psi^{\nu}{}_{\delta}
	Q^{(\alpha\beta} Q^{\gamma\delta)}
	\bigg\}
	\label{eq:A_1-isometry}
\end{align}
with the result obtained by Donnelly in~\cite{Donnelly}:
\begin{gather} \label{eq:A_1-Donnelly}
A_1=\abs{\det B} \bigg[
	\frac{1}{3}\, R
	+\frac{2}{3}\, R_{(\alpha|\mu|\beta)\nu} B^{\mu\alpha} B^{\nu\beta}
	\bigg] \,,
\end{gather}
where
\begin{gather*}
B=(I-d\Phi)^{-1}
\end{gather*}
and all tensors are evaluated at the point $x_0$.

From~\eqref{eq:Q-curved-point} it follows that 
\begin{gather*}
Q^{\alpha\beta}=B^\alpha{}_\mu B^{\beta\mu} \,.
\end{gather*}
Since $d\Phi(x_0)=\Psi$ and $(I-d\Phi)B=I$, we have
\begin{gather*}
\Psi^\alpha{}_\mu B^\mu{}_\beta=B^\alpha{}_\beta-\delta^\alpha{}_\beta
\,,
\end{gather*}
and~\eqref{eq:A_1-isometry} takes the form
\begin{align*}
A_1 &=\abs{\det B} 
	\bigg\{
	\frac{1}{6}\, R
	-\frac{1}{6}\, R_{\alpha\beta} B^\alpha{}_\mu B^{\beta \mu}
	-\frac{1}{3}\, R_{\mu\alpha} (B^\mu{}_\nu-\delta^\mu{}_\nu) B^{\alpha \nu}
	\displaybreak[0]\\ &
	+\frac{1}{6}\, R_{\mu\nu} (B^\mu{}_\xi B^{\nu \xi}-2\,B^{\mu\nu}+\delta^{\mu\nu})
	\displaybreak[0]\\ &
	+R_{\alpha\mu\beta\nu} 
	B^{(\alpha}{}_\xi B^{\beta |\xi|} \big(B^\mu{}_\pi-\delta^\mu{}_\pi \big) \big(B^{\nu)\pi}-\delta^{\nu)\pi}\big) 
	\bigg\} .
\end{align*}
After performing straightforward computations, we get the same expression as~\eqref{eq:A_1-Donnelly}.

\chapter{Conclusions}

Let us summarize the results obtained in this thesis. Our goal was to study the geometry of the fixed point set $\Sigma$ of a smooth mapping $\Phi$ on a Riemannian manifold $(M,g)$ by computing the asymptotic expansion as $t\to 0^+$ of the trace of the deformed heat kernel $U\big(t;x,\Phi(x)\big)$ of the Laplace operator on $M$. The fixed point set $\Sigma$ was assumed to be a smooth compact submanifold of $M$.

Multi-dimensional Gaussian integrals are expressed in terms of full contractions of the symmetrized products of symmetric tensors with a symmetric (2,0)-tensor. In Chapter~\ref{ch:SymmetricTensors} we developed a technique for obtaining symmetrization-free formulas for such contracted products. The described algorithm can be easily implemented for automatic computations.

In Chapter~\ref{ch:FlatManifolds} we considered the case of a flat manifold, $M=\R^2$. We showed the existence of the asymptotic expansion in the form
\begin{gather*}
\int_M \dvol(x)\, U\big(t;x,\Phi(x)\big) \sim
\sum_{k=0}^\infty t^{k-m/2} \int_\Sigma \dvol(y)\, a_k(y) \,,
\end{gather*}
where $m=\dim\Sigma$ and $a_k(y)$ are scalar invariants depending polynomially on derivatives of the function $S(x)=\sigma(x,\Phi(x))$, and the inverse matrix of its Hessian matrix.
We developed a generalized Laplace method for computing the coefficients $a_k$ and computed explicitly the coefficients $a_0$, $a_1$, and $a_2$ for both zero- and one-dimensional fixed point sets $\Sigma$. We used the results from Chapter~\ref{ch:SymmetricTensors} to obtain a symmetrization-free form in the zero-dimensional case.

Finally, in Chapter~\ref{ch:CurvedManifolds} we obtained the asymptotic expansion as $t\to 0^+$ for an arbitrary curved manifold $M$ and a zero-dimensional fixed point set $\Sigma=\set{x_0}$ in the form
\begin{gather*}
\int_M \dvol(x)\, U\big(t;x,\Phi(x)\big) \sim
\sum_{k=0}^\infty t^k A_k \,,
\end{gather*}
where $A_k$ are scalar invariants depending polynomially on covariant derivatives of the curvature of the metric $g$, symmetrized covariant derivatives of the differential $d\Phi$ of the mapping $\Phi$, and the matrix
\begin{gather*}
Q= \left\{ (I-d\Phi)^T g (I-d\Phi) \right\}^{-1} ,
\end{gather*}
evaluated at the point $x_0$.
The coefficients $A_k$ are expressed in terms of symmetrized covariant derivatives of functions of the form $f(x,\Phi(x)$. We developed an algorithm for computation of such symmetrized covariant derivatives and realized it as the Python script listed in Appendix~\ref{ap:Script}. We computed explicitly the coefficients $A_0$, $A_1$, and $A_2$. In a particular case, when the mapping $\Phi$ is an isometry, our result for the coefficient $A_1$ coincides with the former result of Donnelly~\cite{Donnelly}.

\appendix

\chapter{Computation of Symmetrized Covariant Derivatives}

\section{Generating Python Script: Source Code}  \label{ap:Script}

\small
\begin{verbatim}
#####################################################################
# Computes symmetrized covariant derivatives of the function
#   S(x)=\sigma(x,\Phi(x))
# on the fixed point set of the mapping \Phi(x).
#
# Input:
#   sSeparator - a separator which will be inserted between all terms
#        in the output
#   nDerivativeOrder - the order of the derivative to be computed
#
#   Uses indices from lsAlphabet for derivatives and from
#   lsPrimeAlphabet for dummy indices in contractions
#
# Creates in the current directory file named "output.tex" which
# contains symmetrized derivatives of the function S(x) up to
# nDerivativeOrder written in LaTeX notation.
#
# Properties of the world functions, taken into account:
#   [\sigma]=0
#   [\sigma_\alpha]=0
#   [\sigma_{\alpha\beta\gamma}]=0
#   [\sigma_{(\alpha_1\dots\alpha_n)}]=0, n>=3
#   [\sigma_{\mu(\alpha_1\dots\alpha_n)}]=0, n>=2
#
# An existing file with the same name will be overwritten.
#
# Copyright (c) Andrey Novoseltsev, 2005
#####################################################################

import os
from array import array

lsAlphabet=[
r"\alpha",
r"\beta",
r"\gamma",
r"\delta",
r"\eps",
r"\zeta",
r"\eta",
r"\theta",
r"\kappa",
r"\lambda"]

lsPrimeAlphabet=[
r"\mu",
r"\nu",
r"\xi",
r"\pi",
r"\rho",
r"\tau",
r"\phi",
r"\chi",
r"\psi",
r"\omega"]

sSigma=r"\sigma"
sPsi=r"\Psi"

sSeparator="\n"+r"\displaybreak[0]\\ &"+"\n"

nDerivativeOrder=6

def NextDistribution(n,m,k,aI):
    if aI[0]==n-m-k+1:
        return False
    aI0old=aI[0]
    aI[0]=1
    for i in range(1,k):
        if aI[i]!=1:
            aI[i]=aI[i]-1
            aI[i-1]=aI0old+1
            return True



# Computing terms coefficients
dCoefficients={}
dCoefficients[0,0,()]=1

for n in range(1,nDerivativeOrder+1):
    for m in range(n):
        for k in range(1,n-m+1):
            aI=array('H')
            for i in range(1,k):
                aI.append(1)
            aI.append(n-m-k+1)

            while 1==1:
                #Compute the current coefficient
                coef=0
                if m!=0:
                    coef=coef+dCoefficients[m-1,k,tuple(aI)]
                if aI[k-1]==1:
                    coef=coef+dCoefficients[m,k-1,tuple(aI)[0:k-1]]
                for i in range(k):
                    if aI[i]>1:
                        aI[i]=aI[i]-1
                        coef=coef+dCoefficients[m,k,tuple(aI)]
                        aI[i]=aI[i]+1
                dCoefficients[m,k,tuple(aI)]=coef
                #Proceed to the next coefficient
                if not NextDistribution(n,m,k,aI): break
    dCoefficients[n,0,()]=1



# Constructing LaTeX expressions
fOutput=file('output.tex','w+')

fOutput.write(
 "S=0"
+"\n%%%%%%%%%%%%%%%%%%%%%%%%%%%%%%%%%%%%%%%%%%%%\n"
+"S_{"+lsAlphabet[0]+"}=0"
+"\n%%%%%%%%%%%%%%%%%%%%%%%%%%%%%%%%%%%%%%%%%%%%\n"
+"S_{("+lsAlphabet[0]+lsAlphabet[1]+")}="
+"g_{"+lsAlphabet[0]+lsAlphabet[1]+"}"
+"-2\, "+sPsi+"_{("+lsAlphabet[0]+lsAlphabet[1]+")}"
+"+"+sPsi+"_{"+lsPrimeAlphabet[0]+"("+lsAlphabet[0]+"} "
+sPsi+"^{"+lsPrimeAlphabet[0]+"}{}_{"+lsAlphabet[1]+")}")

for n in range(3,nDerivativeOrder+1):
    fOutput.write("\n%%%%%%%%%%%%%%%%%%%%%%%%%%%%%%%%%%%%%%%%%%%%\n")
    fOutput.write("S_{(")
    for i in range(n): fOutput.write(lsAlphabet[i])
    fOutput.write(")}=\n")

    nNumberOfTerms=0
    for t in range(2,n+1):
        for m in range(t-1,-1,-1):
            k=t-m
            aI=array('H')
            for i in range(1,k):
                aI.append(1)
            aI.append(n-m-k+1)

            while 1==1:
                #Check if the current term should be written
                if t==3 or (k==1 and m>=2) or (k>=n-1 and k>=3):
                    if not NextDistribution(n,m,k,aI):
                        break
                    else:
                        continue
                #Write the current term
                if nNumberOfTerms!=0: fOutput.write(sSeparator)
                if k==1 and m==1:
                    fOutput.write("-")
                elif nNumberOfTerms!=0:
                    fOutput.write("+")
                nNumberOfTerms=nNumberOfTerms+1
                if dCoefficients[m,k,tuple(aI)]!=1:
                    fOutput.write(str(dCoefficients[m,k,tuple(aI)])
                      +"\, ")
                if t==2:
                    fOutput.write("g_{")
                    for i in range(k):
                        fOutput.write(lsPrimeAlphabet[i])
                    if m!=0: fOutput.write("(")
                    for i in range(m): fOutput.write(lsAlphabet[i])
                    fOutput.write("} ")
                else:
                    fOutput.write("["+sSigma+"_{")
                    for i in range(k):
                        fOutput.write(lsPrimeAlphabet[i]+"'")
                    if m!=0: fOutput.write("(")
                    for i in range(m): fOutput.write(lsAlphabet[i])
                    fOutput.write("}] ")
                fOutput.write(sPsi+"^{"+lsPrimeAlphabet[0]+"}{}_{")
                if m==0: fOutput.write("(")
                for i in range(m,m+aI[0]):
                    fOutput.write(lsAlphabet[i])
                nextIndex=m+aI[0]
                for i in range(1,k):
                    fOutput.write("} "+sPsi+"^{"+lsPrimeAlphabet[i]
                      +"}{}_{")
                    for j in range(nextIndex,nextIndex+aI[i]):
                        fOutput.write(lsAlphabet[j])
                    nextIndex=nextIndex+aI[i]
                fOutput.write(")}")
                #Proceed to the next coefficient
                if not NextDistribution(n,m,k,aI): break
    fOutput.write("\n% Number of Terms: "+str(nNumberOfTerms))

fOutput.close()
\end{verbatim}
\normalsize

\section{Generating Python Script: Output}

\vspace{-5ex}
\begin{align*}
S&=0
\,,\\
S_{\alpha}&=0
\,,\\
S_{(\alpha\beta)}&=g_{\alpha\beta}-2\, \Psi_{(\alpha\beta)}+\Psi_{\mu(\alpha} \Psi^{\mu}{}_{\beta)}
\,,\\
S_{(\alpha\beta\gamma)}&=
-3\, g_{\mu(\alpha} \Psi^{\mu}{}_{\beta\gamma)}
\displaybreak[0]\\ &
+g_{\mu\nu} \Psi^{\mu}{}_{(\alpha} \Psi^{\nu}{}_{\beta\gamma)}
\displaybreak[0]\\ &
+2\, g_{\mu\nu} \Psi^{\mu}{}_{(\alpha\beta} \Psi^{\nu}{}_{\gamma)}
\,,\\
S_{(\alpha\beta\gamma\delta)}&=
-4\, g_{\mu(\alpha} \Psi^{\mu}{}_{\beta\gamma\delta)}
\displaybreak[0]\\ &
+g_{\mu\nu} \Psi^{\mu}{}_{(\alpha} \Psi^{\nu}{}_{\beta\gamma\delta)}
\displaybreak[0]\\ &
+3\, g_{\mu\nu} \Psi^{\mu}{}_{(\alpha\beta} \Psi^{\nu}{}_{\gamma\delta)}
\displaybreak[0]\\ &
+3\, g_{\mu\nu} \Psi^{\mu}{}_{(\alpha\beta\gamma} \Psi^{\nu}{}_{\delta)}
\displaybreak[0]\\ &
+6\, [\sigma_{\mu'\nu'(\alpha\beta}] \Psi^{\mu}{}_{\gamma} \Psi^{\nu}{}_{\delta)}
\,,\\
S_{(\alpha\beta\gamma\delta\eps)}&=
-5\, g_{\mu(\alpha} \Psi^{\mu}{}_{\beta\gamma\delta\eps)}
\displaybreak[0]\\ &
+g_{\mu\nu} \Psi^{\mu}{}_{(\alpha} \Psi^{\nu}{}_{\beta\gamma\delta\eps)}
\displaybreak[0]\\ &
+4\, g_{\mu\nu} \Psi^{\mu}{}_{(\alpha\beta} \Psi^{\nu}{}_{\gamma\delta\eps)}
\displaybreak[0]\\ &
+6\, g_{\mu\nu} \Psi^{\mu}{}_{(\alpha\beta\gamma} \Psi^{\nu}{}_{\delta\eps)}
\displaybreak[0]\\ &
+4\, g_{\mu\nu} \Psi^{\mu}{}_{(\alpha\beta\gamma\delta} \Psi^{\nu}{}_{\eps)}
\displaybreak[0]\\ &
+10\, [\sigma_{\mu'\nu'(\alpha\beta}] \Psi^{\mu}{}_{\gamma} \Psi^{\nu}{}_{\delta\eps)}
\displaybreak[0]\\ &
+20\, [\sigma_{\mu'\nu'(\alpha\beta}] \Psi^{\mu}{}_{\gamma\delta} \Psi^{\nu}{}_{\eps)}
\displaybreak[0]\\ &
+5\, [\sigma_{\mu'\nu'\xi'(\alpha}] \Psi^{\mu}{}_{\beta} \Psi^{\nu}{}_{\gamma} \Psi^{\xi}{}_{\delta\eps)}
\displaybreak[0]\\ &
+10\, [\sigma_{\mu'\nu'\xi'(\alpha}] \Psi^{\mu}{}_{\beta} \Psi^{\nu}{}_{\gamma\delta} \Psi^{\xi}{}_{\eps)}
\displaybreak[0]\\ &
+15\, [\sigma_{\mu'\nu'\xi'(\alpha}] \Psi^{\mu}{}_{\beta\gamma} \Psi^{\nu}{}_{\delta} \Psi^{\xi}{}_{\eps)}
\displaybreak[0]\\ &
+10\, [\sigma_{\mu'\nu'(\alpha\beta\gamma}] \Psi^{\mu}{}_{\delta} \Psi^{\nu}{}_{\eps)}
\displaybreak[0]\\ &
+10\, [\sigma_{\mu'\nu'\xi'(\alpha\beta}] \Psi^{\mu}{}_{\gamma} \Psi^{\nu}{}_{\delta} \Psi^{\xi}{}_{\eps)}
\,,\\
S_{(\alpha\beta\gamma\delta\eps\zeta)}&=
-6\, g_{\mu(\alpha} \Psi^{\mu}{}_{\beta\gamma\delta\eps\zeta)}
\displaybreak[0]\\ &
+g_{\mu\nu} \Psi^{\mu}{}_{(\alpha} \Psi^{\nu}{}_{\beta\gamma\delta\eps\zeta)}
\displaybreak[0]\\ &
+5\, g_{\mu\nu} \Psi^{\mu}{}_{(\alpha\beta} \Psi^{\nu}{}_{\gamma\delta\eps\zeta)}
\displaybreak[0]\\ &
+10\, g_{\mu\nu} \Psi^{\mu}{}_{(\alpha\beta\gamma} \Psi^{\nu}{}_{\delta\eps\zeta)}
\displaybreak[0]\\ &
+10\, g_{\mu\nu} \Psi^{\mu}{}_{(\alpha\beta\gamma\delta} \Psi^{\nu}{}_{\eps\zeta)}
\displaybreak[0]\\ &
+5\, g_{\mu\nu} \Psi^{\mu}{}_{(\alpha\beta\gamma\delta\eps} \Psi^{\nu}{}_{\zeta)}
\displaybreak[0]\\ &
+15\, [\sigma_{\mu'\nu'(\alpha\beta}] \Psi^{\mu}{}_{\gamma} \Psi^{\nu}{}_{\delta\eps\zeta)}
\displaybreak[0]\\ &
+45\, [\sigma_{\mu'\nu'(\alpha\beta}] \Psi^{\mu}{}_{\gamma\delta} \Psi^{\nu}{}_{\eps\zeta)}
\displaybreak[0]\\ &
+45\, [\sigma_{\mu'\nu'(\alpha\beta}] \Psi^{\mu}{}_{\gamma\delta\eps} \Psi^{\nu}{}_{\zeta)}
\displaybreak[0]\\ &
+6\, [\sigma_{\mu'\nu'\xi'(\alpha}] \Psi^{\mu}{}_{\beta} \Psi^{\nu}{}_{\gamma} \Psi^{\xi}{}_{\delta\eps\zeta)}
\displaybreak[0]\\ &
+18\, [\sigma_{\mu'\nu'\xi'(\alpha}] \Psi^{\mu}{}_{\beta} \Psi^{\nu}{}_{\gamma\delta} \Psi^{\xi}{}_{\eps\zeta)}
\displaybreak[0]\\ &
+24\, [\sigma_{\mu'\nu'\xi'(\alpha}] \Psi^{\mu}{}_{\beta\gamma} \Psi^{\nu}{}_{\delta} \Psi^{\xi}{}_{\eps\zeta)}
\displaybreak[0]\\ &
+18\, [\sigma_{\mu'\nu'\xi'(\alpha}] \Psi^{\mu}{}_{\beta} \Psi^{\nu}{}_{\gamma\delta\eps} \Psi^{\xi}{}_{\zeta)}
\displaybreak[0]\\ &
+48\, [\sigma_{\mu'\nu'\xi'(\alpha}] \Psi^{\mu}{}_{\beta\gamma} \Psi^{\nu}{}_{\delta\eps} \Psi^{\xi}{}_{\zeta)}
\displaybreak[0]\\ &
+36\, [\sigma_{\mu'\nu'\xi'(\alpha}] \Psi^{\mu}{}_{\beta\gamma\delta} \Psi^{\nu}{}_{\eps} \Psi^{\xi}{}_{\zeta)}
\displaybreak[0]\\ &
+[\sigma_{\mu'\nu'\xi'\pi'}] \Psi^{\mu}{}_{(\alpha} \Psi^{\nu}{}_{\beta} \Psi^{\xi}{}_{\gamma} \Psi^{\pi}{}_{\delta\eps\zeta)}
\displaybreak[0]\\ &
+3\, [\sigma_{\mu'\nu'\xi'\pi'}] \Psi^{\mu}{}_{(\alpha} \Psi^{\nu}{}_{\beta} \Psi^{\xi}{}_{\gamma\delta} \Psi^{\pi}{}_{\eps\zeta)}
\displaybreak[0]\\ &
+4\, [\sigma_{\mu'\nu'\xi'\pi'}] \Psi^{\mu}{}_{(\alpha} \Psi^{\nu}{}_{\beta\gamma} \Psi^{\xi}{}_{\delta} \Psi^{\pi}{}_{\eps\zeta)}
\displaybreak[0]\\ &
+5\, [\sigma_{\mu'\nu'\xi'\pi'}] \Psi^{\mu}{}_{(\alpha\beta} \Psi^{\nu}{}_{\gamma} \Psi^{\xi}{}_{\delta} \Psi^{\pi}{}_{\eps\zeta)}
\displaybreak[0]\\ &
+3\, [\sigma_{\mu'\nu'\xi'\pi'}] \Psi^{\mu}{}_{(\alpha} \Psi^{\nu}{}_{\beta} \Psi^{\xi}{}_{\gamma\delta\eps} \Psi^{\pi}{}_{\zeta)}
\displaybreak[0]\\ &
+8\, [\sigma_{\mu'\nu'\xi'\pi'}] \Psi^{\mu}{}_{(\alpha} \Psi^{\nu}{}_{\beta\gamma} \Psi^{\xi}{}_{\delta\eps} \Psi^{\pi}{}_{\zeta)}
\displaybreak[0]\\ &
+10\, [\sigma_{\mu'\nu'\xi'\pi'}] \Psi^{\mu}{}_{(\alpha\beta} \Psi^{\nu}{}_{\gamma} \Psi^{\xi}{}_{\delta\eps} \Psi^{\pi}{}_{\zeta)}
\displaybreak[0]\\ &
+6\, [\sigma_{\mu'\nu'\xi'\pi'}] \Psi^{\mu}{}_{(\alpha} \Psi^{\nu}{}_{\beta\gamma\delta} \Psi^{\xi}{}_{\eps} \Psi^{\pi}{}_{\zeta)}
\displaybreak[0]\\ &
+15\, [\sigma_{\mu'\nu'\xi'\pi'}] \Psi^{\mu}{}_{(\alpha\beta} \Psi^{\nu}{}_{\gamma\delta} \Psi^{\xi}{}_{\eps} \Psi^{\pi}{}_{\zeta)}
\displaybreak[0]\\ &
+10\, [\sigma_{\mu'\nu'\xi'\pi'}] \Psi^{\mu}{}_{(\alpha\beta\gamma} \Psi^{\nu}{}_{\delta} \Psi^{\xi}{}_{\eps} \Psi^{\pi}{}_{\zeta)}
\displaybreak[0]\\ &
+20\, [\sigma_{\mu'\nu'(\alpha\beta\gamma}] \Psi^{\mu}{}_{\delta} \Psi^{\nu}{}_{\eps\zeta)}
\displaybreak[0]\\ &
+40\, [\sigma_{\mu'\nu'(\alpha\beta\gamma}] \Psi^{\mu}{}_{\delta\eps} \Psi^{\nu}{}_{\zeta)}
\displaybreak[0]\\ &
+15\, [\sigma_{\mu'\nu'\xi'(\alpha\beta}] \Psi^{\mu}{}_{\gamma} \Psi^{\nu}{}_{\delta} \Psi^{\xi}{}_{\eps\zeta)}
\displaybreak[0]\\ &
+30\, [\sigma_{\mu'\nu'\xi'(\alpha\beta}] \Psi^{\mu}{}_{\gamma} \Psi^{\nu}{}_{\delta\eps} \Psi^{\xi}{}_{\zeta)}
\displaybreak[0]\\ &
+45\, [\sigma_{\mu'\nu'\xi'(\alpha\beta}] \Psi^{\mu}{}_{\gamma\delta} \Psi^{\nu}{}_{\eps} \Psi^{\xi}{}_{\zeta)}
\displaybreak[0]\\ &
+6\, [\sigma_{\mu'\nu'\xi'\pi'(\alpha}] \Psi^{\mu}{}_{\beta} \Psi^{\nu}{}_{\gamma} \Psi^{\xi}{}_{\delta} \Psi^{\pi}{}_{\eps\zeta)}
\displaybreak[0]\\ &
+12\, [\sigma_{\mu'\nu'\xi'\pi'(\alpha}] \Psi^{\mu}{}_{\beta} \Psi^{\nu}{}_{\gamma} \Psi^{\xi}{}_{\delta\eps} \Psi^{\pi}{}_{\zeta)}
\displaybreak[0]\\ &
+18\, [\sigma_{\mu'\nu'\xi'\pi'(\alpha}] \Psi^{\mu}{}_{\beta} \Psi^{\nu}{}_{\gamma\delta} \Psi^{\xi}{}_{\eps} \Psi^{\pi}{}_{\zeta)}
\displaybreak[0]\\ &
+24\, [\sigma_{\mu'\nu'\xi'\pi'(\alpha}] \Psi^{\mu}{}_{\beta\gamma} \Psi^{\nu}{}_{\delta} \Psi^{\xi}{}_{\eps} \Psi^{\pi}{}_{\zeta)}
\displaybreak[0]\\ &
+15\, [\sigma_{\mu'\nu'(\alpha\beta\gamma\delta}] \Psi^{\mu}{}_{\eps} \Psi^{\nu}{}_{\zeta)}
\displaybreak[0]\\ &
+20\, [\sigma_{\mu'\nu'\xi'(\alpha\beta\gamma}] \Psi^{\mu}{}_{\delta} \Psi^{\nu}{}_{\eps} \Psi^{\xi}{}_{\zeta)}
\displaybreak[0]\\ &
+15\, [\sigma_{\mu'\nu'\xi'\pi'(\alpha\beta}] \Psi^{\mu}{}_{\gamma} \Psi^{\nu}{}_{\delta} \Psi^{\xi}{}_{\eps} \Psi^{\pi}{}_{\zeta)}
\,.
\end{align*}

\section{Coincidence Limits of the World Function} \label{ap:CLofWF}

\subsection{Forth order}

\vspace{-5ex}
\begin{align*}
[\sigma_{\mu'\nu'\xi'\alpha}] &=
	[\sigma_{\alpha\mu'\nu'\xi'}]
	\,,\\
[\sigma_{\alpha\mu'(\nu'\xi')}] &=
	\frac{1}{2}( [\sigma_{\alpha\mu'\nu'\xi'}]+[\sigma_{\alpha\mu'\xi'\nu'}] )
	=[\sigma_{\alpha\mu'\nu'\xi'}]-\frac{1}{2}[\sigma_{\alpha a'} R^{a'}{}_{\mu'\nu'\xi'}]
	\,,\\
[\sigma_{\alpha\mu'\nu'\xi'}] &=
	[\sigma_{\alpha\mu'(\nu'\xi')}] - \frac{1}{2}\, R_{\alpha\mu\nu\xi} 
	= -\frac{1}{3}\, R_{\alpha(\nu|\mu|\xi)} - \frac{1}{2}\, R_{\alpha\mu\nu\xi}
	\\ &
	= -\frac{1}{6}\, (R_{\alpha\nu\mu\xi}+R_{\alpha\xi\mu\nu}) - \frac{1}{2}\, R_{\alpha\mu\nu\xi}
	\\ &
	= \frac{1}{6}\, (R_{\alpha\nu\xi\mu}+R_{\alpha\mu\nu\xi}+R_{\alpha\nu\xi\mu}) + \frac{1}{2}\, R_{\alpha\mu\xi\nu}
	\\ &
	= \frac{1}{3}\, (R_{\alpha\nu\xi\mu}+R_{\alpha\mu\xi\nu})
	=\frac{2}{3}\, R_{\alpha(\mu|\xi|\nu)}
	\,,\\
[\sigma_{\mu'\nu'\alpha\beta}] &=
	-[\sigma_{\mu'\nu'\beta'\alpha}] =
	-\frac{2}{3}\, R_{\alpha(\mu|\beta|\nu)}
	\,,\\
[\sigma_{\mu'\nu'\xi'\pi'}] &=
	-[\sigma_{\mu'\nu'\xi'\pi}] =
	-\frac{2}{3}\, R_{\pi(\mu|\xi|\nu)}
	\,.
\end{align*}

\subsection{Fifth order}

\vspace{-5ex}
\begin{align*}
[\sigma_{\mu'\nu'\xi'\pi'\alpha}] &=
	[\sigma_{\alpha\mu'\nu'\xi'\pi'}]
	\,,\displaybreak[0]\\
[\sigma_{\alpha\mu'(\nu'\xi'\pi')}] &=
	\frac{1}{6}\, (
	[\sigma_{\alpha\mu'\nu'\xi'\pi'}]
	+[\sigma_{\alpha\mu'\nu'\pi'\xi'}]
	+[\sigma_{\alpha\mu'\xi'\nu'\pi'}]
	\displaybreak[0]\\ &
	+[\sigma_{\alpha\mu'\xi'\pi'\nu'}]
	+[\sigma_{\alpha\mu'\pi'\nu'\xi'}]
	+[\sigma_{\alpha\mu'\pi'\xi'\nu'}]
	)
	\displaybreak[0]\\ &=
	\frac{1}{3}\, (
	[\sigma_{\alpha\mu'\nu'\xi'\pi'}]
	+[\sigma_{\alpha\mu'\xi'\nu'\pi'}]
	+[\sigma_{\alpha\mu'\pi'\nu'\xi'}]
	)
	\displaybreak[0]\\ &=
	\frac{2}{3}\, [\sigma_{\alpha\mu'\nu'\xi'\pi'}]
	-\frac{1}{3}\, [(\sigma_{\alpha a'} R^{a'}{}_{\mu'\nu'\xi'})_{;\pi'}]
	\displaybreak[0]\\ &
	+\frac{1}{3}\, [\sigma_{\alpha\mu'\nu'\pi'\xi'}]
	-\frac{1}{3}\, [(\sigma_{\alpha a'} R^{a'}{}_{\mu'\nu'\pi'})_{;\xi'}]	
	\displaybreak[0]\\ &=
	[\sigma_{\alpha\mu'\nu'\xi'\pi'}]
	+\frac{1}{3}\, R_{\alpha\mu\nu\xi;\pi}
	+\frac{1}{3}\, R_{\alpha\mu\nu\pi;\xi}	
	\,,\displaybreak[0]\\
[\sigma_{\alpha\mu'\nu'\xi'\pi'}] &=
	[\sigma_{\alpha\mu'(\nu'\xi'\pi')}]
	-\frac{2}{3} R_{\alpha\mu\nu(\xi;\pi)}
	=
	-\frac{1}{2}\, R_{\alpha(\nu|\mu|\xi;\pi)}
	-\frac{2}{3} R_{\alpha\mu\nu(\xi;\pi)}
	\displaybreak[0]\\ &=
	-\frac{1}{6}\, R_{\alpha\nu\mu(\xi;\pi)}
	-\frac{1}{6}\, R_{\alpha(\xi|\mu\nu|;\pi)}
	-\frac{1}{6}\, R_{\alpha(\xi|\mu|\pi);\nu}
	-\frac{2}{3} R_{\alpha\mu\nu(\xi;\pi)}
	\displaybreak[0]\\ &=
	\frac{5}{6}\, R_{(\mu|\alpha|\nu)(\xi;\pi)}
	-\frac{1}{6}\, R_{(\xi|\alpha|\pi)(\mu;\nu)}
	\,,\displaybreak[0]\\
[\sigma_{\mu'\nu'\xi'\alpha\beta}] &=
	[\sigma_{\mu'\nu'\xi'\alpha}]_{;\beta}
	-[\sigma_{\mu'\nu'\xi'\beta'\alpha}]
	\\ &
	=\frac{2}{3}\, R_{\alpha(\mu|\xi|\nu);\beta}
	+\frac{5}{6}\, R_{\alpha(\mu\nu)(\xi;\beta)}
	-\frac{1}{6}\, R_{\alpha(\xi\beta)(\mu;\nu)}
	\\ &
	=\frac{2}{3}\, R_{\alpha(\mu|\xi|\nu);\beta}
	+\frac{5}{12}\, R_{\alpha(\mu\nu)\xi;\beta}
	+\frac{5}{12}\, R_{\alpha(\mu\nu)\beta;\xi}
	\\ &
	-\frac{1}{12}\, R_{\alpha\xi\beta(\mu;\nu)}
	-\frac{1}{12}\, R_{\alpha\beta\xi(\mu;\nu)}
	\\ &
	=\frac{2}{3}\, R_{(\mu|\xi|\nu)(\alpha;\beta)}
	+\frac{5}{12}\, R_{\xi(\mu\nu)(\alpha;\beta)}
	+\frac{5}{12}\, R_{\alpha(\mu\nu)\beta;\xi}
	-\frac{1}{12}\, R_{(\alpha|\xi|\beta)(\mu;\nu)}
	\\ &
	=\frac{1}{4}\, R_{(\mu|\xi|\nu)(\alpha;\beta)}
	-\frac{5}{12}\, R_{(\mu|\alpha|\nu)\beta;\xi}
	-\frac{1}{12}\, R_{(\alpha|\xi|\beta)(\mu;\nu)}
	\,.
\end{align*}

\subsection{Sixth order}

\vspace{-5ex}
\begin{align*}
[\sigma_{\mu'\nu'(\alpha\beta\gamma\delta)}] &=
	-[\sigma_{\mu'(\alpha\beta\gamma\delta)\nu}] 
	\,,\displaybreak[0]\\
[\sigma_{\mu'\alpha(\beta\gamma\delta\nu)}] &= 
	\frac{1}{4}\, (
	[\sigma_{\mu'\alpha(\beta\gamma\delta)\nu}]
	+[\sigma_{\mu'\alpha(\beta\gamma|\nu|\delta)}]
	+[\sigma_{\mu'\alpha(\beta|\nu|\gamma\delta)}]
	+[\sigma_{\mu'\alpha\nu(\beta\gamma\delta)}]
	)
	\displaybreak[0]\\ &=
	\frac{1}{4}\, [\sigma_{\mu'\alpha(\beta\gamma\delta)\nu}]
	+\frac{1}{4}\, [\sigma_{\mu'\alpha(\beta\gamma|\nu|\delta)}]
	+\frac{1}{2}\, [\sigma_{\mu'\alpha(\beta|\nu|\gamma\delta)}]
	\displaybreak[0]\\ &
	-\frac{1}{4}\, [(\sigma_{\mu' a} R^a{}_{\alpha(\beta|\nu|})_{;\gamma\delta)}]
	\displaybreak[0]\\ &=
	\frac{1}{4}\, [\sigma_{\mu'\alpha(\beta\gamma\delta)\nu}]
	+\frac{3}{4}\, [\sigma_{\mu'\alpha(\beta\gamma|\nu|\delta)}]
	-\frac{1}{2}\, [(\sigma_{\mu'a(\beta} R^a{}_{|\alpha|\gamma|\nu|})_{;\delta)}]
	\displaybreak[0]\\ &
	-\frac{1}{2}\, [(\sigma_{\mu'\alpha a} R^a{}_{(\beta\gamma|\nu|})_{;\delta)}]
	-\frac{1}{4}\, [\sigma_{\mu' a (\beta\gamma}] R^a{}_{|\alpha|\delta)\nu}
	+\frac{1}{4}\, R_{\mu\alpha(\beta|\nu|;\gamma\delta)}
	\displaybreak[0]\\ &=
	[\sigma_{\mu'\alpha(\beta\gamma\delta)\nu}]
	-\frac{3}{4}\, [\sigma_{\mu' a (\beta\gamma}] R^a{}_{|\alpha|\delta)\nu}
	-\frac{3}{4}\, [\sigma_{\mu'\alpha a (\beta}] R^a{}_{\gamma\delta)\nu}
	\displaybreak[0]\\ &
	-\frac{3}{4}\, [\sigma_{\mu'\alpha(\beta|a|}] R^a{}_{\gamma\delta)\nu}
	-\frac{1}{2}\, [\sigma_{\mu'a(\beta\gamma}] R^a{}_{|\alpha|\delta)\nu}
	\displaybreak[0]\\ &
	-\frac{1}{2}\, [\sigma_{\mu'\alpha a (\beta}] R^a{}_{\gamma\delta)\nu}
	-\frac{1}{4}\, [\sigma_{\mu' a (\beta\gamma}] R^a{}_{|\alpha|\delta)\nu}
	+\frac{1}{4}\, R_{\mu\alpha(\beta|\nu|;\gamma\delta)}
	\displaybreak[0]\\ &=
	[\sigma_{\mu'\alpha(\beta\gamma\delta)\nu}]
	-\frac{3}{2}\, [\sigma_{\mu' a (\beta\gamma}] R^a{}_{|\alpha|\delta)\nu}
	-\frac{5}{4}\, [\sigma_{\mu'\alpha a (\beta}] R^a{}_{\gamma\delta)\nu}
	\displaybreak[0]\\ &
	-\frac{3}{4}\, [\sigma_{\mu'\alpha(\beta|a|}] R^a{}_{\gamma\delta)\nu}
	+\frac{1}{4}\, R_{\mu\alpha(\beta|\nu|;\gamma\delta)}
	\,,\displaybreak[0]\\
[\sigma_{\mu'\alpha(\beta\gamma\delta)\nu}] &=
	[\sigma_{\mu'\alpha(\beta\gamma\delta\nu)}]
	+\frac{3}{2}\,[\sigma_{\mu' a (\beta\gamma}] R^a{}_{|\alpha|\delta)\nu}
	+\frac{5}{4}\,[\sigma_{\mu'\alpha a (\beta}] R^a{}_{\gamma\delta)\nu}
	\displaybreak[0]\\ &
	+\frac{3}{4}\,[\sigma_{\mu'\alpha(\beta|a|}] R^a{}_{\gamma\delta)\nu}
	-\frac{1}{4}\,R_{\mu\alpha(\beta|\nu|;\gamma\delta)}
	\displaybreak[0]\\ &
	=
	-\frac{3}{5}\, R_{\mu(\beta|\alpha|\gamma;\delta\nu)}
	-\frac{7}{15}\, R_{\mu(\beta| a |\gamma} R^a{}_{\delta|\alpha|\nu)}
	-\frac{1}{2}\, R_{\mu(\beta| a |\gamma} R^a{}_{|\alpha|\delta)\nu}
	\displaybreak[0]\\ &
	-\frac{5}{24}\, R_{\mu a \alpha(\beta} R^a{}_{\gamma\delta)\nu}
	-\frac{5}{24}\, R_{\mu(\beta|\alpha a|} R^a{}_{\gamma\delta)\nu}
	-\frac{5}{8}\, R_{\mu\alpha a (\beta} R^a{}_{\gamma\delta)\nu}
	\displaybreak[0]\\ &
	-\frac{1}{8}\, R_{\mu (\beta|\alpha a |} R^a{}_{\gamma\delta)\nu}
	-\frac{1}{8}\, R_{\mu a \alpha(\beta} R^a{}_{\gamma\delta)\nu}
	-\frac{3}{8}\, R_{\mu\alpha(\beta| a |} R^a{}_{\gamma\delta)\nu}
	\displaybreak[0]\\ &
	-\frac{1}{4}\,R_{\mu\alpha(\beta|\nu|;\gamma\delta)}
	\displaybreak[0]\\ &
	=
	-\frac{3}{5}\, R_{\mu(\beta|\alpha|\gamma;\delta\nu)}
	-\frac{1}{4}\,R_{\mu\alpha(\beta|\nu|;\gamma\delta)}
	-\frac{7}{15}\, R_{\mu(\beta| a |\gamma} R^a{}_{\delta|\alpha|\nu)}
	\displaybreak[0]\\ &
	-\frac{1}{2}\, R_{\mu(\beta| a |\gamma} R^a{}_{|\alpha|\delta)\nu}
	-\frac{1}{3}\, R_{\mu a \alpha(\beta} R^a{}_{\gamma\delta)\nu}
	-\frac{1}{3}\, R_{\mu(\beta|\alpha a|} R^a{}_{\gamma\delta)\nu}
	\displaybreak[0]\\ &
	-\frac{5}{8}\, R_{\mu\alpha a (\beta} R^a{}_{\gamma\delta)\nu}
	-\frac{3}{8}\, R_{\mu\alpha(\beta| a |} R^a{}_{\gamma\delta)\nu}
	\,,\displaybreak[0]\\
[\sigma_{\mu'(\alpha\beta\gamma\delta)\nu}] &=
	-\frac{3}{20}\, R_{\mu(\alpha\beta|\nu|;\gamma\delta)}
	-\frac{1}{4}\, R_{\mu(\alpha\beta|\nu|;\gamma\delta)}
	-\frac{7}{60}\, R_{\mu(\alpha| a |\beta} R^a{}_{\gamma\delta)\nu}
	\displaybreak[0]\\ &
	-\frac{1}{2}\, R_{\mu(\alpha| a |\beta} R^a{}_{\gamma\delta)\nu}
	-\frac{1}{3}\, R_{\mu(\alpha\beta| a |} R^a{}_{\gamma\delta)\nu}
	-\frac{5}{8}\, R_{\mu(\alpha| a |\beta} R^a{}_{\gamma\delta)\nu}
	\displaybreak[0]\\ &
	-\frac{3}{8}\, R_{\mu(\alpha\beta| a |} R^a{}_{\gamma\delta)\nu}
	\displaybreak[0]\\ &
	=
	-\frac{2}{5}\, R_{\mu(\alpha\beta|\nu|;\gamma\delta)}
	-\frac{149}{120}\, R_{\mu(\alpha| a |\beta} R^a{}_{\gamma\delta)\nu}
	-\frac{17}{24}\, R_{\mu(\alpha\beta| a |} R^a{}_{\gamma\delta)\nu}
	\displaybreak[0]\\ &
	=
	\frac{2}{5}\, R_{\mu(\alpha|\nu|\beta;\gamma\delta)}
	+\frac{8}{15}\, R_{ a (\alpha|\mu|\beta} R^a{}_{\gamma|\nu|\delta)}
	\,,\displaybreak[0]\\
[\sigma_{\mu'\nu'\xi'(\alpha\beta\gamma)}] &=
	[\sigma_{\mu'\nu'(\alpha\beta\gamma)}]_{;\xi}
	-[\sigma_{\mu'\nu'(\alpha\beta\gamma)\xi}]
	\displaybreak[0]\\ &=
	[\sigma_{\mu'\nu'(\alpha\beta\gamma)}]_{;\xi}
	-[\sigma_{\mu'(\alpha\beta\gamma)\xi}]_{;\nu}
	+[\sigma_{\mu'(\alpha\beta\gamma)\xi\nu}]
	\,,\displaybreak[0]\\
[\sigma_{\mu'(\alpha\beta\gamma\xi)\nu}] &=
	\frac{1}{4}\, (
	[\sigma_{\mu'(\alpha\beta\gamma)\xi\nu}]
	+[\sigma_{\mu'(\alpha\beta|\xi|\gamma)\nu}]
	+[\sigma_{\mu'(\alpha|\xi|\beta\gamma)\nu}]
	+[\sigma_{\mu'\xi(\alpha\beta\gamma)\nu}]
	)
	\displaybreak[0]\\ &=
	\frac{1}{4}\, [\sigma_{\mu'(\alpha\beta\gamma)\xi\nu}]
	+\frac{1}{4}\, [\sigma_{\mu'(\alpha\beta|\xi|\gamma)\nu}]
	+\frac{1}{2}\, [\sigma_{\mu'(\alpha|\xi|\beta\gamma)\nu}]
	\displaybreak[0]\\ &=
	\frac{1}{4}\, [\sigma_{\mu'(\alpha\beta\gamma)\xi\nu}]
	+\frac{3}{4}\, [\sigma_{\mu'(\alpha\beta|\xi|\gamma)\nu}]
	-\frac{1}{2}\, [(\sigma_{\mu' a} R^a{}_{(\alpha\beta|\xi|})_{;\gamma)\nu}]
	\displaybreak[0]\\ &=
	[\sigma_{\mu'(\alpha\beta\gamma)\xi\nu}]
	-\frac{3}{2}\, [(\sigma_{\mu' a (\alpha} R^a{}_{\beta\gamma)\xi})_{;\nu}]
	\displaybreak[0]\\ &
	-\frac{1}{2}\, [\sigma_{\mu' a (\alpha|\nu|}] R^a{}_{\beta\gamma)|\xi}
	+\frac{1}{2}\, R_{\mu(\alpha\beta|\xi|;\gamma)\nu}
	\displaybreak[0]\\ &=
	[\sigma_{\mu'(\alpha\beta\gamma)\xi\nu}]
	-2\, [\sigma_{\mu' a (\alpha|\nu|}] R^a{}_{\beta\gamma)\xi}
	+\frac{1}{2}\, R_{\mu(\alpha\beta|\xi|;\gamma)\nu}
	\,,\displaybreak[0]\\
[\sigma_{\mu'(\alpha\beta\gamma)\xi\nu}] &=
	[\sigma_{\mu'(\alpha\beta\gamma\xi)\nu}]
	+2\, [\sigma_{\mu' a (\alpha|\nu|}] R^a{}_{\beta\gamma)\xi}
	-\frac{1}{2}\, R_{\mu(\alpha\beta|\xi|;\gamma)\nu}
	\displaybreak[0]\\ &=
	\frac{2}{5}\, R_{\mu(\alpha|\nu|\beta;\gamma\xi)}
	+\frac{8}{15}\, R_{ a (\alpha|\mu|\beta} R^a{}_{\gamma|\nu|\xi)}
	\displaybreak[0]\\ &
	+\frac{4}{3}\, R_{(\mu| a |\nu)(\alpha} R^a{}_{\beta\gamma)\xi}
	-\frac{1}{2}\, R_{\mu(\alpha\beta|\xi|;\gamma)\nu}
	\,,\displaybreak[0]\\
[\sigma_{\mu'\nu'\xi'(\alpha\beta\gamma)}] &=
	-\frac{5}{12}\, R_{(\alpha|\mu|\beta|\nu|;\gamma)\xi}
	-\frac{1}{12}\, R_{\mu(\alpha|\nu|\beta;\gamma)\xi}
	+\frac{1}{12}\, R_{(\alpha| \mu \xi|\beta;\gamma)\nu}
	\displaybreak[0]\\ &
	-\frac{5}{12}\, R_{(\alpha| \mu |\beta|\xi|;\gamma)\nu}
	+\frac{2}{5}\, R_{\mu(\alpha|\nu|\beta;\gamma\xi)}
	+\frac{8}{15}\, R_{ a (\alpha|\mu|\beta} R^a{}_{\gamma|\nu|\xi)}
	\displaybreak[0]\\ &
	+\frac{4}{3}\, R_{(\mu| a |\nu)(\alpha} R^a{}_{\beta\gamma)\xi}
	-\frac{1}{2}\, R_{\mu(\alpha\beta|\xi|;\gamma)\nu}
	\displaybreak[0]\\ &=
	-\frac{1}{2}\, R_{\mu(\alpha|\nu|\beta;\gamma)\xi}
	+\frac{2}{5}\, R_{\mu(\alpha|\nu|\beta;\gamma\xi)}
	\displaybreak[0]\\ &
	+\frac{8}{15}\, R_{ a (\alpha|\mu|\beta} R^a{}_{\gamma|\nu|\xi)}
	+\frac{4}{3}\, R_{(\mu| a |\nu)(\alpha} R^a{}_{\beta\gamma)\xi}
	\displaybreak[0]\\ &=
	-\frac{1}{2}\, R_{\mu(\alpha|\nu|\beta;\gamma)\xi}
	+\frac{1}{10}\, R_{\mu(\alpha|\nu|\beta;\gamma)\xi}
	+\frac{1}{10}\, R_{\mu(\alpha|\nu|\beta;|\xi|\gamma)}
	\displaybreak[0]\\ &
	+\frac{1}{10}\, R_{\mu(\alpha|\nu\xi|;\beta\gamma)}
	+\frac{1}{10}\, R_{\mu\xi\nu(\alpha;\beta\gamma)}
	+\frac{2}{15}\, R_{ a (\alpha|\mu|\beta} R^a{}_{\gamma)\nu\xi}
	\displaybreak[0]\\ &
	+\frac{2}{15}\, R_{ a (\alpha|\mu|\beta} R^a{}_{|\xi\nu|\gamma)}
	+\frac{2}{15}\, R_{ a (\alpha|\mu\xi} R^a{}_{\beta|\nu|\gamma)}
	+\frac{2}{15}\, R_{ a \xi\mu(\alpha} R^a{}_{\beta|\nu|\gamma)}
	\displaybreak[0]\\ &
	+\frac{4}{3}\, R_{(\mu| a |\nu)(\alpha} R^a{}_{\beta\gamma)\xi}
	\,.
\end{align*}
Symmetrization in $(\mu,\nu,\xi)$ is assumed in all terms of the following expression,
\begin{align*}
[\sigma_{(\mu'\nu'\xi')(\alpha\beta\gamma)}] &=
	-\frac{2}{5}\, R_{\mu(\alpha|\nu|\beta;\gamma)\xi}
	+\frac{1}{10}\, R_{\mu(\alpha|\nu|\beta;|\xi|\gamma)}
	-\frac{4}{3}\, R_{a \mu(\alpha|\nu|} R^a{}_{\beta|\xi|\gamma)}
	\displaybreak[0]\\ &=
	-\frac{2}{5}\, R_{\mu(\alpha|\nu|\beta;|\xi|\gamma)}
	+\frac{2}{5}\, R_{ a (\alpha|\nu|\beta} R^a{}_{|\mu\xi|\gamma)}
	+\frac{2}{5}\, R_{\mu a \nu (\alpha} R^a{}_{\beta|\xi|\gamma)}
	\displaybreak[0]\\ &
	+\frac{1}{10}\, R_{\mu(\alpha|\nu|\beta;|\xi|\gamma)}
	-\frac{4}{3}\, R_{a \mu(\alpha|\nu|} R^a{}_{\beta|\xi|\gamma)}
	\displaybreak[0]\\ &=
	-\frac{3}{10}\, R_{\mu(\alpha|\nu|\beta;|\xi|\gamma)}
	-\frac{4}{3}\, R_{a \mu(\alpha|\nu|} R^a{}_{\beta|\xi|\gamma)}
	\,.
\end{align*}

\begin{Bibliographyno}

\bibitem{Avramidi1991}
	Avramidi, I. G.,
    A Covariant Technique for the Calculation of the One-Loop Effective Action,
    \textit{Nuclear Phys. B} \textbf{355} (1991), no. 3, 712--754.

\bibitem{Avramidi1999}
	Avramidi, I. G.,
    Covariant Techniques for Computation of the Heat Kernel,
    \textit{Rev. Math. Phys.} \textbf{11} (1999), 947--980, {\tt arXiv:hep-th/9704166}.

\bibitem{Avramidi2000}
	Avramidi, I. G.,
    \textit{Heat Kernel and Quantum Gravity},
	Springer-Verlag, Berlin (2000).

\bibitem{AvramidiNotes}
	Avramidi, I. G.,
    \textit{Lecture Notes on Mathematical Physics}, (2000).
    Retrieved March~2, 2005, from\\ http://www.nmt.edu/$\sim$iavramid/notes/mathphys.pdf

\bibitem{Berline}
	Berline, N., Getzler E., Vergne M.,
	\textit{Heat Kernels and Dirac Operators},
	Springer-Verlag, Berlin (2004).

\bibitem{Donnelly}
	Donnelly, H.,
	Spectrum and the Fixed Point Sets of Isometries. I,
	\textit{Math. Ann.} \textbf{224} (1976), no. 2, 161--170.

\bibitem{Estrada}
	Estrada, R. and Fulling, S. A.,
    Distributional Asymptotic Expansions of Spectral Functions and of the Associated Green Kernels,
	\textit{Electron. J. Differential Equations}, \textbf{7} (1999), {\tt arXiv:funct-an/9710003}.

\bibitem{Erdelyi}
	Erd\'elyi, A., Magnus, W., Oberhettinger, F., and Tricomi, F. G.,
    \textit{Higher Transcendental Functions. Vols. I, II},
    McGraw-Hill Book Company, Inc., New York (1953).
    
\bibitem{Fedoryuk}
	Fedoryuk, M. V.,
	\textit{Asymptotics: Integrals and Series},
	``Nauka'', Moscow (1987) [In Russian].
	
\bibitem{Frankel}
	Frankel, T.,
	\textit{The Geometry of Physics},
	Cambridge University Press, Cambridge (2004).	
	
\bibitem{Gilkey}
	Gilkey, P. B.,
	\textit{Invariance Theory, the Heat Equation, and the Atiyah-Singer Index Theorem},
	Chemical Rubber Company, Boca Raton, Florida (1995).

\bibitem{Grigoryan}
	Grigor'yan, A.,
	Estimates of Heat Kernels on Riemannian Manifolds,
	\textit{Spectral theory and geometry},
	Cambridge Univ. Press, Cambridge (1999), 140--225.

\bibitem{Kirsten}
	Kirsten, K.,
	\textit{Spectral Functions in Mathematics and Physics},
	Chapman \& Hall/CRC, Boca Raton, Florida (2001).
	
\bibitem{Kobayashi}
	Kobayashi, S. and Nomizu, K.,
	\textit{Foundations of Differential Geometry. Vol~I},
	Interscience Publishers, a division of John Wiley \& Sons, New York-London (1963).

\bibitem{Kreyszig}
	Kreyszig, E.,
	\textit{Introductory Functional Analysis with Applications},
	John Wiley \& Sons, New York (1978).
	
\bibitem{Oersted}
	Oersted, B.,
	\textit{Equivariant zeta functions and conformal geometry},
	MSRI video archive.
	Retrieved November 12, 2004 from\\
	http://elasmo.kaist.ac.kr/ln/msri/2001/spectral/oersted/1/

\bibitem{Rosenberg}
	Rosenberg, S.,
	\textit{The Laplacian on a Riemannian Manifold},
	Cambridge University Press, Cambridge (1997). 
	
\bibitem{Synge}
	Synge, J. L.,
	\textit{Relativity: The General Theory},
    North-Holland Publishing Co., Amsterdam (1960).
    
\bibitem{Vassilevich}
	Vassilevich, D. V.,
	Heat Kernel Expansion: User's Manual,
	\textit{Phys. Rep.}, \textbf{388}~(2003), no. 5--6, 279--360, {\tt arXiv:hep-th/0306138}.

\end{Bibliographyno}

\end{document}